\documentclass[final]{siamltex}
\usepackage{showlabels}
\usepackage{xcolor}
\usepackage{amsmath}
\usepackage{amssymb}
\usepackage{mathrsfs}
\usepackage{graphicx}
 \usepackage{bm}
\usepackage{color}
\usepackage{float}
\usepackage{epstopdf}
\usepackage{caption}
\usepackage{subcaption}
\usepackage[breaklinks,colorlinks=true,linkcolor=blue,citecolor=red,backref=page]{hyperref}

\DeclareMathAlphabet{\itbf}{OML}{cmm}{b}{it}
\DeclareMathAlphabet\mathbfcal{OMS}{cmsy}{b}{n}

\renewcommand{\hat}{\widehat}
\renewcommand{\tilde}{\widetilde}
\def\RR{\mathbb{R}}
\def\bx{{{\itbf x}}}

\def\by{{{\itbf y}}}

\def\bu{{{\itbf u}}}

\def\bbeta{\boldsymbol{\eta}}
\def\bg{{\itbf g}}
\def\be{{\itbf e}}

\def\bU{{\itbf U}}
\def\bY{{\itbf Y}}

\def\bv{{\itbf v}}

\def\bV{{\itbf V}}
\def\bR{{\itbf R}}

\def\bI{{\itbf I}}
\def\bS{{\itbf S}}

\def\bM{{\itbf M}}
\def\bD{{\itbf D}}

\def\cP{\mathcal{P}}

\def\cA{{\mathcal A}}

\def\bLa{{\boldsymbol{\Lambda}}}

\def\FWI{{\scalebox{0.5}[0.4]{FWI}}}

\def\RM{{\scalebox{0.5}[0.4]{ROM}}}
\def\bAR{\cA^{{\scalebox{0.5}[0.4]{ROM}}}}

\def\om{\omega}
\def\la{\lambda}
\def\12{{\frac{1}{2}}}

\newtheorem{rem}[theorem]{Remark}
\newtheorem{algorithm}{Algorithm}
\newtheorem{lem}{Lemma}
\newtheorem{thm}{theorem}
\newtheorem{cor}{Corollary}

\newcommand{\bc}{\textcolor{black}}

\begin{document}

\title{When data driven reduced order modeling meets full waveform inversion} \author{Liliana Borcea\footnotemark[1] \and Josselin
  Garnier\footnotemark[2] \and Alexander V. Mamonov\footnotemark[3] \and J\"{o}rn Zimmerling\footnotemark[4]}

\maketitle


\renewcommand{\thefootnote}{\fnsymbol{footnote}}
\footnotetext[1]{Department of Mathematics, University of Michigan,
  Ann Arbor, MI 48109. {\tt borcea@umich.edu}}
\footnotetext[2]{CMAP, CNRS, Ecole Polytechnique, Institut Polytechnique de Paris, 91120 Palaiseau, France.  {\tt
    josselin.garnier@polytechnique.edu}}
\footnotetext[3]{Department of Mathematics, University of Houston, TX 77204-3008. {\tt mamonov@math.uh.edu}}
\footnotetext[4]{ Uppsala Universitet, Department of Information Technology, Division of Scientific Computing, 75105 Uppsala, Sweden. {\tt jorn.zimmerling@it.uu.se}}
\markboth{L. BORCEA, J. GARNIER, A.V. Mamonov, J. Zimmerling}{Waveform inversion}

\begin{abstract}
Waveform inversion is concerned with estimating a heterogeneous medium, modeled by variable coefficients of wave equations, 
using sources that emit probing signals and receivers that record the generated waves. It is an old and intensively studied inverse problem with a  wide range of applications, but the existing inversion methodologies are still far from satisfactory. The typical mathematical formulation is a nonlinear least squares data fit optimization  and the difficulty stems  from the non-convexity of the objective function that displays numerous local minima at which local optimization approaches stagnate.  This pathological behavior has at least  three unavoidable causes: (1) The mapping from the unknown coefficients to the wave field is nonlinear and complicated. (2) 
The sources and receivers typically lie on a single side of the medium, so only  backscattered waves are measured. (3) The probing signals are band limited and with high frequency content.  There is a lot of activity in the computational science and engineering communities that seeks to mitigate 
the difficulty of estimating the medium by data fitting. In this paper we present a different point  of view, based on  reduced order models (ROMs) of two operators that control the wave propagation. The ROMs are called data driven because they are computed directly from the measurements, without any knowledge of the wave field inside the inaccessible medium. This computation is non-iterative and uses  
standard numerical linear algebra methods. The resulting ROMs capture features of the physics of wave propagation  in a complementary way and have surprisingly  good approximation properties that facilitate  waveform inversion. \\

In this arxiv version two important typos are corrected when compared to the published \cite{PublishedSirev} version. The typo was in Equation~\eqref{eq:DF2} in Theorem 3 and carried over into Corollary~\ref{lem.3}.
The proofs are correct.
\end{abstract}
\begin{keywords}
Inverse wave scattering, data driven, reduced order modeling, optimization.
\end{keywords}

\begin{AMS}
65M32, 41A20
\end{AMS}
\section{Introduction to waveform inversion. Paper motivation and outline}
The estimation of a heterogeneous medium from the time history of the wave field recorded at a few accessible  locations
is important in  medical diagnostics via ultrasound, nondestructive 
evaluation of aging concrete in bridges and buildings, testing of aircraft fuselage, radar imaging, underwater acoustics, seismology, geophysical 
exploration and so on. It is an inverse problem for a wave equation or a system of such equations (acoustic, electromagnetic or elastic),
where the heterogeneous medium is modeled by unknown variable coefficients. The goal is to determine these coefficients from measurements of the wave field, the solution of the wave equation, with forcing that is typically, but not always, controlled by the user.

To keep the presentation simple, we consider acoustic waves in a medium with constant mass density and unknown variable wave speed $c(\bx)$, 
a piecewise smooth, non-negative and bounded function. The wave field is modeled by the acoustic pressure $p(t,\bx)$, the solution of the wave equation
\begin{equation}
\left[\partial_t^2 - c^2(\bx) \Delta \right] p(t,\bx) = S(t,\bx), \qquad t \in \RR, \quad \bx \in \Omega \subset \RR^d,
\label{eq:I1}
\end{equation}
in dimension $d = 2$ or $3$, with forcing $S(t,\bx)$ and homogeneous initial condition 
$
p(t,\bx) = 0$ at $ t \ll 0$ i.e., at negative time, outside the temporal support of $S(t,\bx)$. 
We assume  a bounded and simply connected domain $\Omega$ with smooth enough perfectly reflecting boundary $\partial \Omega$. This may be a real boundary of a closed cavity, but in general it is a fictitious boundary introduced for convenience of the analysis and computations. Note that the waves propagate at finite speed, so if we place $\partial \Omega$ far enough from the spatial support of $S(t,\bx)$, there are no boundary effects over the finite duration of the data gather experiment and  we can choose any homogeneous boundary conditions.  

The forcing in \eqref{eq:I1} is  commonly controlled by the user and is modeled as a point source at $\bx_s \in \Omega$ that emits a probing signal $f(t)$,
\begin{equation}
S(t,\bx) = f(t) \delta_{\bx_s}(\bx),
\label{eq:I5}
\end{equation}
where $\delta_{\bx_s}(\bx)$ is the Dirac delta  at $\bx_s$. Typically, there are $m_s$ such sources, which probe the medium one at a time, so to keep track of the generated wave, the solution of \eqref{eq:I1} with right-hand side \eqref{eq:I5}, we denote it by 
$p^{(s)}(t,\bx)$. The inverse problem is:  \emph{Determine $c(\bx)$ from the measurements at receiver locations $\by_r$, for $r = 1, \ldots, m_r$, organized in the $m_r \times m_s$ time-dependent matrix $ \boldsymbol{\mathcal{M}}(t)$ with entries}
\begin{equation}
\mathcal{M}_{r,s}(t) = p^{(s)}(t,\by_r), \qquad s = 1, \ldots, m_s, ~~ r = 1, \ldots, m_r, ~~ t \in (t_{\rm min},t_{\rm max}).
\label{eq:arrayResp}
\end{equation}

The duration and frequency content of the probing signal matter and are application specific. The  generic model  of the emitted signal is 
$ 
f(t) =  B\varphi(Bt) \cos (\om_o t),
$
where $B \varphi(Bt)$ is an envelope function supported at  $t \in [-1/B,1/B]$, with Fourier transform $\hat \varphi(\om/B)$ that is large at angular frequency $\om$ satisfying $|\om| \le O(B)$.  Thus,  $B$ is called the bandwidth. The modulation by the cosine shifts the frequency content of $f(t)$ to the intervals $|\om \pm \om_o| \le O(B)$, according to the formula 
\begin{equation}
\hat f(\om) = \int_{\RR} dt \, f(t) e^{i \om t} = \frac{1}{2} \left[ \hat \varphi \left(\frac{\om-\om_o}{B}\right) + \hat \varphi \left(\frac{\om+\om_o}{B}\right) \right],
\label{eq:deffh}
\end{equation}
so $\om_o$ is called the central frequency. Basic resolution studies suggest that signals with short $O(1/B)$ support and with high frequency content  produce sharper estimates of the medium \cite{cheney2009fundamentals,bleistein2001multidimensional}. This explains  the 
wide use of pulses with modulation  frequency $\om_o \gg B$,  but there are exceptions.  In long-range imaging with radar \cite{curlander1991synthetic}  pulses are not used because sources have limited instantaneous power of emission and, due to geometrical spreading,  the waves that  reach the receivers are weaker than  ambient noise. Chirped signals (long signals whose frequency increases with time) are used instead,
 because they allow an  increase in the delivered net power and better signal-to-noise ratios. The inversion methodology with pulses or chirps is the same, due to the 
simple ``pulse compression"  data processing  $f(-t) \star_t p^{(s)}(t,\by_r)$, where $f(-t)$ is the time-reversed emitted signal and $\star_t$ denotes  time convolution \cite{curlander1991synthetic}.  Mathematically, by linearity of the wave equation, this processing is equivalent to working with the wave generated by the signal 
\begin{equation}
F(t) = f(-t) \star_t f(t),
\label{eq:defF}
\end{equation}
which turns out to have much shorter duration than the chirp i.e.,  it is a pulse \cite{gilman2017transionospheric}. 

We refer the interested reader to the mathematical literature on the uniqueness and stability of the inverse problem   (see \cite{belishev2007recent,stefanov2005stable,yu2001global} and others).  Many of these studies are for different types of data, like the Dirichlet to Neumann map, which are not available in the applications that we have in mind. Furthermore, uniqueness 
cannot be expected  to hold in the strict sense, even in the best circumstances, because the measurements are band limited: Variations of  $c(\bx)$ on  scales that are much smaller than the central wavelength cannot be determined. This issue is  addressed in practice by a proper choice of the space $\mathcal{W}$ in which the search speed $w(\bx)$ lies and, by adding a regularization penalty to the optimization.

The generic  formulation of the inverse problem is  a nonlinear least squares data fit optimization: $  \min_{w \in \mathcal{W}} \mathcal{O}^\FWI(w) + \mbox{ regularization}$, with
\begin{equation}
\mathcal{O}^\FWI(w) = 
 \sum_{r = 1}^{m_r} \sum_{s = 1}^{m_s} \int_{t_{\rm min}}^{t_{\rm max}} dt \,  |\mathcal{M}_{r,s}(t) - p^{(s)}(t,\by_r;w)|^2.
\label{eq:FWI}
 \end{equation}
Here $ \boldsymbol{\mathcal{M}}(t)$ is the data matrix (\ref{eq:arrayResp}), $w(\bx) \mapsto p^{(s)}(t,\by_r;w)$ is called the forward map, and $p^{(s)}(t,\bx;w)$ denotes  the (numerical) solution of the  wave equation  with  wave speed  $w(\bx)$. \bc{Note that to avoid complicated notation, we omit the true wave speed from the list of arguments of the wave field and measurements. Thus, in our notation $p^{(s)}(t,\bx)$ is the pressure field for the true speed $c(\bx)$ and $p^{(s)}(t,\bx;w)$ is the pressure field for the search speed $w(\bx)$.}

\begin{figure*}[t]
\begin{center}
\begin{tabular}{ccc}
Wave speed $c$ in (m/s) &
Log of  $\mathcal{O}^\FWI$\\
\includegraphics[width=0.29\textwidth]{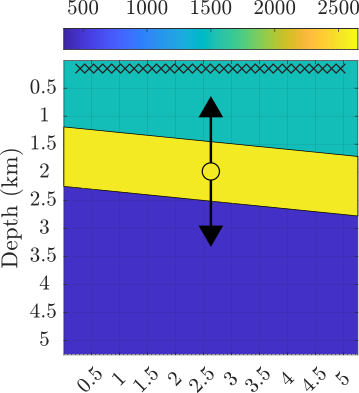} &
\raisebox{0.in}{\includegraphics[width=0.39\textwidth]{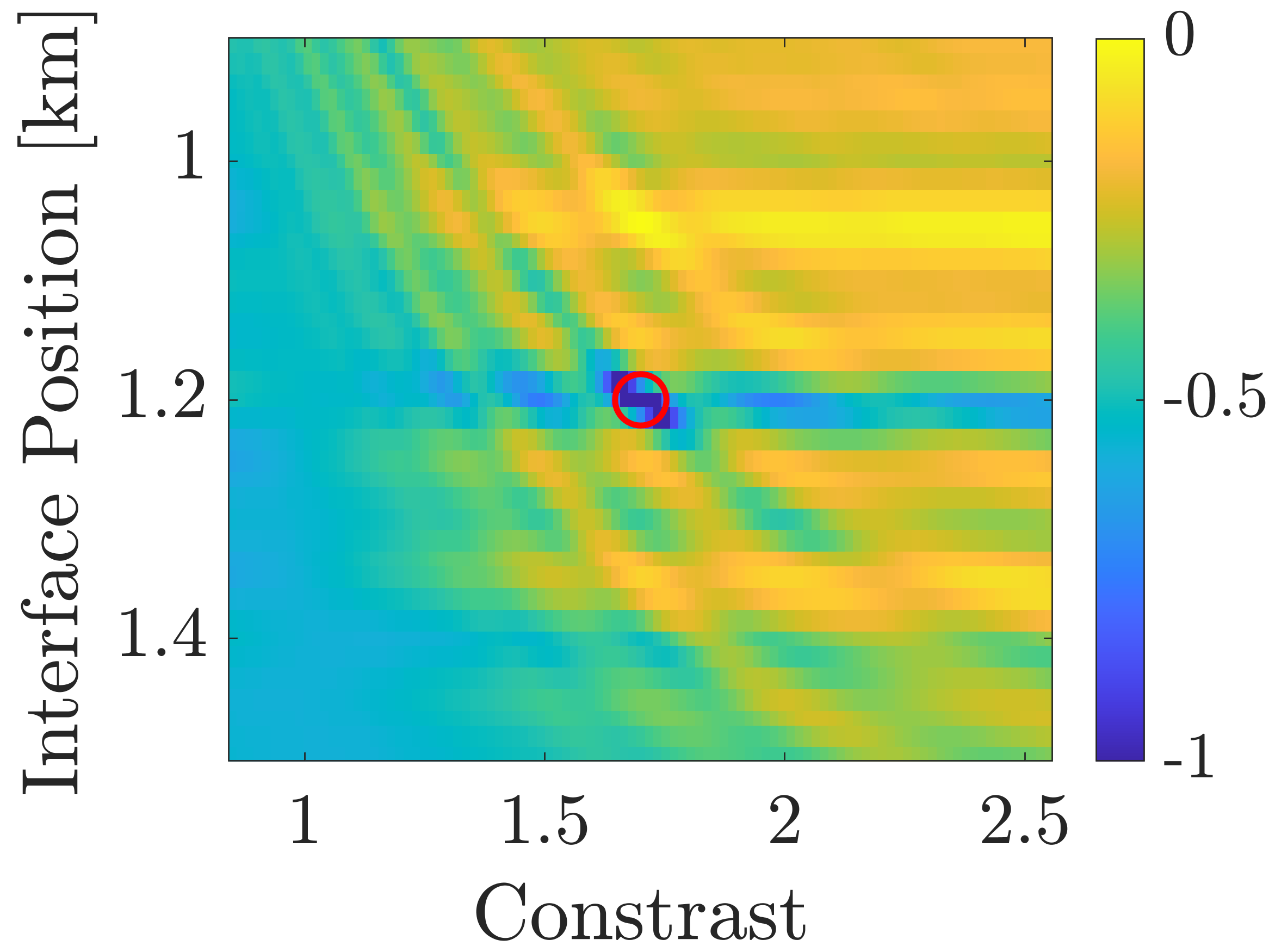}} \end{tabular}
\end{center}
\vspace{-0.1in}
\caption{Illustration of  a two-parameter search  applied to the true piecewise constant model $c(\bx)$
(left plot), in a closed rectangular cavity. The sources and receivers are co-located and shown as \bc{black $\times$}. The probing signal contains frequencies in the interval $[2,10]$Hz. The search parameters are the depth of the slanted fast layer (interface position), 
that varies over the range indicated by the \bc{black arrows}, and the ratio of the wave speed inside and above the layer (contrast). 
The FWI objective function (\ref{eq:FWI}) is shown in the right plot. {The true parameters are 
indicated  by $\textcolor{magenta}{\bigcirc}$.} }\label{fig:FWI_topo}
\end{figure*}

The optimization formulation defined by \eqref{eq:FWI} has been coined by the acronym FWI (Full Waveform Inversion) in the geophysics community \cite{virieux2009overview}, so we shall refer to it as such. 
The FWI  formulation applies to an arbitrary data acquisition geometry but the result of the inversion is strongly influenced by it. For example, it is easier to work with sources and receivers placed all around $\Omega$, so that both the waves transmitted through the medium and those that are backscattered can be measured.  Unfortunately, the options for placing the sources and receivers are usually limited to one side of $\Omega$, so  only the backscattered waves can be measured.  It is well known that the FWI objective function (\ref{eq:FWI}) for such measurements exhibits numerous local minima at search speeds $w(\bx)$ that bear no relation to $c(\bx)$ \cite{virieux2009overview}. See Fig. \ref{fig:FWI_topo} for an illustration. Any iterative local optimization method \cite{gill2019practical} applied to such an objective function will likely fail to obtain a good approximation of $c(\bx)$.

 The pathological behavior of the FWI objective function seen in Fig. \ref{fig:FWI_topo} is called ``cycle-skipping" and arises when the timing between the measured waves and the ones predicted by the forward model  differs by more than half the cycle of oscillation.  The main causes  are the complicated  nonlinear forward mapping  and the lack of low frequencies in the  probing signal. Waves are backscattered by the rough part of $c(\bx)$, the ``reflectivity", which corresponds to the jump discontinuities at the slanted layer 
in Fig.  \ref{fig:FWI_topo}. Reflections at these discontinuities and at the  top boundary, located just above the sources/receivers, cause multiple arrivals (echoes) in the recorded data. When moving the layer up and down, the errors in the arrival times of these echoes exceed half a cycle, repeatedly, so we get the multiple horizontal stripes (highs and lows)  in Fig.  \ref{fig:FWI_topo}. The smooth part of $c(\bx)$,  the  ``kinematics", determines the travel time of the transmitted waves and is another cause of  cycle-skipping in Fig. \ref{fig:FWI_topo}. The curved stripes in the figure appear because when misestimating $c(\bx)$ inside the layer, we get wrong arrival times of the waves that penetrate there, scatter at the bottom of the layer and then reach the receivers.

Cycle-skipping can sometimes be mitigated by adding prior information on $c(\bx)$ or by starting with a good initial guess. For the latter, the following observation is useful: It  is easier to  determine the kinematics from very low frequency data, because we can have larger errors of  travel times that remain within the long cycle of oscillation \cite{bunks1995multiscale,virieux2009overview}.
High frequencies are better for estimating the reflectivity, if the kinematics is known and the multiple scattering effects are not too strong \cite{symes2008migration,bleistein2001multidimensional}. This has motivated studies like \cite{chen1997inverse,borges2017high,bunks1995multiscale} that obtain progressively better 
estimates of $c(\bx)$, starting from the lower frequencies.  Unfortunately, low enough frequencies are rarely available and a reasonable initial guess does not guarantee success. Indeed, Fig. \ref{fig:FWI_topo} shows that  cycle-skipping may occur even near the true $c(\bx)$.  

The outstanding question is: How can we improve the objective function \eqref{eq:FWI} so that numerically feasible local optimization methods can give good estimates of $c(\bx)$, irrespective of the starting point? One approach is to get rid of the  $L^2([t_{\rm min},t_{\rm max}])$  metric and replace it with a better one, such as the Wasserstein metric 
from optimal transport theory \cite{EngquistFroese,yang2018application}. The convexity of the data misfit function in this metric has been shown for a few simple 
models in \cite{engquist2022optimal,mahankali2020convexity} and  it holds for  the example in Fig. \ref{fig:FWI_topo}. However, we will present later in this paper another example where convexity does not hold. 
Another approach, known as  ``modeling operator extension" \cite{huang2018source} or ``extended FWI" \cite{herrmann2013,warner2016},  introduces in a systematic way additional degrees of freedom in the optimization, and then drives the results towards a physically meaningful result. An explicit analysis of such a method, that minimizes the data misfit over both $c(\bx)$ and $f(t)$ in a very simple setting, can be found in \cite{symes2022error}. 
There are various other ideas that have been tried, including machine learning  \cite{corte2020deep,ding2022coupling}. Nevertheless, the state of the field remains far from satisfactory and theoretical guarantees that one method or another will work even in simple settings are rare.

Our goal in this paper is to show that tools from reduced order modeling can be used to improve the waveform inversion methodology.  Model order reduction is a popular topic in computational science, concerned mostly with reducing the computational complexity of a given dynamical system (the model), for purposes like design and control \cite{antoulas2001survey,schilders2008model,benner2015survey}.  The idea is to 
obtain a low-dimensional, computationally inexpensive reduced order model (ROM), that approximates the response of the true model over a 
range of input parameters.  \bc{In waveform inversion the  model is defined by the wave equation, which we know, and the 
wave speed $c(\bx)$, which is unknown.  Thus, we are interested in data driven ROMs  that are computed directly from the measurements of the wave field, without knowing the wave speed model. }

 Data driven reduced order modeling 
is a rapidly growing field that combines ideas from optimization, numerical analysis and projection based reduced order modeling. Much of it  
is concerned with using a sparse set of  ``snapshots" of solutions of a time-dependent partial differential equation  to either learn the equation \cite{brunton2019data,brunton2016discovering}  or to approximate the time evolution of its solution \cite{hesthaven2022reduced,herkt2013convergence,kunisch2010optimal,lieu2006reduced}. Both goals are of interest in waveform inversion,  but none of  these methods can be used because they assume knowledge of the snapshots at all points in $\Omega$.

\begin{figure*}[t]
\begin{center}
\includegraphics[width=0.28\textwidth]{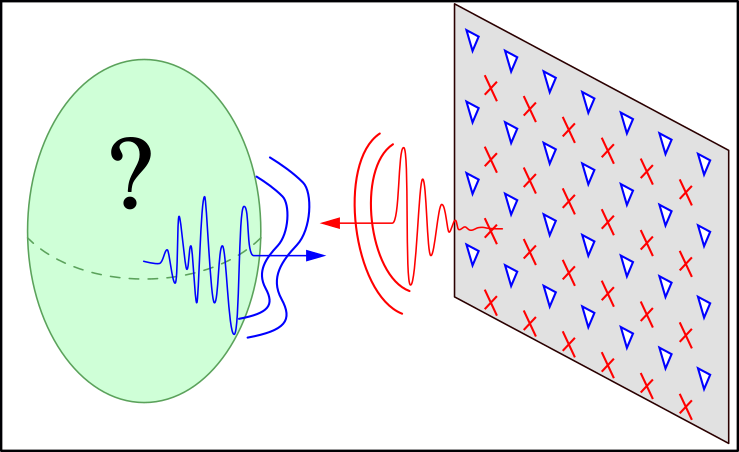}
\includegraphics[width=0.28\textwidth]{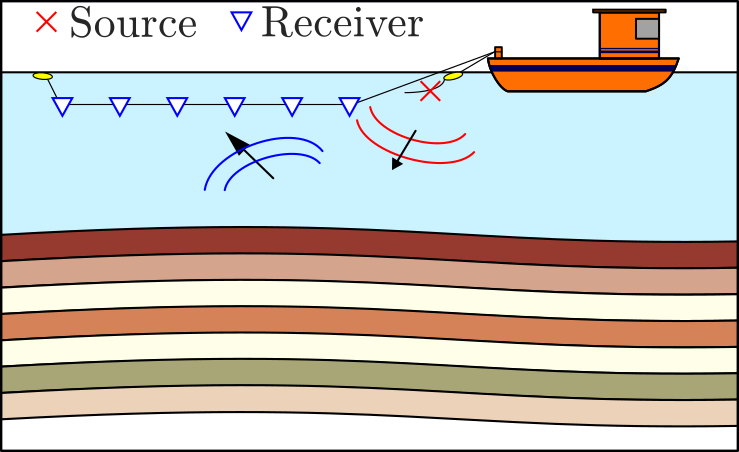}
\includegraphics[width=0.28\textwidth]{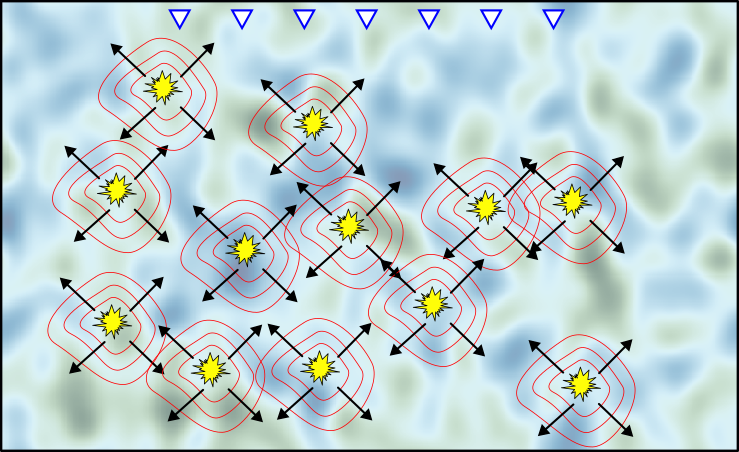}
\end{center}
\vspace{-0.in}
\caption{Illustration of  three data acquisition setups: With an active array (left), where the sources and receivers are co-located; with a towed-streamer (middle); \bc{with a passive array of receivers (blue triangles) and uncontrolled random sources (yellow) dispersed throughout the medium (right).} }
\label{fig:drawing}
\end{figure*}

The data driven ROMs used in this paper 
look like standard (reduced basis) projection ROMs \cite{hesthaven2022reduced,brunton2019data}, because they 
are  Galerkin projections of two  operators related to $-c^2(\bx) \Delta$, on the space 
spanned by the snapshots 
$
\{p^{(s)}(t_j,\bx), \bx \in \Omega\}_{s=1}^{m_s}$ at some time instants $t_j$, for $j = 0 , \ldots, n-1$. However, there are fundamental differences: 
\begin{enumerate}
\item The snapshots and therefore the approximation space are unknown. The ROMs can, however, be computed directly from the measurements  $\boldsymbol{\mathcal{M}}(t)$.
\item The ROMs are matrices with structure designed to capture, at the algebraic level,  the causal physics of wave propagation. They also 
have superior  approximation properties that prove useful in inversion.
\end{enumerate}
The starting point of our ROM construction is the mapping of the measurements $\boldsymbol{\mathcal{M}}(t)$ defined in \eqref{eq:arrayResp} to a new data matrix $\bD(t)$, whose components evaluated on the time grid $\{t_j\}_{j \ge 0}$ can be expressed as inner products of  snapshots.  The entries of the ``mass" and ``stiffness" matrices in the Galerkin approximations are determined by the same inner products, so they, and the ROMs,  can be computed from $\boldsymbol{\mathcal{M}}(t)$. 
The mapping  $\boldsymbol{\mathcal{M}}(t) \mapsto \bD(t)$ can be carried out without any knowledge of $c(\bx)$, but  it does require 
having  co-located sources and receivers: 
$\by_s = \bx_s, $ for $s = 1, \ldots, m_s = m_r = m.$ We assume henceforth such a setup (see the left plot in Fig. \ref{fig:drawing} for an illustration) and refer to the collection of  sources/receivers as an ``active array", which measures  the ``array response matrix" $\boldsymbol{\mathcal{M}}(t)$  
defined in \eqref{eq:arrayResp}.  This is an $m \times m$ symmetric matrix, due to source-receiver reciprocity. 

Note that the active array setup  is common in radar imaging and phased array ultrasonics, 
but it is not feasible in other fields, like geophysics.  Nevertheless, there are 
data acquisition setups  used in such fields, where it is possible to obtain a good approximation of the matrix $\boldsymbol{\mathcal{M}}(t)$ used by our construction. An example is the synthetic aperture,  
towed hydrophone streamer data acquisition  used in marine seismic surveys (middle plot in Fig. \ref{fig:drawing}), where $\boldsymbol{\mathcal{M}}(t)$
is approximated using source-receiver reciprocity and interpolation, as explained in  \cite{borcea2022waveform}. Another example 
discussed later in the paper considers  opportunistic data acquisition setups with ``passive arrays" of receivers, where the wave excitation comes from uncontrolled noise sources \cite{garnier2016passive} (see the illustration in the right plot of Fig. \ref{fig:drawing}).

The presentation of our ROM-based waveform inversion methodology is organized in  three sections: 

Section~\ref{sect:ROMs} gives the data driven computation of two ROMs: 
The first ROM is for the ``wave propagator operator", named so  
because it propagates the wave field from the instants $t_{j-1}$ and $ t_j$ to the next instant $ t_{j+1}$, 
for $j \ge 1$. It was first introduced in \cite{druskin2016direct}, where the resulting ROM was interpreted as a three-point spatial finite difference scheme on a special grid called ``optimal". Such grids have been used in the past for obtaining superconvergent approximations of Neumann to Dirichlet maps \cite{druskin2002three,asvadurov2003optimal} and for 
solving various inverse problems  \cite{borcea2005continuum,borcea2010circular,borcea2008electrical}. However, the optimal grids 
are understood only in one dimension, which is why the results in \cite{druskin2016direct} are limited to estimating the reflectivity of layered media.
The extension of the propagator ROM to higher dimensions and its analysis that establishes the connection to Galerkin projections is given in 
\cite{borcea2020reduced}, although the ROM was used before in \cite{DtB,borcea2019robust,druskin2018nonlinear} for the purpose of imaging with waves i.e., localizing the support of the reflectivity of a medium with known kinematics. The second ROM is for the operator $-c^2(\bx) \Delta$
and was introduced recently in \cite{borcea2022waveform}. We explain how the two ROMs are related, discuss their approximation properties 
and compare them with the time domain reduced basis ROMs in the computational literature.

Section~\ref{sect:WFI} is concerned with the use of the ROMs for waveform inversion. Roughly, this  involves the chain of maps 
$ \boldsymbol{\mathcal{M}}(t) \mapsto \bD(t)\ \mapsto \mbox{ROM} \mapsto \mbox{estimate of } c(\bx),
$
where the first two maps are the ones discussed in the previous section and the last map is computed using iterative optimization. \bc{We describe two distinct ideas for estimating $c(\bx)$: The first one is to use a  ROM estimate of the wave field inside $\Omega$, called the ``estimated internal wave".  This wave field depends on both the search and true wave speeds, it can be computed from the measurements, and it can be used to linearize approximately the map $c(\bx) \mapsto   \bD(t)$ \cite{borcea2022internal}. Here we introduce a better use of the estimated internal wave, based on  the observation that  this  wave fits the data, by construction, but unlike the true internal wave, it does not solve the wave equation.  Thus, we estimate $c(\bx)$ by minimizing the residual in the wave equation. The second idea 
was introduced  recently in \cite{borcea2022waveform}. It also gets rid of data fitting  and minimizes the ROM  misfit instead. }

Section~\ref{sect:passive} gives the generalization of the inversion methodology to the opportunistic data acquisition setup with random sources.

We end with a few concluding remarks and open questions in section \ref{sec:conc}.  

\section{The data driven ROMs}
\label{sect:ROMs}
The description of the computation of the ROMs involves four main steps: The first step, described in section \ref{sect:ROM.1}, is the data mapping $\boldsymbol{\mathcal{M}}(t) \mapsto \bD(t)$.  The second step, described in section \ref{sect:ROM.ts}, defines  the wave snapshots and the 
propagator operator that governs their evolution. 
The third step, given in section \ref{sect:ROM.2}, introduces two Galerkin approximation schemes: one for the time stepping of the snapshots and one for the wave equation. The approximations are in an unknown space, but  the Galerkin coefficients satisfy equations that  are data driven. The fourth step, given in section \ref{sect:ROM.4},  uses these equations to obtain our two ROMs. We give in section \ref{sect:ROM.5} a brief comparison  of these ROMs to the standard (reduced basis) projection ones found in the computational science literature.  \bc{We end in section \ref{sect:regul} with a brief discussion of regularization of the ROM computation.}

\subsection{Data transformation}
\label{sect:ROM.1}
We will use functional calculus for  a symmetrized wave operator. To define it, we must specify  the boundary conditions. Any homogeneous boundary conditions that ensure a trivial kernel of the operator will do, so to fix ideas we consider henceforth the Dirichlet boundary conditions
$
p^{(s)}(t,\bx) = 0,$  at $\bx \in \partial \Omega.$ We assume 
that $c(\bx)$ equals the known constant $\bar{c}$ near the array \bc{i.e.,  within the distance traveled by the waves during the short 
duration of the probing pulse $f(t)$.}
This assumption holds in most inversion setups and it simplifies the presentation. It is possible to extend the  results  to cases where it does not hold, but there are no significant insights brought by such an extension.

The operator $-c^2(\bx) \Delta$ with homogeneous Dirichlet boundary conditions 
is  positive definite and self-adjoint with respect to the inner product weighted by $c^{-2}(\bx)$.  We prefer to work with the $L^2(\Omega)$ inner product,  so we use the similarity transformation 
\begin{equation}
p^{(s)}(t,\bx) \mapsto P^{(s)}(t,\bx) = \frac{\overline{c}}{c(\bx)} p^{(s)}(t,\bx), \quad \bc{t \in \RR, ~~ \bx \in \Omega},
\label{eq:Sim1}
\end{equation}
that acts as the identity at the points in the array
\begin{equation}
P^{(s)}(t,\bx_r) = p^{(s)}(t,\bx_r) = \mathcal{M}_{r,s}(t), \qquad s,r = 1, \ldots, m.
\label{eq:Sim2}
\end{equation}
The wave operator $\partial_t^2 - c^2(\bx) \Delta$ is transformed by \eqref{eq:Sim1} to $\partial_t^2 + \cA$, 
where 
\begin{equation}
\cA = - c(\bx) \Delta \left[ c(\bx) \cdot \right],
\label{eq:D4}
\end{equation}
is positive definite and self-adjoint, with compact resolvent. Its spectrum  \cite[Section 5.3]{kato} consists of a sequence of positive 
eigenvalues $\{\theta_j\}_{j \ge 1}$  and the eigenfunctions $\{y_j(\bx)\}_{j \ge 1}$ form an 
orthonormal basis of $L^2(\Omega)$. Functional calculus on  $\cA$ is defined as usual: If $\Phi:\mathbb{C} \mapsto \mathbb{C}$ is a continuous function, then 
$\Phi(\cA)$ is the self-adjoint operator with the same eigenfunctions as $\cA$ and the eigenvalues $\{\Phi(\theta_j)\}_{j \ge 1}$.

Our data  mapping is stated in the next lemma. Its purpose is twofold: First, it gives a Duhamel type principle, which maps the forcing in the wave equation to an initial condition. Second, it leads to an inner product expression of the data entries, which is used in the computation of the  ROMs.

\vspace{0.05in} 
\begin{lem} 
\label{lem.1}
Define the new data matrix $\bD(t)$ by the mapping 
\begin{equation}
\boldsymbol{\mathcal{M}}(t) \mapsto \bD(t)  =  \boldsymbol{\mathcal{M}}^f(t) + \boldsymbol{\mathcal{M}}^{f}(-t), \quad \mbox{with} ~\boldsymbol{\mathcal{M}}^f(t) = -f'(-t) \star_t \boldsymbol{\mathcal{M}}(t),\label{eq:D5}
\end{equation}
\bc{for $t > 0$.}
Its entries have the integral (inner) product expression
\begin{align}
D_{r,s}(t) &= 
\int_{\Omega} d \bx \, u_0^{(r)}(\bx) u^{(s)}(t,\bx), \label{eq:defDrs} \qquad s,r = 1,\ldots, m,
\end{align}
where 
\begin{equation}
u^{(s)}(t,\bx) = \cos \big(t \sqrt{\cA} \, \big) u^{(s)}_0(\bx) 
\label{eq:D7}
\end{equation}
is the solution of the homogeneous wave equation 
\begin{align}
(\partial_t^2 + \cA) u^{(s)}(t,\bx) &= 0, \qquad t > 0, ~~ \bx \in \Omega, \label{eq:D8}
\end{align}
with initial state
\begin{align}
u^{(s)}(0,\bx) &= u_0^{(s)}(\bx) =  \left|\hat f\big(\sqrt{\cA} \, \big) \right|\delta_{\bx_s}(\bx), \qquad 
\partial_t u^{(s)}(0,\bx) = 0.
\label{eq:D7_0}
\end{align}
\end{lem}

\vspace{0.05in}
\textbf{Proof:} The proof is given in \cite[Appendix A]{borcea2020reduced} but we include it here, with more detailed explanation,  for convenience of the reader. Due to the convolution by $f'(-t)$ in \eqref{eq:D5}, the $(r,s)$ entries 
of $\boldsymbol{\mathcal{M}}^f(t)$ are $\mathcal{M}_{r,s}^f(t) = P^{(s),f}(t,\bx_r),$
where 
\begin{align}
P^{(s),f}(t,\bx) = -f'(-t) \star_t P^{(s)}(t,\bx) = F(t) \star_t  H(t) \cos \big( t \sqrt{\cA} \,  \big) \delta_{\bx_s}(\bx), \label{eq:PsEq}
\end{align}
solves the wave equation  with right hand side $F'(t) \delta_{\bx_s}(\bx)$. Recall definition \eqref{eq:defF} of $F(t)$. Here $H(t)$  is the Heaviside step function and the second term in the convolution is the time derivative of the causal Green's function of the wave operator $\partial_t^2 + \cA$. Functional calculus gives
\begin{equation}
\cos \big( t \sqrt{\cA} \,  \big) \delta_{\bx_s}(\bx) =  \int_0^\infty d \theta \, 
\cos(t \sqrt{\theta}) \rho(\theta,\bx,\bx_s), \label{eq:funcCalc}
\end{equation}
where $\rho(\theta,\bx,\bx_s)= \sum_{j \ge 1} \delta(\theta-\theta_j) y_j(\bx_s) y_j(\bx) $ is the discrete spectral measure density associated with $\cA$.

Substituting \eqref{eq:funcCalc} into \eqref{eq:PsEq} and writing the convolution  via the Fourier transform, we obtain
\begin{equation}
P^{(s),f}(t,\bx) = \int_0^\infty d \theta \left[ \frac{1}{2} \cos(t \sqrt{\theta}) \hat F(\sqrt{\theta}) + 
\int_{\RR} \frac{d \om}{2 \pi} \frac{i \om \hat F(\om)}{(\theta -\om^2)}e^{-i \om t} \right] \rho(\theta,\bx,\bx_s).
\label{eq:D5p}
\end{equation}
The new $(r,s)$ data entries in \eqref{eq:D5} are:  
$
D_{r,s}(t) = P^{(s),f}(t,\bx_r) + P^{(s),f}(-t,\bx_r).$
The motivation for adding the  negative $t$ term  is that   the even wave
\begin{align}
P_e^{(s),f}(t,\bx) &= P^{(s),f}(t,\bx) + P^{(s),f}(-t,\bx) =  \int_0^\infty d \la \cos(t \sqrt{\theta}) \hat F(\sqrt{\theta}) \rho(\theta,\bx,\bx_s)\nonumber \\ &= \cos \big(t \sqrt{\cA} \, \big) \hat F\big(\sqrt{\cA} \, \big) \delta_{\bx_s}(\bx), \qquad \bc{t >0, ~~ \bx \in \Omega},
\label{eq:D6}
\end{align}
has a simpler  expression  than \eqref{eq:D5p}, because
\begin{align*}
\int_{\RR} \frac{d \om}{2 \pi} \frac{ \om \hat F(\om)}{(\theta -\om^2)} \left(e^{-i \om t} + e^{i \om t} \right)  = \int_{\RR} \frac{d\om}{\pi} \frac{\om |\hat f(\om)|^2}{(\theta - \om^2)} \cos(\om t)  = 0.
\end{align*}
Moreover,  \eqref{eq:D6} solves  the homogeneous wave equation (\ref{eq:D8}) with  initial condition
\begin{equation*}
P_e^{(s),f}(0,\bx) = \hat F\big(\sqrt{\cA} \, \big) \delta_{\bx_s}(\bx), \qquad \partial_t P_e^{(s),f}(0,\bx) = 0.
\end{equation*}
It remains to connect  \eqref{eq:D6}  to the wave defined in \eqref{eq:D7}. To do so, we take the square root of 
$\hat F(\om) = |\hat f(\om)|^2$  and use that functions of $\cA$ commute. We get that  
\begin{align}
P_e^{(s),f}(t,\bx) &= \left|\hat f\big(\sqrt{\cA} \, \big)\right| \underbrace{\cos \big(t \sqrt{\cA} \, \big) \left| \hat f\big(\sqrt{\cA} \, \big)\right| \delta_{\bx_s}(\bx) }_{u^{(s)}(t,\bx)}, \qquad \bc{t >0, ~~ \bx \in \Omega}.
\label{eq:PtoU}
\end{align}
The statement of the lemma follows after substituting \eqref{eq:PtoU}  into the expression of $D_{r,s}(t)$  and using that $\cA$ is self-adjoint
\begin{align}
D_{r,s}(t) &= P_e^{(s),f}(t,\bx_r) = \int_{\Omega} d \bx \, \delta_{\bx_r} (\bx)  \left|\hat f\big(\sqrt{\cA} \, \big)\right| u^{(s)}(t,\bx) 
\label{eq:D6p} \\&= 
\int_{\Omega} d \bx \,  \underbrace{\left|\hat f \big(\sqrt{\cA} \, \big)\right| \delta_{\bx_r} (\bx)}_{u_0^{(r)}(\bx)} u^{(s)}(t,\bx) .  \nonumber \qquad  \Box
\end{align}

Note that the transformation \eqref{eq:D5} can be carried out without any knowledge of the medium, if the measurements 
start at $t_{\rm min} = -t_F$ (or earlier), where $(-t_F,t_F)$ is the support of $F(t)$.  If the measurements start later, at $t_{\rm min} > - t_F$, we need the assumption 
that $c(\bx) = \bar{c}\, $ is known near the array, to compute the missing measurements in the interval $(-t_F,t_{\rm min})$. Near means within the distance $\bar c (t_{\rm min} + t_F)$ from the sources/receivers.

\subsection{The snapshots and the propagator}
\label{sect:ROM.ts} 
To define the ROMs, we will  use the snapshots of the wave \eqref{eq:D7} on a time grid with step $\tau$. The choice of $\tau$ matters, as explained in \cite[Section 6 \& Appendix A]{borcea2021reduced} and in \cite{borcea2022internal,borcea2022waveform}. 
For our purpose, it suffices to say that it 
should be $\tau \sim 2 \pi/\om_o,$ where $\sim$ means equal, up to a finite constant. A larger $\tau$ does not sample the wave field according to the Nyquist criterium and gives poor approximation properties of the ROM. A smaller $\tau$ leads to poor conditioning of the 
mass matrix defined below and requires regularization. \bc{We comment briefly on regularization in 
section \ref{sect:regul}. More details are in  
\cite[Appendix E]{borcea2022waveform}.}

Let us gather the snapshots  for all the $m$ sources in the row vector fields 
\begin{equation}
\bu_j(\bx) = \left( u^{(1)}(j \tau,\bx), \ldots, u^{(m)}(j \tau,\bx) \right), \qquad j \ge 0, ~~ \bc{\bx \in \Omega},
\label{eq:S1}
\end{equation}
called henceforth the ``vector snapshots". According to \eqref{eq:D7}, they are given by 
\begin{equation}
\bu_j(\bx) = \cos \big(j \tau \sqrt{\cA} \, \big) \bu_0(\bx), \qquad j \ge 0, ~~ \bc{\bx \in \Omega},
\label{eq:S2}
\end{equation}
and thanks to the trigonometric identity 
$
\cos [(j+1)\alpha] + \cos[|j-1|\alpha] = 2 \cos(\alpha) \cos(j \alpha),$  $\forall \alpha \in \RR,
$
they satisfy the exact time stepping scheme
\begin{equation}
\bu_{j+1}(\bx) = 2 \cP \bu_{j}(\bx) - \bu_{|j-1|}(\bx), \qquad  \cP = \cos \big(\tau \sqrt{\cA} \, \big),
\label{eq:defProp}
\end{equation}
for $j \ge 0$ and $\bc{\bx \in \Omega}$, driven by the propagator operator $\cP$.

The ROMs will be computed from the \bc{$m \times m$  data matrices $\bD_j = \bD(j\tau)$, defined   in \eqref{eq:D5} and evaluated at the time instants 
$t_j = j \tau$, for $j = 0, \ldots, 2n-1$}. According to Lemma \ref{lem.1}, these have the integral expression
\begin{equation}
\bD_j =  \int_{\Omega} d \bx \, \bu_0^T(\bx) \bu_j(\bx) = \int_{\Omega} d \bx \, \bu_0^T(\bx) \cos \big( j \tau \sqrt{\cA}\, \big) \bu_0(\bx),
\label{eq:dataSnap}
\end{equation}
with kernel given by the Chebyshev polynomial $T_j$ of the first kind  \cite{gradshteyn2014table} of the propagator operator:
$
 \cos(j \tau \sqrt{\cA}\, ) = \cos( j \arccos \cP) = T_j ( \cP).
$

\subsection{Data driven Galerkin approximations}
\label{sect:ROM.2}
We consider two Galerkin approximations of wave propagation, 
\bc{in the finite-dimensional subspace $\mathscr{U}$ of $L^2(\Omega)$ generated by the $nm$ snapshots (for all $m$ sources and $n$ time  instants):}
\begin{equation}
\bc{\mathscr{U} = \mbox{range} \left\{  (u^{(s)}(j\tau,\bx))_{\bx \in \Omega} ,\, j=0,\ldots,n-1, \, s=1,\ldots,m \right\} . }
\label{eq:G1}
\end{equation}
We assume that  
$\mbox{dim}(\mathscr{U}) = nm$. This holds  in general if $\tau \sim 2 \pi/\om_o$ and the separation between the sources is $\sim 2\pi \bar{c}/\om_o$. If the dimension of $\mathscr{U}$ is smaller than $nm$, then the ROM construction requires regularization \bc{(section \ref{sect:regul} and  \cite[Appendix E]{borcea2022waveform})}.

\bc{Let us denote by 
\begin{equation}
\bU(\bx) =  \big(\bu_0(\bx), \ldots, \bu_{n-1}(\bx)\big), \quad \bx \in \Omega,
\label{eq:G1b}
\end{equation}
the $nm$ dimensional row vector field that contains all the snapshots.}
The first Galerkin approximation is for the time stepping equation \eqref{eq:defProp}. It approximates the vector snapshots \eqref{eq:S1}  
by 
\begin{equation}
\bu_j^{\rm Gal}(\bx) = \bU(\bx) \bg_j, \qquad j \ge 0, ~~\bc{\bx \in \Omega},
\label{eq:G2}
\end{equation}
using the Galerkin coefficient matrices $\bg_j \in \RR^{nm\times m}$ defined so that the residual is orthogonal to the space $\mathscr{U}$.
To write this explicitly, we use an orthonormal basis of $\mathscr{U}$, stored in the $nm$ dimensional row vector field 
\begin{equation}
\label{eq:orthonV}
\bV(\bx) = \big(\bv_0(\bx), \ldots, \bv_{n-1}(\bx) \big), ~~\bc{\bx \in \Omega}, \quad \mbox{with} ~~\int_{\Omega} d \bx\,  \bV^T(\bx)  \bV(\bx) = \bI_{nm}.
\end{equation} 
Here $\bI_{nm}$ is the $nm \times nm$ identity matrix. As we did for $\bU(\bx)$, we organize the entries of $\bV(\bx)$  in the $m$ dimensional row vectors $\bv_j(\bx)$, associated with the instants $j \tau$, for $j = 0, \ldots, n-1$.  This makes sense because our basis  
has the ``causal" property
\begin{equation}
\bc{ (\bv_j(\bx))_{\bx \in \Omega} \in \mbox{range} \left\{  (\bu_{j'}(\bx))_{\bx \in \Omega} , \, j'=0,\ldots,j\right\},   }
 \label{eq:causalVs}
\end{equation}
\bc{for any $j = 0, \ldots, n-1 $.}
It is obtained by the  Gram-Schmidt orthogonalization of the components of $\bU(\bx)$, 
\begin{equation}
\bU(\bx) = \bV(\bx) \bR, \qquad \bc{\bx \in \Omega},
\label{eq:G17}
\end{equation}
where $\bR$ is an invertible $nm \times nm$ matrix with block upper triangular structure. 

The orthogonality of the residual to $\mathscr{U}$ can now be written as  
\begin{equation}
\int_\Omega d \bx\,  \bV^T(\bx) \left[ \bU(\bx) \bg_{j+1} + \bU(\bx) \bg_{|j-1|} - 2 \cP \bU(\bx) \bg_j \right] = {\bf 0}, \qquad \forall j \ge 0,
\label{eq:G3}
\end{equation}
and solving for $\bV(\bx)$ from \eqref{eq:G17} we get the following time stepping scheme for the Galerkin coefficients, 
\begin{equation}
\bR^{-T} \left[\bM (\bg_{j+1} + \bg_{|j-1|}) - 2 \bS \bg_j \right] = 0, \qquad j \ge 0.
\label{eq:G4}
\end{equation}
Here  $\bR^{-T}$ denotes the transpose of $\bR^{-1}$,   and $\bM$ and $\bS$ are the $nm \times nm$ ``mass" (Gramian) and ``stiffness" matrices
\begin{equation}
\bM = \int_\Omega d \bx\,  \bU^T(\bx) \bU(\bx), \qquad \bS = \int_\Omega d \bx\,  \bU^T(\bx) \cP \bU(\bx).
\label{eq:G5}
\end{equation}

The second Galerkin approximation is for the wave equation \eqref{eq:D8}. To distinguish it from the first, we use the tilde in the notation 
of the approximate wave field 
\begin{equation}
\tilde \bu^{\rm Gal}(t,\bx) = \bU(\bx) \tilde\bg(t), \qquad t \ge 0, ~~ \bc{\bx \in \Omega},
\label{eq:G7}
\end{equation}
where the Galerkin coefficients $\tilde \bg(t)$ are $nm \times m$ time dependent matrices. These are defined 
so that when substituting \eqref{eq:G7} into the wave equation (\ref{eq:D8}), the residual is orthogonal to the space $\mathscr{U}$ i.e., 
\begin{equation}
\int_\Omega d \bx\,  \bV^T(\bx) \left[ \bU(\bx) \tilde \bg''(t) + \cA \bU(\bx) \tilde \bg(t) \right] = {\bf 0}, \qquad t \ge 0,
\label{eq:G8}
\end{equation}
where $\bg''(t) $ denotes the second derivative of $\bg(t) $.
Using again equation \eqref{eq:G17}, we obtain  the semi-discretized wave equation in the Galerkin framework
\begin{equation}
\bR^{-T} \left[\bM  \tilde \bg''(t) + \tilde \bS \tilde \bg(t)\right] = {\bf 0}, \qquad t \ge 0,
\label{eq:G9}
\end{equation}
where the mass matrix $\bM$ is the same as in \eqref{eq:G5}, but the stiffness matrix is 
\begin{equation}
\tilde \bS = \int_{\Omega} d \bx \, \bU^T(\bx) \cA \bU(\bx).
\label{eq:G10}
\end{equation}

The Galerkin approximations described above are standard, except for the multiplication by $\bR^{-T}$ of equations 
\eqref{eq:G4} and \eqref{eq:G9}. Since $\bR$ is nonsingular\footnote{\bc{Recall that here we assume linearly independent 
snapshots i.e., $\rm{dim}{\mathscr{U}} = nm$. Otherwise, we use regularization to obtain an invertible approximation of $\bR$, as explained in section \ref{sect:regul}.}}, this multiplication does not change the solution of these equations. We only use it to get a good algebraic structure of our ROMs, as explained in the next section. 
The question is, how can we use these Galerkin approximations when we do not know the snapshots in $\bU(\bx)$? The next two theorems show that, in fact, equations \eqref{eq:G4} and \eqref{eq:G9} are array data driven. 

\vspace{0.05in}\begin{thm}
\label{thm.1}
Let $\be_j $ denote the  $nm \times m$ column blocks of the identity matrix i.e., 
 $\bI_{nm} = \big(\be_0, \be_1, \ldots, \be_{n-1} \big)$. The  Galerkin coefficients in approximation \eqref{eq:G2} satisfy
\begin{equation}
\label{eq:G11}
\bg_j = \be_j, \qquad  j = 0, \ldots, n-1.
\end{equation}
The $m \times m$ blocks of the mass and stiffness matrices defined in \eqref{eq:G5} are given by the first $2n-1$ data matrices \eqref{eq:dataSnap} as follows:
\begin{align}
\hspace{-0.15in}\bM_{i,j} &= \int_{\Omega} d \bx \, \bu_i^T(\bx) \bu_j(\bx) = \frac{1}{2}\big(\bD_{i+j} + \bD_{|i-j|} \big),  \label{eq:G12} \\
\hspace{-0.15in}\bS_{i,j} &= \int_{\Omega} d \bx \, \bu_i^T(\bx) \cP \bu_j(\bx)= \frac{1}{4}\big( \bD_{i+j+1} + \bD_{|i-j+1|}+\bD_{|i+j-1|} + \bD_{|i-j-1|} \big),
\label{eq:G13}\end{align}
for $i,j = 0, \ldots, n-1$. The orthonormal basis $\bV(\bx)$ cannot be computed without knowing $\bU(\bx)$, but the block upper triangular matrix $\bR$ in its definition \eqref{eq:G17} is data driven. It is the block Cholesky factor of the mass matrix:
$ ~
\bR^T \bR = \bM.
$
\end{thm}

\vspace{0.05in}Note that there is an ambiguity in the definition of the block Cholesky factorization. The difference between various factorization algorithms is in the computation of the diagonal blocks of $\bR$, which involves taking the square root of a symmetric, positive 
definite $m \times m$ matrix. Any algorithm will do, as long as it is used consistently throughout the inversion procedure. We  use [16, Algorithm 5.2], which takes the square root using the spectral decomposition.

\vspace{0.05in}
\textbf{Proof of Theorem \ref{thm.1}:}  \bc{Equation \eqref{eq:G11} seems natural, because 
the first $n$ snapshots define the approximation space. However, for the result to hold, we also need that
these snapshots satisfy exactly the time stepping equation, which is indeed the case. The proof of \eqref{eq:G11} is as follows:}
Since $\mbox{dim} (\mathscr{U}) = nm$, the matrices $\bR$ and $\bM$ are nonsingular.
Definition \eqref{eq:G1b} of $\bU(\bx)$ and equation \eqref{eq:defProp} give that, if equation 
\eqref{eq:G11} holds, then the residual is 
\begin{align*}
\bU(\bx) (\bg_{j+1} +  \bg_{|j-1|}) - 2 \cP \bU(\bx) \bg_j &=  \bU(\bx) (\be_{j+1} + \be_{|j-1|}) - 2 \cP \bU(\bx) \be_j   \\
&= \bu_{j+1}(\bx) + \bu_{|j-1|} (\bx)- 2 \cP \bu_j(\bx) = 0,
\end{align*}
for $j = 0, \ldots, n-1$. Obviously, this residual satisfies equation \eqref{eq:G3}, which is equivalent to \eqref{eq:G4}. But equation \eqref{eq:G4}
with initial conditions $\bg_0 = \be_0$ and $\bg_{-1} = \bg_1$ 
 has a unique solution, since $\bR^{-T} \bM$ is invertible, so the  Galerkin coefficients must satisfy \eqref{eq:G11}. Note that we use these initial 
 conditions to ensure that the Galerkin approximation is exact at $t = 0$ and that it is even in time.

To prove \eqref{eq:G12}, we use the expression \eqref{eq:S2} of the snapshots in the definition of the $(i,j)$ block of $\bM$,
\begin{align*}
\bM_{i,j} &= \int_{\Omega} d \bx \, \left\{ \cos \big(i \tau \sqrt{\cA} \, \big) \bu_0(\bx) \right\}^T \cos \big(j\tau \sqrt{\cA} \, \big) \bu_0(\bx) \\
&= \int_{\Omega} d \bx \, \bu_0^T(\bx) \cos \big(i \tau \sqrt{\cA} \, \big)  \cos \big(j\tau \sqrt{\cA} \, \big) \bu_0(\bx) \\
&= \int_{\Omega} d \bx \, \bu_0^T(\bx) \frac{1}{2} \left\{ \cos \big[(i+j) \tau \sqrt{\cA} \, \big] + \cos \big[|i-j| \tau \sqrt{\cA} \, \big]    \right\} \bu_0(\bx),
\end{align*}
where the second equality is because $\cA$ is self-adjoint and the last equality is by the trigonometric identity: 
$
2\cos(i \alpha)\cos(j \alpha) = \cos\big[(i+j)\alpha\big] + \cos\big[|i-j|\alpha\big],$ $ \forall \, \alpha \in \RR.$
The result \eqref{eq:G12} follows from  equation \eqref{eq:dataSnap}.

The calculation of $\bS$ is similar, because the time stepping scheme gives
\[
\cP \bu_j(\bx) = \frac{1}{2} \left[ \bu_{j+1}(\bx) + \bu_{|j-1|}(\bx) \right], \qquad \bc{j \ge 0, ~~ \bx \in \Omega}.
\]
Substituting this into the definition of $\bS_{i,j}$ and using the calculation above we get \eqref{eq:G13}. 
Finally,  we deduce from the Gram-Schmidt orthogonalization 
\eqref{eq:G17} and the definition \eqref{eq:G5} of the mass matrix that 
\begin{equation}
\bM = \int_{\Omega} d \bx \, \bU^T(\bx) \bU(\bx) = \bR^T \int_{\Omega} d \bx \, \bV^T(\bx) \bV(\bx) \bR = \bR^T \bR. \quad \Box
\label{eq:Chol}
\end{equation}

The second Galerkin approximation is on the same space $\mathscr{U}$, using  the same basis in $\bV(\bx)$  
and the same mass matrix $\bM$. The next theorem describes the other terms in the semi-discrete wave equation \eqref{eq:G9}.

\vspace{0.05in}
\begin{thm}
\label{thm.2}
The time dependent Galerkin coefficients in the approximation \eqref{eq:G7} satisfy the initial conditions
\begin{equation}
\tilde \bg(0) = \be_0, \qquad \tilde \bg'(0) = {\bf 0}.
\label{eq:G14}
\end{equation}
The $m \times m$ blocks of the stiffness matrix $\tilde \bS$ are given by 
\begin{equation}
\tilde \bS_{i,j} = \int_{\Omega} d \bx \, \bu_i^T(\bx) \cA \bu_j(\bx) = - \frac{1}{2} \left[ \ddot \bD_{i+j} + \ddot \bD_{|i-j|} \right], 
\qquad i,j = 0, \ldots, n-1,
\label{eq:G15}
\end{equation}
where $\ddot \bD_j$ denotes the second derivative of $\bD(t)$ evaluated at $t = j \tau$.
\end{thm}

\vspace{0.05in}
\textbf{Proof:} Equation \eqref{eq:G14} ensures that the Galerkin approximation \eqref{eq:G7} satisfies exactly the initial conditions
\begin{equation}
\label{eq:ROMAini}
\tilde \bu^{\rm Gal} (0,\bx) = \bU(\bx) \be_0 = \bu_0(\bx), \quad \partial_t \tilde \bu^{\rm Gal} (0,\bx) = \bU(\bx)  \bg'(0) = {\bf 0}, 
\qquad \bc{\bx \in \Omega}.
\end{equation}
The expression \eqref{eq:G15} of the stiffness matrix follows from the wave equation
\[
\partial_t^2 \bu(t = j \tau, \bx) = - \cA \bu_j(\bx),\qquad \forall j \ge 0, ~~ \bc{\bx \in \Omega},\]
and the definition of $\bD(t)$ given in Lemma \ref{lem.1}. $~\Box$

\subsection{The data driven  ROMs}
\label{sect:ROM.4}
We are now ready to define our two ROMs from the data driven Galerkin equations \eqref{eq:G4} and \eqref{eq:G9}.

\subsubsection{The ROM propagator}
\label{sect:ROMProp}
\bc{This ROM is derived from equation \eqref{eq:G4}, and defines a discrete time dynamical system, with state at time instant $t_j = j \tau$ given by the $j^{\rm th}$ ROM snapshot, which is the 
$nm \times m$ matrix }
\begin{equation}
\bu_j^\RM = \int_{\Omega} d \bx \, \bV^T(\bx) \bu_j^{\rm Gal}(\bx) = \int_{\Omega} d \bx \, \bV^T(\bx) \bU(\bx) \bg_j = 
\bR \bg_j, \qquad j \ge 0.
\label{eq:G20}
\end{equation}
The second equality in this equation is by definition \eqref{eq:G2}  and the third is due to the Gram-Schmidt equation \eqref{eq:G17}. 
The ROM snapshots are data driven: The first $n$ of them are just the $nm \times m$ block columns of  $\bR$, due to 
equation \eqref{eq:G11},
\begin{equation}
\bu_j^\RM = \bR \be_j, \qquad j = 0, \ldots, n-1.
\label{eq:G21}
\end{equation}
The others are obtained by time stepping in  equation \eqref{eq:G4}, for $j \ge n-1$.

\vspace{0.05in}
\begin{rem} The first ROM snapshot has nonzero entries only in the first $m \times m$ block. This is the algebraic way of capturing that 
 $\bu(t=0,\bx)$ is supported near the array. As $j$ increases, the row blocks of $\bu_j^\RM$ fill in sequentially, thus capturing the 
progressive advancement of the wave away from the array. This physical interpretation is easier to visualize in the one-dimensional 
case, where we can transform the space coordinate to time using the travel time map $x \mapsto \int_0^x d x'/c(x')$. We work in higher dimensions,
where such a travel time transformation cannot be done. Nevertheless, the interpretation remains formally true, in the sense that 
we can associate to each instant $j \tau$ a  maximum ``depth coordinate" reached by the wave, while the dependence on the 
other spatial coordinates is comprised in the $m \times m$ blocks of \eqref{eq:G20}.
\end{rem}

\vspace{0.05in}
The time stepping  scheme for the ROM snapshots is obtained from equations \eqref{eq:G4}, \eqref{eq:Chol}  and the definition 
\eqref{eq:G20},
\begin{equation}
\bu_{j+1}^\RM = 2 \cP^\RM \bu_j^\RM - \bu_{|j-1|}^\RM, \qquad j \ge 0.
\label{eq:G22}
\end{equation}
It is controlled by the ROM propagator $\cP^\RM$, the symmetric $nm \times nm$ matrix defined as the Galerkin projection of the propagator operator  \eqref{eq:defProp},
\begin{equation}
\cP^\RM = \int_{\Omega} d \bx \, \bV^T(\bx) \cP \bV(\bx).
\label{eq:G23}
\end{equation}
This ROM propagator is data driven thanks to equation \eqref{eq:G17} and definition \eqref{eq:G5},
\begin{equation}
\cP^\RM = \bR^{-T}  \int_{\Omega} d \bx \, \bU^T(\bx) \cP \bU(\bx)\bR^{-1} = \bR^{-T} \bS \bR^{-1}.
\label{eq:G24}
\end{equation} 

\subsubsection{The ROM of the wave operator}
\label{sect:ROMA}
The ROM derived from the Galerkin approximation \eqref{eq:G7}  describes the evolution of  the $nm \times m$ valued  ROM wave 
\begin{equation}
\tilde \bu^\RM(t) = \int_{\Omega} d \bx \, \bV^T(\bx) \tilde \bu^{\rm Gal}(t,\bx) = \int_{\Omega} d \bx \, \bV^T(\bx)  \bU(\bx) \tilde \bg(t) = \bR \tilde \bg(t),
\label{eq:ROMA_sn}
\end{equation}
where we used again equation \eqref{eq:G17}. 
This wave satisfies the initial conditions 
\begin{equation}
\tilde \bu^\RM(0) = \bR \be_0, \qquad \frac{d}{dt} \tilde \bu^\RM(0) = {\bf 0},
\label{eq:G25}
\end{equation}
thanks to equation \eqref{eq:G14},  
and it evolves according to the system of ODEs 
\begin{equation}
\frac{d^2}{dt^2} \tilde \bu^{\RM}(t) + \bAR \tilde \bu^{\RM}(t) = {\bf 0}, \qquad t \ge 0.
\label{eq:G26}
\end{equation}
The wave operator $\partial_t^2 + \cA$ is now replaced by the operator $\frac{d^2}{dt^2} + \bAR$, 
where 
\begin{equation}
\bAR = \int_{\Omega} d \bx \, \bV^T(\bx) \cA \bV(\bx),
\label{eq:G27}
\end{equation}
is the Galerkin projection of $\cA$, a symmetric $nm \times nm$ matrix. Again, this can be computed from the data thanks to the Gram-Schmidt formula 
\eqref{eq:G17}, 
\begin{equation}
\bAR = \bR^{-T} \int_{\Omega} d \bx \, \bU^T(\bx) \cA \bU(\bx) \bR^{-1} =
\bR^{-T} \tilde \bS \bR^{-1}.
\label{eq:G28}
\end{equation}
Here we used definition \eqref{eq:G10} of $\tilde \bS$, which can be obtained from the measurements as stated in Theorem \ref{thm.2}.

\subsection{Comparison of the ROMs}
\label{sect:ROM.5}
How do our ROMs compare with the time domain reduced basis ROMs found in the computational science literature
\cite{benner2015survey,hesthaven2022reduced,lieu2006reduced}, which also use snapshots to define the projection space? 

The philosophy behind the reduced basis ROMs stems from  principle component analysis in multivariate statistics \cite{jollife2016principal}
and the Karhunen-Lo\`{e}ve decomposition in stochastic processes modeling 
\cite{karhunen1947,loeve1978}.
They explore a dynamical system, the wave equation for us,  from knowledge of snapshots of the solution. These snapshots are usually discretized in space and the idea is to 
extract from them a set of uncorrelated vectors, the so-called proper orthogonal decomposition (POD) modes. This can be done via 
singular value decomposition and the reduced basis is given by the modes corresponding to the significant singular values. 

The POD reduced basis ROMs are not useful for waveform inversion, because they require more data than what we can measure at the array. We are also not interested in compressing information. We want to learn how the waves propagate, using the two array data driven ROMs, and then figure out how to estimate the wave speed $c(\bx)$ from them. For the latter task the algebraic structure of the ROM is very important. It is  easier to compute a ROM that approximates the forward map. In particular,  in our Galerkin approximations we could have used any  basis of  $\mathscr{U}$ and we would have obtained a good ROM for approximating this map. But for the inverse problem it is important to use the basis $\bV(\bx)$ which gives the causal algebraic structure of the ROMs. It is not difficult to see from  
Theorems \ref{thm.1} and \ref{thm.2} and equation \eqref{eq:causalVs}  that if we restrict the data to the time instants 
$j = 0, \ldots, 2 J-1$, with $J < n$, then the ROMs will be sensitive to the parts of the medium reached by the wave field  up to time $J \tau$. Thus, we can use the ROMs to estimate $c(\bx)$ in a ``layer peeling" fashion, first near the array and then deeper inside the medium, by increasing $J$ gradually. This helps significantly  the inversion.

Can we say more about the algebraic structures of the ROMs? We can deduce from the ROM time stepping scheme \eqref{eq:G22}
evaluated at $j = 0, \ldots, n-1$ and from equation \eqref{eq:G21} that the ROM propagator $\cP^\RM$ is a  block tridiagonal matrix. A proof of this fact is in \cite[Appendix C]{borcea2020reduced}. The second ROM matrix $\bAR$  does not have the same sparse structure, although its entries decay away from the main diagonal, as proved in  \cite[Appendix D]{borcea2022waveform}. 

The ROM propagator described in section \ref{sect:ROMProp} has surprisingly good approximation properties, meaning that even though it is defined via projection on the space 
spanned by the first $n$ vector snapshots, it fits the data $\bD_j$ for $j$ up to $2n-1$.
This is proved in \cite[Appendix B]{borcea2020reduced}. We 
give here a reformulation of this result, in a slightly weaker form, that is easier to explain and is used in the next section.

\vspace{0.05in}
\begin{thm}
\label{lem.2}
The propagator ROM snapshots satisfy the data fit relations 
\begin{align}
\bD_j &= \int_{\Omega} d \bx \, \bu_0^T(\bx) \bu_j(\bx) = (\bu_0^{\RM})^T \bu_j^\RM, \label{eq:DF1}
\\
\bD_{j+n-1} &= \int_{\Omega} d \bx \, [2 \bu_{n-1}^T(\bx) \bu_j(\bx)  - \bu_0^T(\bx) \bu_{n-1-j}(\bx)] 
\nonumber\\
{}&= 2 (\bu_{n-1}^{{\RM}})^T\bu_j^\RM -(\bu_0^\RM)^T \bu_{n-1-j}^\RM, 
\label{eq:DF2}
\end{align}
for $j = 0, \ldots, n-1$.
\end{thm}

\vspace{0.05in}
\textbf{Proof:} The first equality in \eqref{eq:DF1} is from equation \eqref{eq:dataSnap}. The second equality is 
because by definition  \eqref{eq:G1b} of $\bU(\bx)$, the Gram-Schmidt equation \eqref{eq:G17} and the expression 
\eqref{eq:G21} of the first $n$ ROM snapshots, we get 
\begin{equation}
\bu_j(\bx) = \bU(\bx) \be_j = \bV(\bx) \bR \be_j = \bV(\bx) \bu_j^\RM, \qquad \bc{\bx \in \Omega}, \quad j = 0, \ldots, n-1,
\label{eq:DF3}
\end{equation}
and therefore, 
\begin{align}
\nonumber
\bD_j &=  \int_{\Omega} d \bx \big[ \bV(\bx) \bu_0^\RM\big]^T \bV(\bx) \bu_j^\RM  = (\bu_0^\RM)^T \int_{\Omega} d \bx  \bV^T(\bx) \bV(\bx) \bu_j^\RM  \\
&\stackrel{\eqref{eq:orthonV}}{=}
(\bu_0^{\RM})^T  \bu_j^\RM.   \label{eq:DF4}
\end{align}

To prove equation \eqref{eq:DF2}, we recall the calculation of the mass matrix given in Theorem \ref{thm.1}. 
Setting $i = n-1$ and $j = 0, \ldots, n-1$  in equation \eqref{eq:G12},  we get 
\begin{equation}
\bD_{j+n-1} = 2 \int_{\Omega} d \bx \, \bu_{n-1}^T(\bx) \bu_j(\bx) - \bD_{n-1-j},
\end{equation}
and the result follows from \eqref{eq:DF3} and \eqref{eq:DF1}. $~~\Box$

The wave operator ROM described in section \ref{sect:ROMA} does not fit the data exactly, it only approximates it. 
This is because the  Galerkin wave approximation \eqref{eq:G7} from which it is derived is exact only at $t = 0$, 
as stated in \eqref{eq:ROMAini}, but not at $t = j \tau$ for $1 \le j \le n-1$.  To get an exact fit we would have 
to include the snapshots of the second time derivative of the wave field in the definition of the approximation space.  
Nevertheless,  we will see in the next section that 
$\bAR$  is useful for estimating $c(\bx)$, in spite of the inexact data fit.

\subsection{\bc{Regularization of the ROM computation}}
\label{sect:regul}
\bc{The ideal case of noiseless data and linearly independent snapshots that define the approximation space in equation \eqref{eq:G1}  guarantees a positive definite 
mass matrix $\bM$ whose square root $\bR$ can be inverted.  In the presence of noise and/or too small time steps $\tau$ and sensor distance separation, the mass matrix $\bM$ will likely be ill-conditioned, singular or indefinite. Therefore, regularization is needed to construct the ROMs. }

\bc{Our approach to regularization is based on spectral projection: Let $\{\la_j\}_{j=1}^{nm}$ be the eigenvalues of $\bM$, in decreasing order, and $\{\by_j\}_{j=1}^{nm}$ the corresponding eigenvectors.  Let $\epsilon$ be a positive threshold
and $r$ be the smallest natural number such that $\la_j < \epsilon$ for $rm < j \le nm$. Define the projected mass  matrix 
\begin{equation}
\bLa^\epsilon = (\bY^\epsilon)^T \bM \bY^{\epsilon} \in \mathbb{R}^{rm \times rm}, \qquad \bY^{\epsilon} = \left(\by_1, \ldots, \by_{rm}\right) \in \mathbb{R}^{nm \times rm}.
\end{equation}
This matrix  is obviously diagonal. 
The threshold $\epsilon$ should be small, while ensuring a reasonable condition number $\la_1/\epsilon$ of $\bLa^\epsilon$. 
}

\bc{
Note that  $\bLa^\epsilon$ cannot be our 
regularized mass matrix because it does not have the correct block Hankel+Toeplitz structure, which  captures the causal wave propagation. To determine the orthogonal transformation that restores this structure, we recall that the ROM propagator must be block tridiagonal. Thus, instead of \eqref{eq:G24}, we compute
\begin{equation}
\boldsymbol{\Pi}^\epsilon  = (\bLa^\epsilon)^{-1/2} (\bY^\epsilon)^T \bS \bY^\epsilon (\bLa^\epsilon)^{-1/2}\in \mathbb{R}^{rm \times rm},
\end{equation}
and then use the block Lanczos algorithm \cite{golubVanLoan} to generate the orthogonal matrix ${\itbf Q}^\epsilon \in \RR^{rm \times rm}$ that gives the block tridiagonal regularized ROM propagator
\begin{equation}
\cP^{\epsilon,\RM} = ({\itbf Q}^\epsilon)^T  
\boldsymbol{\Pi}^\epsilon 
{\itbf Q}^\epsilon \in \mathbb{R}^{rm \times rm}.
\end{equation}
The regularized mass matrix and its square root are then defined by 
\begin{equation}
\bM^\epsilon = ({\itbf Q}^\epsilon)^T \bLa^\epsilon {\itbf Q}^\epsilon = (\bR^\epsilon)^T \bR^\epsilon \in \mathbb{R}^{rm \times rm}.
\end{equation}   
The regularized wave operator ROM follows similarly \cite[Appendix E]{borcea2022waveform}.}

\section{Waveform inversion}
\label{sect:WFI}
\bc{We give two ideas for estimating $c(\bx)$ from our ROMs. The first one,
described in section \ref{sect:INV_INTW},  uses an estimate of the wave field at points inside the medium. This estimate is based on  the ROM snapshots of the  propagator ROM and is given in section \ref{sect:INTW}. The second idea, described in section \ref{sect:INV_AROM}, formulates the inversion as 
a minimization of the wave operator ROM  misfit. We summarize the two approaches in section \ref{sect:algo} and compare them briefly. A more detailed comparison can be seen in the  numerical simulations in section \ref{sect:Numerics}.}

\subsection{The estimated internal wave}
\label{sect:INTW}
The chain of mappings from the wave speed $c(\bx)$ to the propagator ROM is
\begin{equation}
c(\bx) \mapsto \big\{\bD_j, ~j  =0, \ldots, 2n-1\big\} \mapsto \big\{\cP^{\RM}, ~\bu_j^\RM, ~ j \ge 0\big\}.
\label{eq:ctoROM}
\end{equation}
Both maps in this chain are nonlinear, although we  explained that the second map, from the data to the ROM, 
can be computed with linear algebra tools. The nonlinear steps in the ROM computation are the Cholesky factorization of the 
mass matrix $\bM$, which depends linearly on the data, and the inverse of its Cholesky square root $\bR$. Inverting the chain of mappings 
\eqref{eq:ctoROM} 
to estimate $c(\bx)$ seems as hard as inverting the first map, which is what data fitting does, so how can we use the ROM?

We describe here and in the next section an approach that is inspired by ``hybrid inverse problems"  like photo-acoustic tomography, transient elastography, etc. 
 In these multi-physics imaging modalities, the propagation of a primary wave inside the inaccessible domain, a.k.a. the ``internal wave",  is monitored with high time and space resolution by a second type of wave. The point is that knowledge of the  internal wave simplifies considerably  the inversion for the unknown coefficients of the governing PDE \cite{nachman2007conductivity,bal2013reconstruction,arridge2012imaging}. 
Hybrid approaches are mostly limited to medical applications,
because they involve delicate and accurate user controlled apparatus for transmitting several types of waves and measuring all around the body of interest. We propose to use the ROM for estimating the internal waves, without any additional measurements.

The main idea for the estimation comes from equation \eqref{eq:DF3}, which relates \bc{the data driven ROM snapshots $\bu_j^\RM$, which are matrices,  to the true vector snapshots $\bu_j(\bx) $ defined by (\ref{eq:S1}),  which are $\bx$ dependent fields, for $\bx \in \Omega$}. This relation involves the orthonormal basis 
$\bV(\bx)$, which cannot be computed, but has properties that are useful for inversion. So far, these  properties are only partially understood: It is proved in \cite[Appendix A]{borcea2021reduced} using explicit calculations in a layered medium and in a waveguide that $\bV(\bx)$ is insensitive to the reflectivity (rough part of $c(\bx)$) but depends on the kinematics  (smooth part of $c(\bx)$). Extensive numerical simulations carried out in \cite{DtB,druskin2018nonlinear,borcea2021reduced,borcea2022waveform} suggest that this property extends to more general settings as well. 

In many, but not all applications, the kinematics  is only mildly perturbed, so these results motivated originally our definition of the estimated internal wave snapshots,
\begin{equation}
\bu_j^{\rm est}(\bx;w) = \bV(\bx;w) \bu_j^\RM = \bV(\bx;w) \bR \be_j, \quad j = 0, \ldots, n-1, ~~ \bc{\bx \in \Omega},
\label{eq:EstIW}
\end{equation}
where $\bV(\bx;w)$ is the orthonormal basis computed for the search wave speed $w(\bx)$. Why is this better 
than the linearization of the forward mapping used by any iterative inversion algorithm?  That would approximate the
internal wave snapshots by 
\begin{equation}
\hspace{-0.06in}\bu_j(\bx;w)  = \cos \big[ j \tau \sqrt{\cA(w)} \big] \bu_0(\bx) = \bV(\bx;w) \bR(w) \be_j, \quad j = 0, \ldots, n-1, ~~ \bc{\bx \in \Omega,}
\label{eq:Born}
\end{equation}
where $\cA(w) = - w(x) \Delta \big[ w(x) \cdot \big]$ and $\bR(w)$ is the Cholesky square root of the mass matrix 
computed from the simulated data at wave speed $w(\bx)$. The advantage of \eqref{eq:EstIW} over \eqref{eq:Born} is 
that it is consistent with the measurements, as stated  next.

\vspace{0.05in}
\begin{cor}
\label{lem.3}
The estimated internal wave snapshots \eqref{eq:EstIW} fit  the data 
\begin{align*}
\bD_j &= \int_{\Omega} d \bx \, \bu_0(\bx)^T \bu_j(\bx) = \int_{\Omega} d \bx \, \bu_0^{\rm est} (\bx)^T \bu_j^{\rm est}(\bx), 
\\
\bD_{j+n-1} &=  \int_{\Omega} d \bx \,  [2 \bu_{n-1}^T(\bx) \bu_j(\bx)  - \bu_0^T(\bx) \bu_{n-1-j}(\bx)] \\ {}&= \int_{\Omega} d \bx \, [2 \bu_{n-1}^{\rm est}(\bx)^T \bu_j^{\rm est}(\bx) - \bu_0^{\rm est}(\bx)^T \bu_{n-1-j}^{\rm est}(\bx)],
\end{align*}
for $j = 0, \ldots, n-1$, whereas $\bu_j(\bx;w)$ defined in \eqref{eq:Born} do not.
\end{cor}

\vspace{0.05in}
\textbf{Proof:} The proof follows easily from that of Theorem \ref{lem.2}, because equation \eqref{eq:DF4} holds if we replace 
$\bV(\bx)$ by any orthonormal basis, like $\bV(\bx;w)$. The data fit is ensured by the ROM snapshots $\bu_j^\RM = \bR \be_j$,
for $j = 0, \ldots, n-1$, irrespective of the basis. Obviously, if we replace $\bR$ by $\bR(w)$, we no longer fit the data, 
so the snapshots \eqref{eq:Born} are not consistent with the measurements. $~~\Box$

\begin{figure}[t]
\centering
\includegraphics[width=0.38\textwidth]{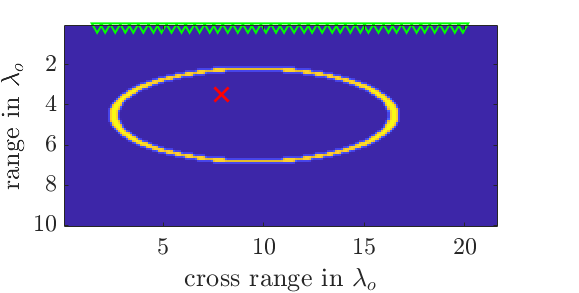}
\includegraphics[width=0.38\textwidth]{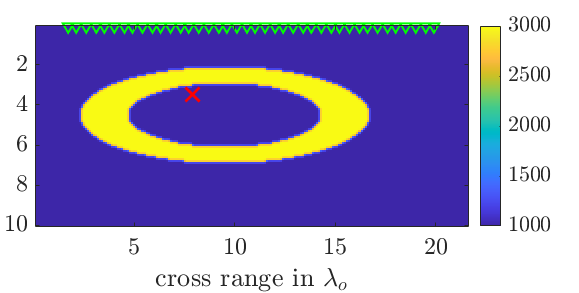}\\
\vspace{0.075in}
\includegraphics[width=0.189\textwidth]{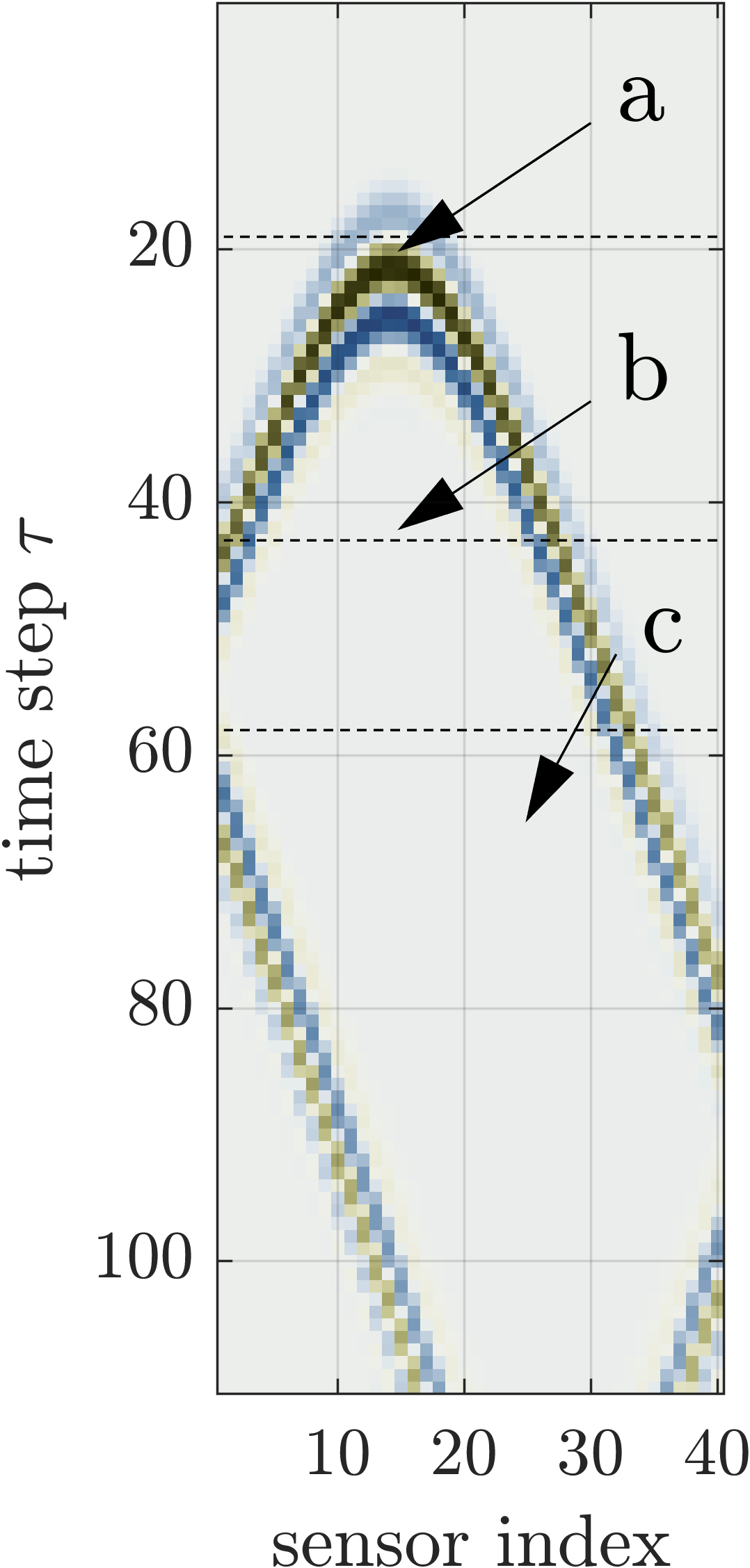}
\hspace{-0.05in} 
\includegraphics[width=0.165\textwidth]{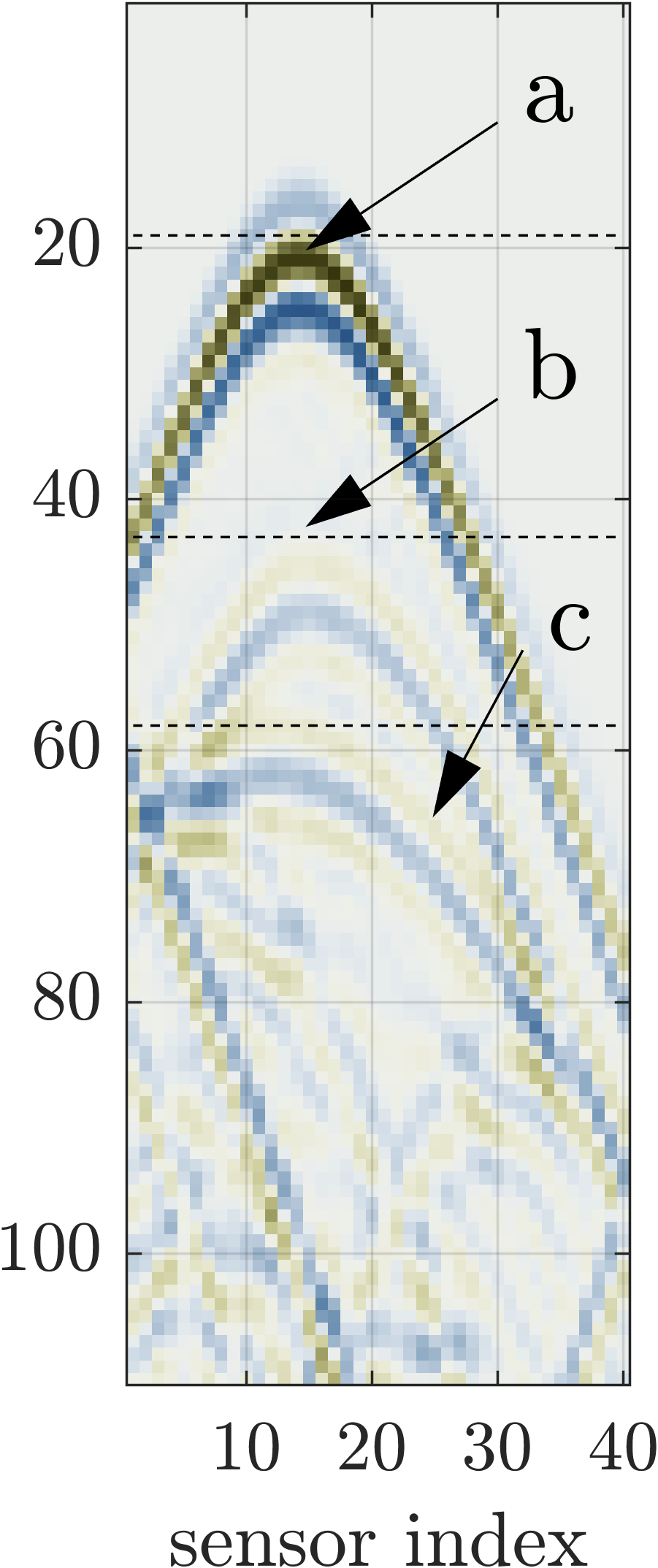}
\hspace{-0.05in} 
\includegraphics[width=0.165\textwidth]{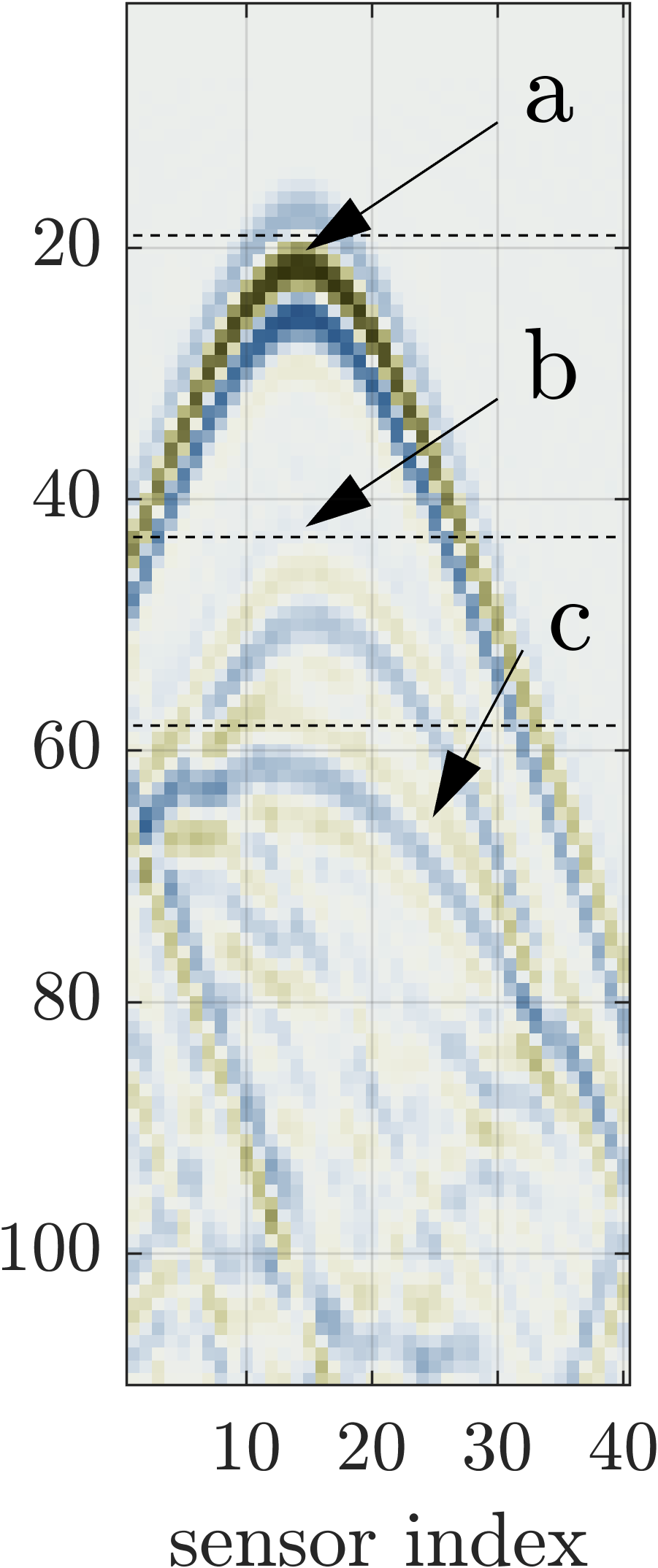}
\hspace{-0.05in} 
\includegraphics[width=0.165\textwidth]{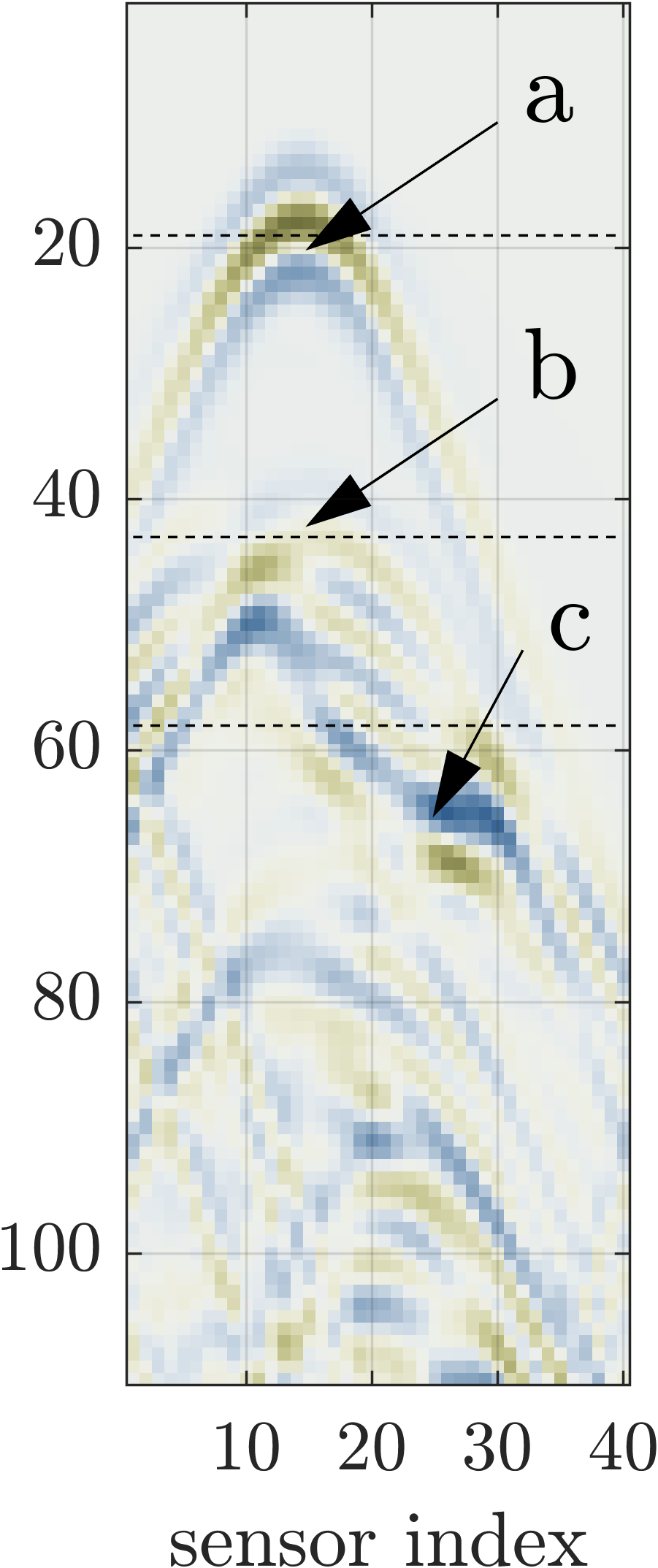}
\hspace{-0.05in} 
\includegraphics[width=0.165\textwidth]{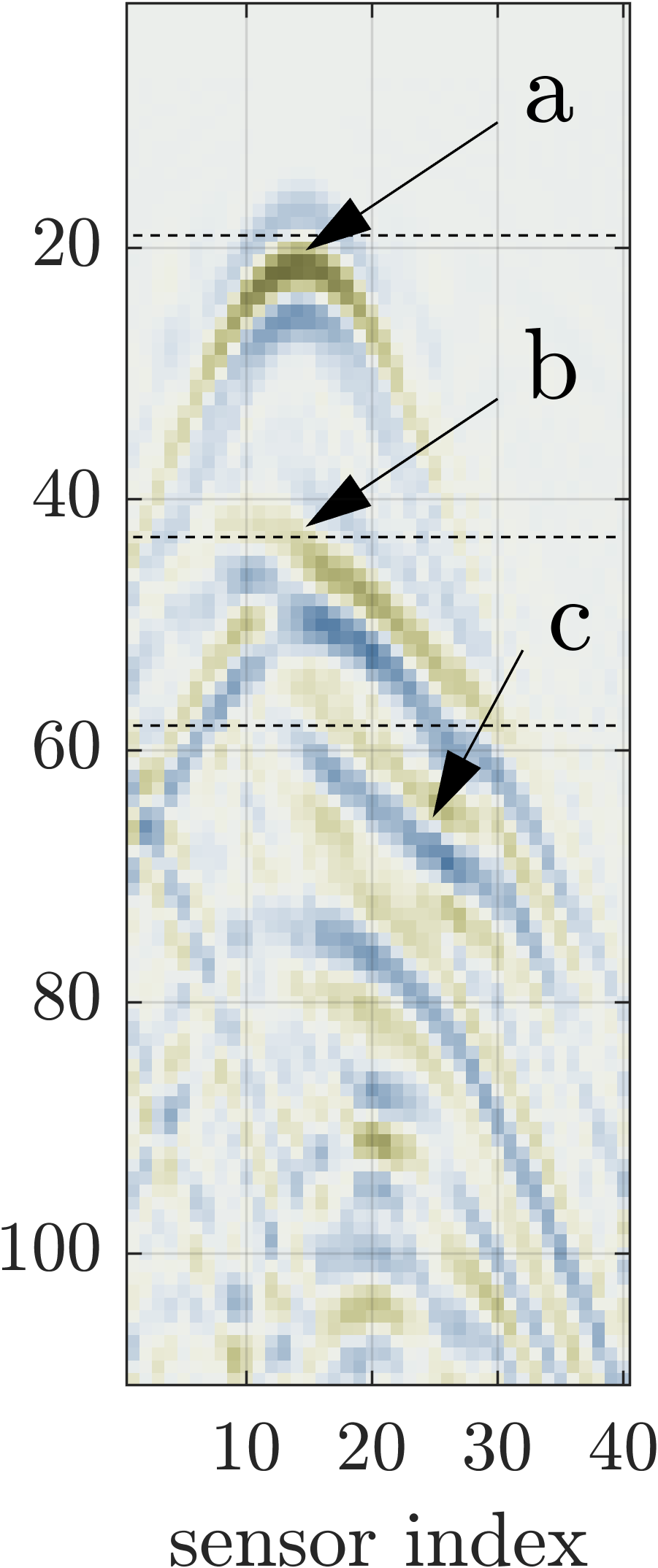}
    	\caption{Illustration of the internal waves at  $\bx$ indicated with \bc{red $\times$} in two media containing a ring shaped inclusion (top plots). The active array is linear and placed at the top of the plots, with sensors indicated by green triangles. The wave speed is shown in the colorbar in units of m/s. 
	The axes are cross-range (measured along the array) and range (measured orthogonal to the array), in units   of the central wavelength.  The bottom plots show the internal waves as  functions of the sensor index $s$ (abscissa) 
	and time in units of $\tau$ (ordinate). From left to right: $u^{(s)}(j\tau,\bx;w=\bar c)$ calculated with the constant wave speed $ \bar{c} = 1$km/s; the true wave $u^{(s)}(j\tau,\bx)$ in the first medium, the estimated $u^{(s),{\rm est}}(j\tau,\bx;w= \bar c)$ in the first medium; 
the true wave $u^{(s)}(j\tau,\bx)$ in the second medium and the estimated $u^{(s),{\rm est}}(j\tau,\bx;w= \bar c)$ in the second medium.	}
     	\label{fig:InternalSolInc}
\end{figure}

\vspace{0.05in}
To illustrate this result, we display  in Fig. \ref{fig:InternalSolInc} the internal waves at a point $\bx$ inside two media with a fast ring shape inclusion with wave speed $3$km/s, in a background of wave speed $\bar{c} = 1$km/s. The difference between the media is the thickness of the ring: In the first medium the ring is thin, so the kinematics is only slightly perturbed. In the second medium the ring is thick, so the perturbation of the kinematics is more significant. The true internal waves are in the second and fourth bottom plots. They display three types of arrivals: (1) The direct arrival (the hyperbola denoted by $a$) of the wave that travels from the array to $\bx$, through the inclusion.  (2) The echoes that have bounced between the inclusion and the top boundary, like those denoted by $b$ and $c$ in the plot. (3) The echoes from the side boundaries of the domain $\Omega$. The bottom left plot shows the wave \eqref{eq:Born} calculated with   $w \equiv \bar c$. It is quite different from the true waves, as it contains only the direct arrival and the echoes from the side boundaries. The third and fifth bottom plots show our estimate \eqref{eq:EstIW} of the internal wave in the two media, for the same $w \equiv \bar c$. They are quite similar to the true waves and 
contain all the arrivals. This is because the information about these arrivals is in  $\bR$, which is responsible for  the data fit stated in Corollary \ref{lem.3}. The purpose of $\bV(\bx;w)$ in 
\eqref{eq:EstIW} is to map these arrivals from the algebraic (ROM) space to the physical space. Because we do not use the right kinematics and this significantly affects $\bV(\bx;w = \bar c)$ in the case of the thick ring, we can see that the arrivals are delayed in the bottom fifth plot compared to the true internal wave in the bottom fourth plot.

\subsection{Inversion with the estimated internal wave}
\label{sect:INV_INTW}
It is well known that the forward map, which relates $c(\bx)$ to the wave measured at the array, is given by the Lippmann-Schwinger integral equation. Due to the transformations described in section \ref{sect:ROM.1}, this equation has the following form derived in \cite[Proposition 4.1]{borcea2022internal} 
\begin{equation}
D_{r,s}(t)-D_{r,s}(t;w) = \int_0^t dt' \int_\Omega d \bx \, \frac{[c^2(\bx)-w^2(\bx)]}{c(\bx) w(\bx)} u^{(s)}(t',\bx) \partial_{t'} u^{(r)}(t',\bx;w),
\label{eq:LS}
\end{equation}
\bc{where we recall that in our notation convention we add $w$ to the arguments of the fields computed with the search speed $w(\bx)$ and suppress $c(\bx)$ in the arguments of the fields corresponding to the true speed.}

Equation \eqref{eq:LS} shows that the mapping $c(\bx) \mapsto \bD(t)$ is nonlinear not only because of the way $c(\bx)$ enters the first factor of the integrand, but also because the internal wave $u^{(s)}(t',\bx)$ depends on $c(\bx)$ in a very complicated way. The first inversion approach introduced in \cite{borcea2022internal} is to  replace this wave by the estimate \eqref{eq:EstIW}. 
The iterative optimization based on the resulting mapping works better than that of FWI, as long as the kinematics is not too perturbed \cite[Section 5]{borcea2022internal}. Otherwise, the iteration can get stuck, like that of FWI.

Since the improvement brought by the estimated internal wave in data fitting is marginal, is there a better way to use it? The new 
idea that we now describe takes the following point of view: Corollary \ref{lem.3} states that the estimated internal wave is consistent 
with the measurements by construction, so why should we seek to fit the data at all? What the estimate \eqref{eq:EstIW} is not guaranteed to do is to be an approximate solution of the wave equation. If that were the case, then $w(\bx)$ would be close to  the solution $c(\bx)$ of the  inverse problem\footnote{This is assuming that $c(\bx)$ can be determined uniquely for the given measurement setup and the inversion is stable due to  proper regularization and parametrization of the search speed.}.

How can we enforce that $\bu_j^{\rm est}(\bx;w)$ be close to the solution $\cos \big[ j \tau \sqrt{\cA(w)} \, \big] \bu_0(\bx)$ of the wave equation for $j = 0, \ldots, n-1?$ First, we note that $\bu_0(\bx)$ is known because it depends on the medium near the array, where the wave speed equals the given constant $\bar c$.  Second, we obtain from  definitions \eqref{eq:EstIW} and \eqref{eq:Born} that the quadratic misfit between the estimated internal wave $\bu_j^{\rm est}(\bx;w)$ and the solution $\cos \big[ j \tau \sqrt{\cA(w)} \, \big] \bu_0(\bx)$ is
\begin{align}
\int_{\Omega} d \bx \sum_{j=0}^{n-1}\big\| \bu_j^{\rm est}(\bx;w) - \cos \big[ j \tau \sqrt{\cA(w)} \, \big] \bu_0(\bx) \big\|_{\rm F}^2 
&=  \int_{\Omega} d \bx \big\| \bV(\bx,w) \big[\bR - \bR(w)\big] \big\|_{\rm F}^2 \nonumber \\
&= \|\bR - \bR(w)\|_{\rm F}^2,
\label{eq:fitR}
\end{align}
because $ \int_{\Omega} d \bx  \bV^T(\bx,w)\bV(\bx,w)= \bI_{nm}$ for any $w$.
Thus, we can just minimize 
the misfit of the Cholesky square roots of the data driven mass matrices. 

We consider a modification of this objective function, given by 
\begin{equation}
\mathcal{O}(w) = \| \bR(w)^{-1} \bR - \bI_{nm} \|_{\rm F}^2.
\label{eq:objR}
\end{equation}
This is motivated by the following observation: Variations in the medium that are closer to the array give much stronger echoes than 
those that are further away. The information about the latter is encoded in the smaller entries of $\bM$ and consequently, the smaller 
singular values of $\bR$. To balance these contributions, we emphasize the weaker events by taking the inverse 
of $\bR(w)$ in \eqref{eq:objR}. \bc{Here we ignore that $\bR(w)$ may be ill-conditioned, but this can be addressed by regularization as explained in 
section \ref{sect:regul} and in \cite[Appendix E]{borcea2022waveform}. }

\subsection{Inversion with the ROM wave operator}
\label{sect:INV_AROM}
The second objective function for ROM based waveform inversion was introduced recently in \cite{borcea2022waveform}. It measures the 
wave operator ROM misfit 
\begin{equation}
\tilde{\mathcal{O}}(w) = \| \bAR - \bAR(w) \|_{\rm F}^2,
\label{eq:objA}
\end{equation}
where $\bAR $ is obtained from the data by (\ref{eq:G28}) using the matrices $\bR$ and $\tilde\bS$ built in Theorems~\ref{thm.1}-\ref{thm.2}. 
The matrix
$\bAR(w)$ is obtained from the synthetic data simulated with the wave speed $w(\bx)$.

Note that in (\ref{eq:objA}) we do not penalize the data misfit and, as explained in section \ref{sect:ROM.5}, the ROM wave $\tilde \bu^{\RM}(t)$ does not reproduce exactly the measurements. Thus, unlike the approach described in the previous section, based on the estimated internal wave, we cannot expect  here an automatic data fit. Nevertheless, the data are fit approximately, because the ROM comes from a Galerkin approximation of the wave equation.
Our motivation for using the objective function \eqref{eq:objA} is threefold: 

\vspace{0.05in} \noindent 1. The operator $\cA$ depends in a simple (quadratic) way on the unknown $c(\bx)$, unlike the propagator operator $\cP$. Indeed,
fitting the propagator is not a good idea, as illustrated in section \ref{sect:Numerics}.

\vspace{0.05in}  \noindent  2. The projection basis $\bV(\bx)$ depends on $c(\bx)$, which makes the analysis of \eqref{eq:objA} very complicated. However,
the existing theoretical results  and numerous numerical simulations suggest that the change of this basis with $c(\bx)$ is slower than that of $\cA$. In other words, we expect that $\bAR$ is approximately quadratic in the wave speed $c(\bx)$. This is obviously 
not the case for the data matrices used in the FWI objective function.

\vspace{0.05in}  \noindent 3.  Numerical evidence \cite{borcea2022waveform,DtB} suggests that the $m$ components of the $\bv_j(\bx)$ block of $\bV(\bx)$ are localized around the maximum depth reached by the wave field at time $j \tau$. Thus, definition \eqref{eq:G27} suggests that the entries in 
$\bAR$ depend mostly on local averages of $c(\bx)$, as is usual in standard Galerkin approximation schemes.

\subsection{\bc{The inversion algorithms. Pros and cons}}
\label{sect:algo}
\bc{Here we give the pseudo-algorithms that summarize the steps in our two  inversion approaches and then, compare them. Both algorithms require a 
parametrization of the search velocity space, 
\begin{equation}
w(\bx) = \bar{c} + \sum_{j=1}^N \eta_j \phi_j(\bx), \qquad \boldsymbol{\eta} = (\eta_1, \ldots, \eta_N) \in \RR^N, 
\label{eq:searchw}
\end{equation}
using some user defined basis functions $\{\phi_j(\bx)\}_{j=1}^N$, 
for $\bx \in \Omega$.}

\vspace{0.05in}
\begin{algorithm}\bc{\textbf{\emph{(Inversion with the internal wave)}}}
\label{AlgR}

\vspace{0.04in} \noindent \bc{\textbf{Input:} Data matrices $\{\bD_j\}_{j=0}^{2n-1}$ calculated from the measurements as 
in Lemma \ref{lem.1} and the reference speed $\bar{c}$.}

 \vspace{0.04in} \noindent \bc{1. Compute $\bM$ with block entries given in equation \eqref{eq:G13}. 
  If $\bM$ is indefinite or singular, use regularization as described in section 
\ref{sect:regul}.} 

 \vspace{0.04in} \noindent \bc{2.  Compute the block Cholesky square root $\bR$ of $\bM$ or its regularized version. }
 
 \vspace{0.04in} \noindent \bc{3.  Starting with the initial vector $\boldsymbol{\eta} = \boldsymbol{\eta}^{(0)}= {\bf 0}$ in equation \eqref{eq:searchw}, proceed:}
 
\vspace{-0.1in} \bc{
 \begin{itemize}
 \itemsep 0.03in
\item  For update index $j \ge 1$ calculate $w(\bx)$ as in \eqref{eq:searchw}, with $\boldsymbol{\eta} = \boldsymbol{\eta}^{(j-1)}.$
\item Calculate $\bR(w)$ following the same procedure of calculating $\bR$. 
\item Compute $\boldsymbol{\eta}^{(j)}$ as a Gauss-Newton update for minimizing the objective function \eqref{eq:objR} with 
the user's choice of the regularization penalty on $\boldsymbol{\eta}$. 
\item Go to the next iteration or stop when the user defined convergence criterion has been met. 
\end{itemize}
}

\vspace{0.04in} \noindent \bc{\textbf{Output:}  The estimate of the wave speed given by \eqref{eq:searchw} with $\boldsymbol{\eta}$ calculated at step 3.}
\end{algorithm}

\vspace{0.1in}
\begin{algorithm}\bc{\textbf{\emph{(Inversion with the wave operator ROM)}}}
\label{AlgA}

\vspace{0.04in} \noindent \bc{\textbf{Input:} Data matrices $\{\bD_j\}_{j=0}^{2n-1}$ calculated from the measurements as 
in Lemma \ref{lem.1} and the reference speed $\bar{c}$.}

 \vspace{0.04in} \noindent \bc{1. Compute $\{\ddot \bD_j\}_{j=0}^{2n-2}$ using, for example, the Fourier transform.}

 \vspace{0.04in} \noindent \bc{2. Compute $\bM$ with the block entries given in equation \eqref{eq:G13} and 
 $\tilde \bS$ with the block entries \eqref{eq:G15}. If needed, regularized matrices can be computed as explained in 
 \cite[Appendix E]{borcea2022waveform}.}

 \vspace{0.04in} \noindent \bc{3. Compute the block Cholesky square root of $\bM$ or its regularized version and 
 the wave operator ROM $\bAR$ using the right hand side in \eqref{eq:G28}.}
 
 \vspace{0.04in} \noindent \bc{4.  Starting with the initial vector $\boldsymbol{\eta} = \boldsymbol{\eta}^{(0)}= {\bf 0}$ in equation \eqref{eq:searchw}, proceed:}
 
\vspace{-0.1in} \bc{
 \begin{itemize}
 \itemsep 0.03in
\item  For update index $j \ge 1$ calculate $w(\bx)$ as in \eqref{eq:searchw}, with $\boldsymbol{\eta} = \boldsymbol{\eta}^{(j-1)}.$
\item Calculate $\bAR(w)$ following the same procedure of calculating $\bAR$. 
\item Compute $\boldsymbol{\eta}^{(j)}$ as a Gauss-Newton update for minimizing the objective function \eqref{eq:objA} with 
the user's choice of the regularization penalty on $\boldsymbol{\eta}$. 
\item Go to the next iteration or stop when the user defined convergence criterion has been met. 
\end{itemize}
}
\vspace{0.04in} \noindent \bc{\textbf{Output:}  The estimate of the wave speed given by \eqref{eq:searchw} with $\boldsymbol{\eta}$ calculated at step 4.}
\end{algorithm}

\vspace{0.05in}
\bc{When we compare the two approaches, we note that Algorithm \ref{AlgR} is easier to use, as it involves only the mass matrix and its block-Cholesky square root. Moreover, if the objective function \eqref{eq:objR} is small, then the data are fit implicitly
by Corollary~\ref{lem.3}. Algorithm \ref{AlgA} is slightly more expensive computationally, because it requires the computation 
of the second derivatives of the data matrices and the wave operator ROM. Moreover, the data are not fit as well as in the 
first approach. Nevertheless, the numerical results show that both approaches give comparable estimates of the wave speed. In fact,  in a few numerical experiments (see Fig. \ref{fig:Camembert_Inv}), we see a slightly better result with Algorithm \ref{AlgA}.}

\subsection{Numerical results and comparison of the inversion approaches}
\label{sect:Numerics}
We refer to appendix \ref{ap:A} for the details of our numerical simulations. Here we compare the inversion approaches by
visualizing the objective functions in a two-dimensional search space (section \ref{sect:OBJ}) and by showing the inversion results for three well-known challenging media (section \ref{sect:Inv}). \bc{We discuss  in section \ref{sect:compnum} 
the computational cost and compare it to that of FWI.} 
\subsubsection{Visualization of the objective functions}
\label{sect:OBJ}
\begin{figure}[t]
     \centering
     \begin{subfigure}[b]{0.277\textwidth}
         \centering
         \includegraphics[width=\textwidth]{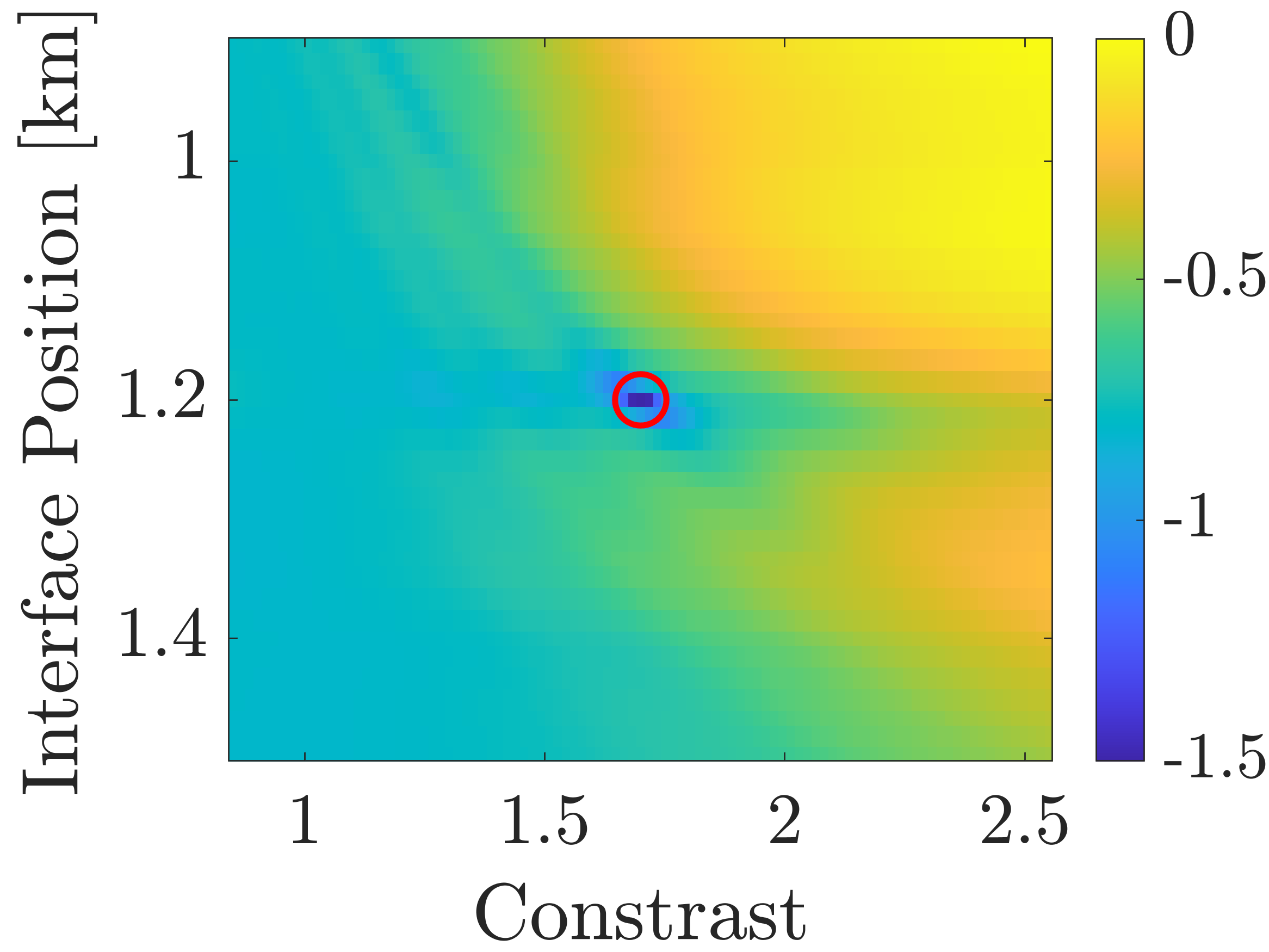}
         \caption{$\log {\cal O}(w)$}
         \label{fig:LayerObjRinvR}
     \end{subfigure}
	\hfil
     \begin{subfigure}[b]{0.23\textwidth}
         \centering
         \includegraphics[width=\textwidth]{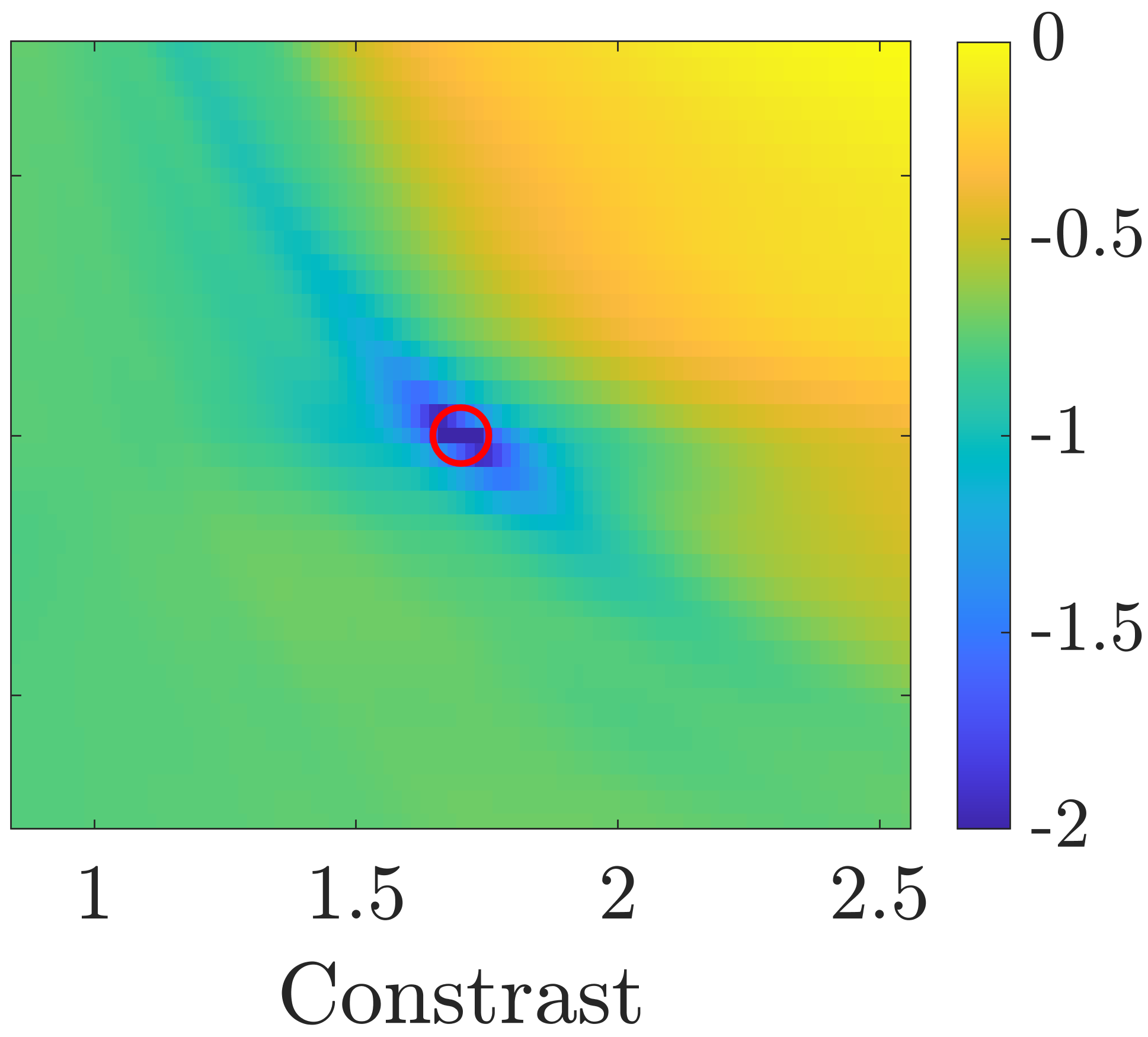}
         \caption{$\log \tilde{\cal O}(w)$}
         \label{fig:LayerObjROM}
     \end{subfigure}
     	\hfil
     \begin{subfigure}[b]{0.23\textwidth}
         \centering
         \includegraphics[width=\textwidth]{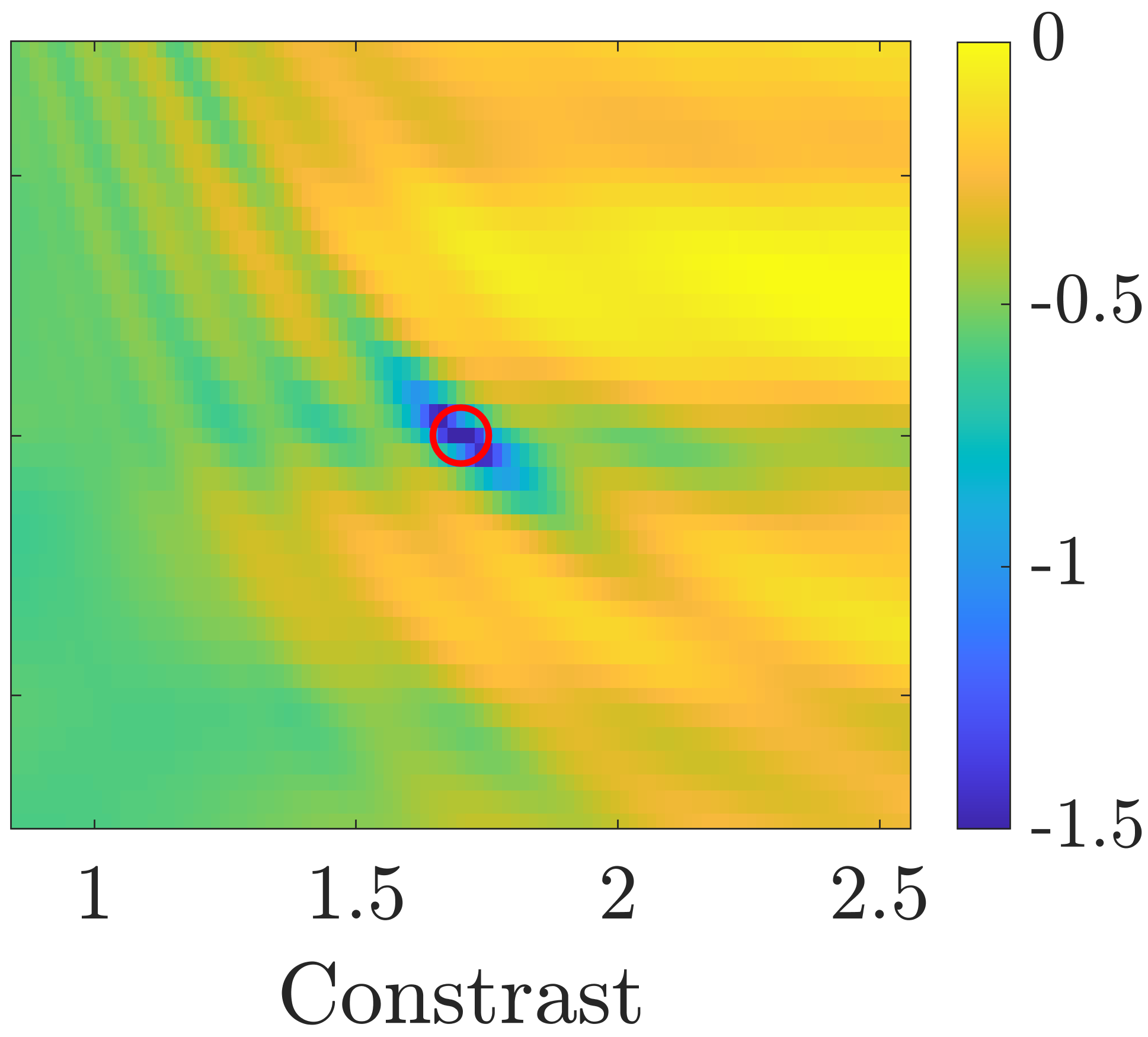}
         \caption{$\cP^\RM$ misfit}
         \label{fig:LayerObjPro}
     \end{subfigure}
	\hfil
         \begin{subfigure}[b]{0.23\textwidth}
         \centering
         \includegraphics[width=\textwidth]{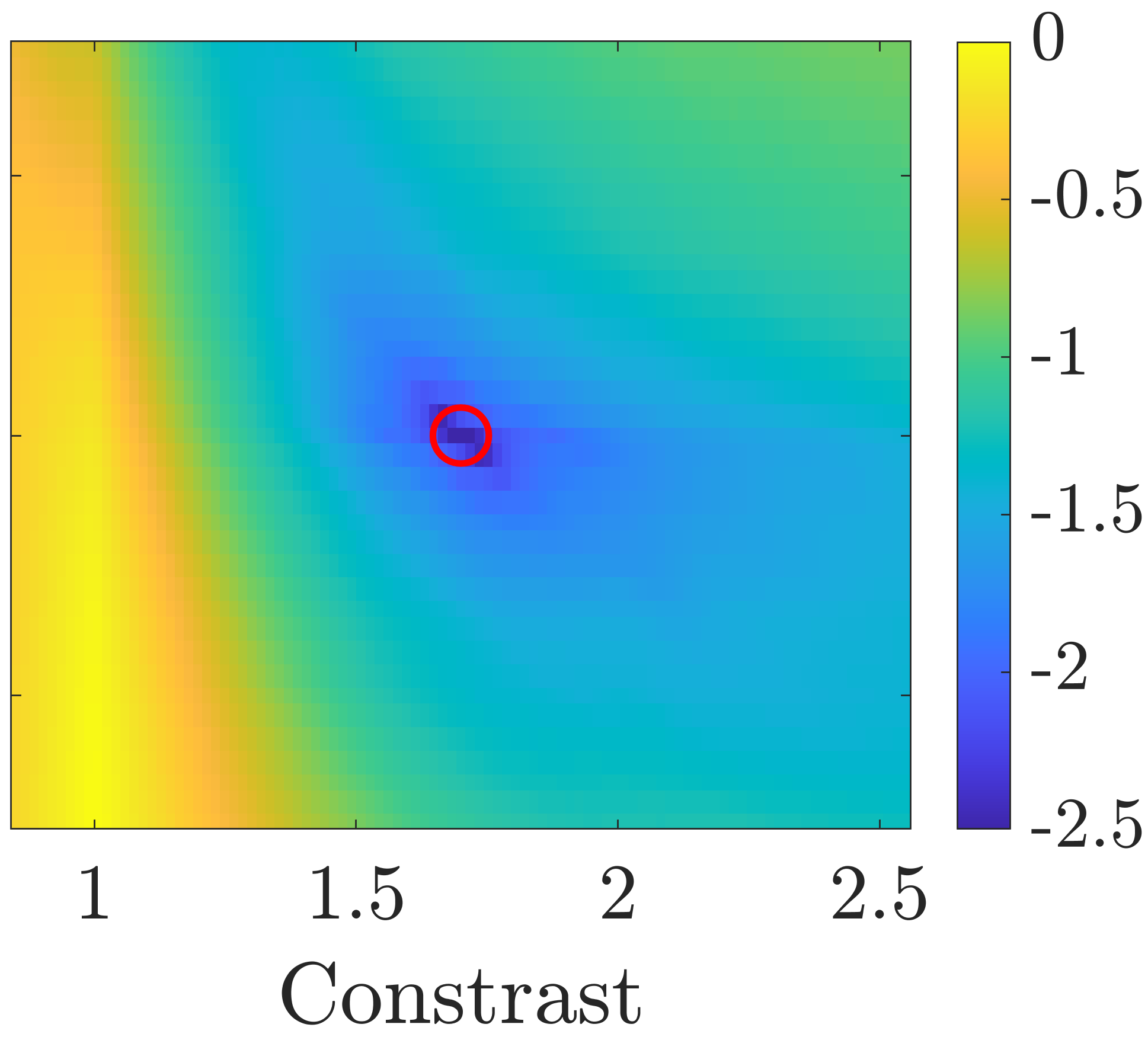}
         \caption{Wasserstein}       
     \end{subfigure}
     \vspace{-0.2in}
     \caption{Display of the log of the objective functions for the setup in Fig. \ref{fig:FWI_topo}. The axes are the same as in the right plot there. The true model is indicated by $ \textcolor{magenta}{\bigcirc}$.}
 \label{fig:LayerObjWas}
\end{figure}

One way to compare the various inversion approaches is to visualize their objective functions using  a two-dimensional search space. We
begin in Fig. \ref{fig:LayerObjWas} with the objective functions for the setup in Fig. \ref{fig:FWI_topo}, where the multiple minima in the FWI objective function were clearly visible. We note that both ROM objective functions \eqref{eq:objR} (plot (a)) and \eqref{eq:objA} (plot (b)) have a single minimum at the correct location. In plot (c) we display the log of the propagator ROM misfit 
$\| \cP^{\RM}- \cP^{\RM} (w) \|_{\rm F}^2$  to show that it is 
not a good choice for the inversion. This is not surprising, because the sought after $c(\bx)$ appears in a very complicated way in the 
expression of the propagator operator defined in \eqref{eq:defProp}. The plot (d) shows  that if instead of quantifying the data misfit in the  $L^2([t_{\rm min},t_{\rm max}])$ norm, as in FWI, we  use the Wasserstein
metric, computed as described in \cite{EngquistFroese}, we get a much better objective function. Indeed, the plot shows a single minimum at the correct location.

\begin{figure}[h]
\centering
\begin{subfigure}[b]{0.41\textwidth}
\includegraphics[width=\textwidth]{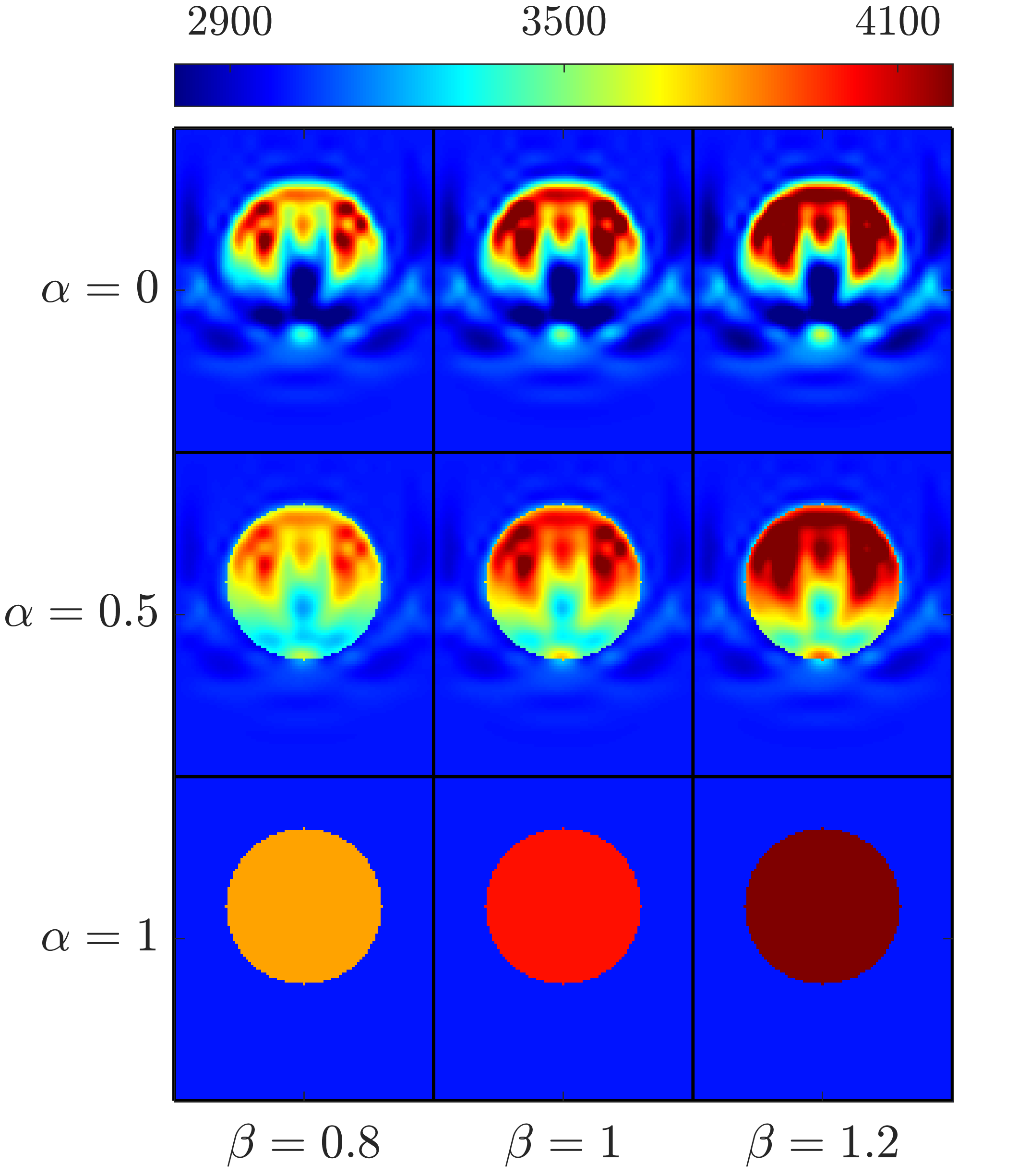}
\label{fig:Cheesesearch}
\end{subfigure}
\hspace{-0.1in}
\raisebox{1.3in}{\begin{subfigure}[b]{0.55\textwidth}
\begin{tabular}{cc}
$(a) \log \mathcal{O}^\FWI(w)$ & (b) Wasserstein\\
\includegraphics[width=0.51\textwidth]{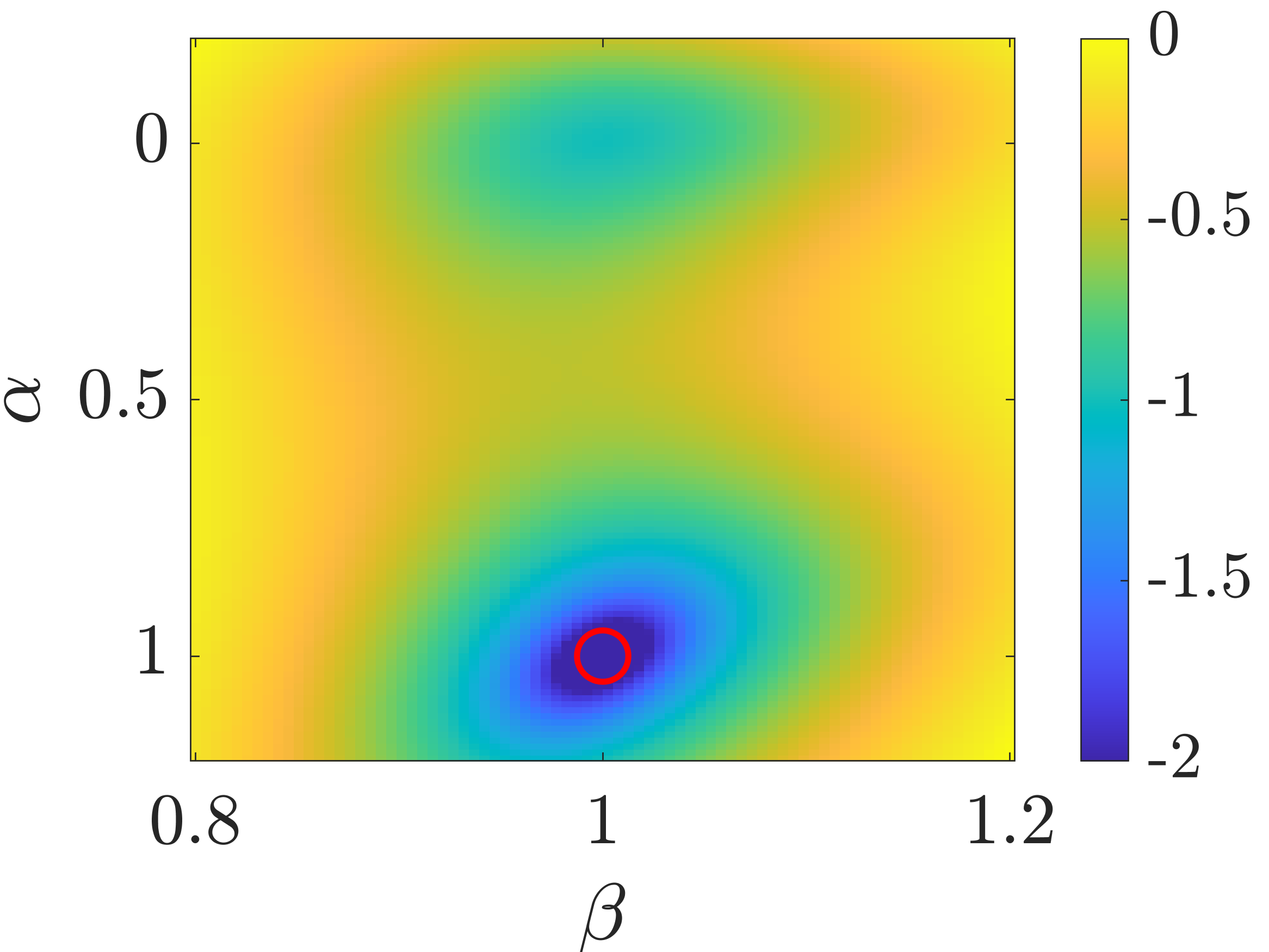} & \includegraphics[width=0.51\textwidth]{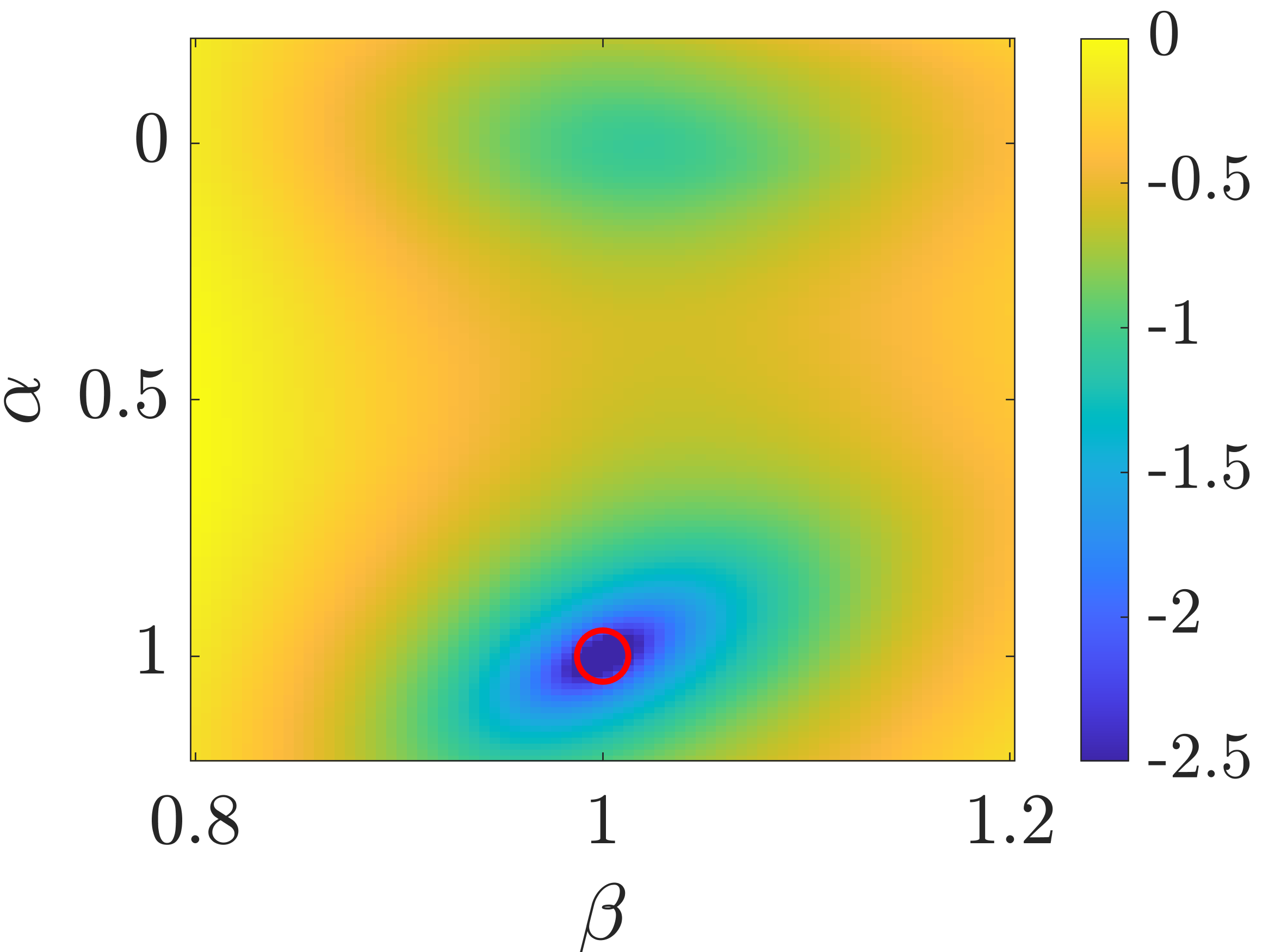} \\
$(c) \log \mathcal{O}(w)$ & (d) $\log \tilde{\mathcal{O}}(w)$\\
\includegraphics[width=0.51\textwidth]{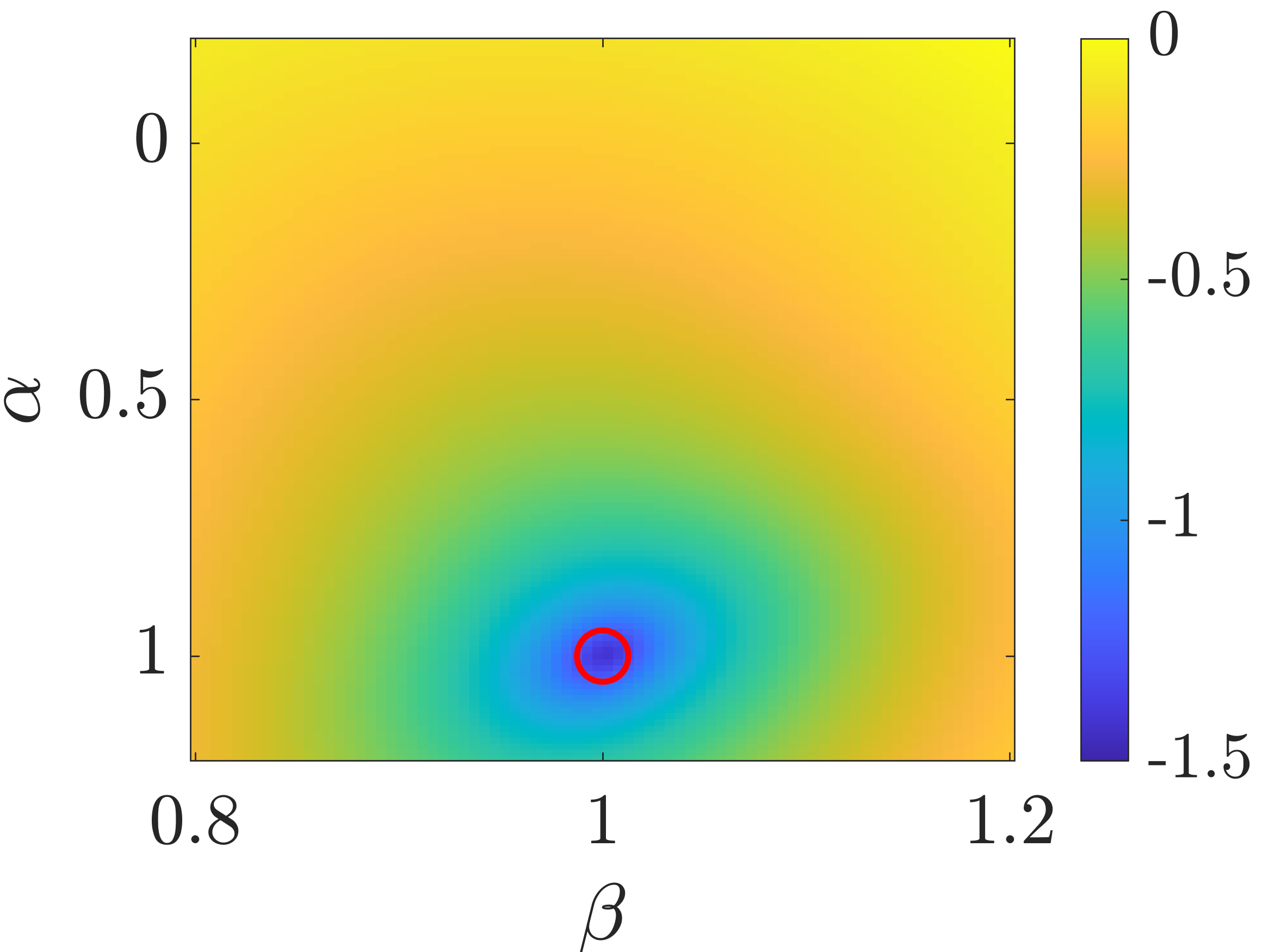} & \includegraphics[width=0.51\textwidth]{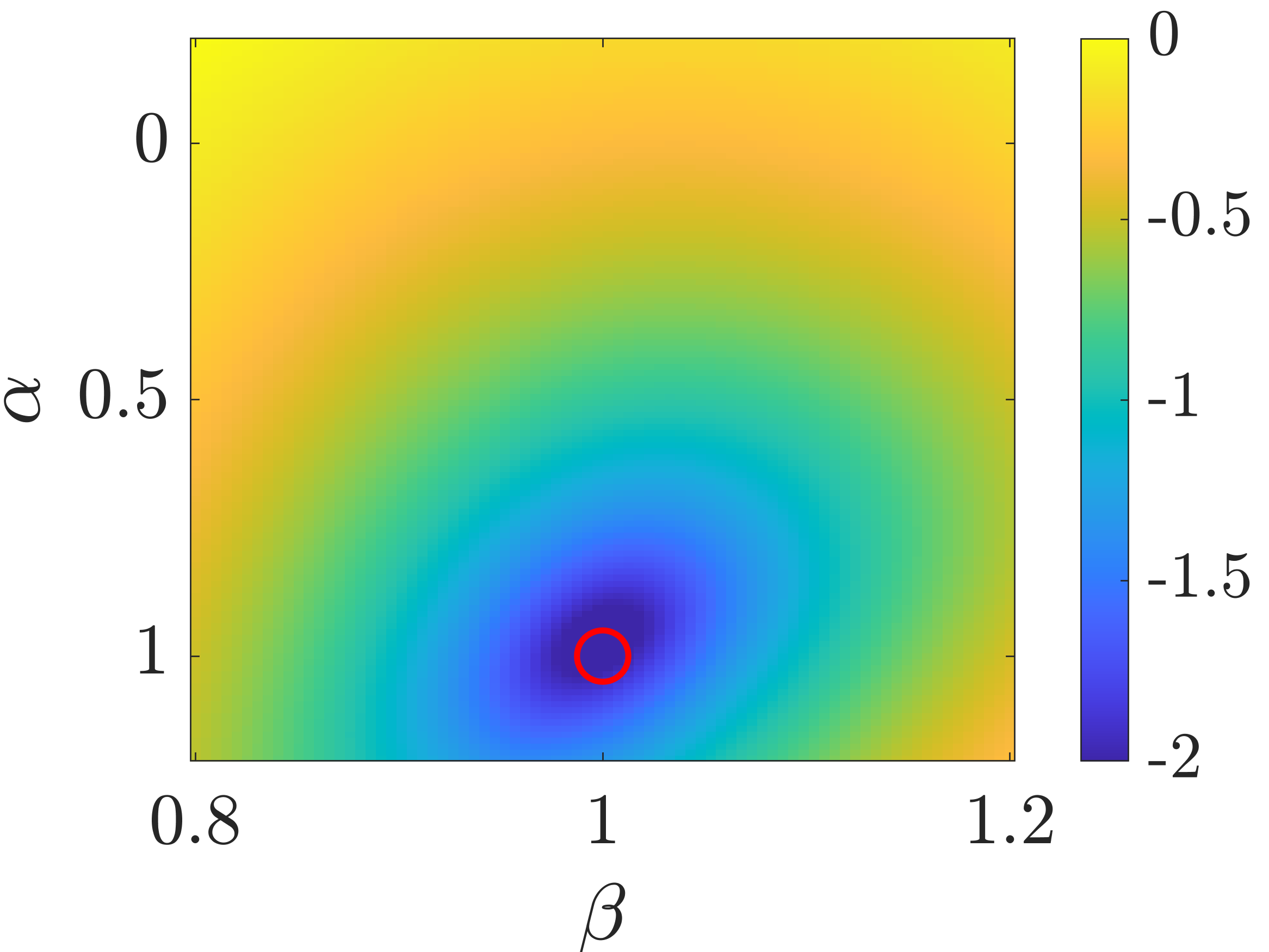}
\end{tabular}
\end{subfigure}}
\vspace{-0.1in}
 \caption{Left plot: Illustration of the search speed \eqref{eq:CamembertSearch} for $9$ choices of $\alpha$ and $\beta$. Right plots: Log of objective functions parametrized by $\beta$ (abscissa) and 
        $\alpha$ (ordinate). The true model is indicated by $\textcolor{magenta}{\bigcirc}$. }
        \label{fig:Cheese}
\end{figure}

To show the advantage of using the ROM based objective functions, we consider next a challenging example, coined as 
the ``Camembert" model in the geophysics community \cite{gauthier1986two}. It consists of a piecewise constant $c(\bx)$, 
modeling a disk shaped fast inclusion embedded in a homogeneous medium with wave speed $\bar{c} = 3$km/s (left plot in Fig. \ref{fig:Camembert_Inv}).  The name is because FWI, starting from the initial guess $w(\bx) = \bar{c}$, can recover well only the top of the inclusion, and the estimate $c^\FWI(\bx)$ of $c(\bx)$ ``melts away" from there, as would Camembert cheese. The details of the inversion are in the next section and appendix \ref{ap:A1}. Here we display  in 
Fig. \ref{fig:Cheese} the objective functions, in a two-dimensional search space with wave speed 
\begin{equation}
w(\bx;\alpha,\beta) =
\bar{c} +  \beta \left\{ (1-\alpha) \left[ c^\FWI(\bx) - \bar{c} \right] + \alpha \left[ c(\bx) - \bar{c} \right] \right\}
\label{eq:CamembertSearch}
\end{equation}
parametrized by  $\alpha \in [-0.21,0.2]$ and $\beta\in [0.8, 1.2]$. Note how $\alpha$ interpolates between the FWI estimate and the true wave speed. The parameter $\beta$ is used to vary the contrast.  Fig. \ref{fig:Cheese} shows that the objective functions based on data fitting (plot (a) for $L^2$-misfit and plot (b) for Wasserstein misfit) have two minima. One of them is at the true speed $c(\bx)$ i.e., $\alpha = 1$ and $\beta = 1$,  and the other one is at the incorrect $c^\FWI(\bx)$ i.e., $\alpha = 0$ and $\beta = 1$. The ROM objective functions (plots (c) and (d)) have a single minimum at $\alpha = 1$ and $\beta = 1$.

\subsubsection{Inversion results}
\label{sect:Inv}
We present  results for three known challenging media. All the results are in two dimensions, for a linear array of co-located sources and receivers. Carrying out the inversion in three dimensions involves no conceptual change, but naturally, the computational cost is higher. This is discussed in \bc{section \ref{sect:compnum}} and \cite{borcea2022waveform}.

The minimization of the objective functions is carried out with the Gauss-Newton method,  in an $N$ dimensional search space $\mathcal{W}$ with $w(\bx)$ of the form \eqref{eq:searchw}, spanned by $N$ Gaussian basis functions $\phi_l(\bx) = \frac{1}{2 \pi \sigma^\perp \sigma} \exp \left[ - \frac{(x^\perp-x_l^\perp)^2}{2 (\sigma^\perp)^2} - 
\frac{(x-x_l)^2}{2 \sigma^2}\right]$ with peaks at the points $\bx_l = (x^\perp_l,x_l)$ on a uniform grid discretizing the 
inversion domain $\Omega_{\rm inv} \subset \Omega$. Here  $x^\perp_l$ denotes the cross-range coordinate (along the array) and $x_l$ is the range coordinate,  in the direction orthogonal to the array. The standard deviations of the Gaussian functions in these directions are $\sigma^\perp$ and $\sigma$.

We do not discuss noise effects on the inversion, but these have been addressed in 
\cite{borcea2022waveform,borcea2022internal}. 
Even in the absence of noise,  diffraction limits the resolution of the inversion to about half the central wavelength.  In our simulations we over parametrize $\mathcal{W}$,  so we add the Tikhonov regularization penalty $\mu \|\bbeta\|_2^2$ to the objective functions, with $\mu$  chosen adaptively during the iteration, as explained in \cite[Appendix D]{borcea2022waveform}.

The inversion is carried out in a layered peeling fashion, by time windowing the data, as explained in section \ref{sect:ROM.5}. If we use $N_t$ 
windows, then each one is of the form $t_{\rm min} \le t \le t_{q,\rm max}$, where $t_{q,\rm max} = t_{\rm min} + \frac{q}{N_t}[t_{\rm max} - t_{\rm min}]$ and  $q = 1, \ldots, N_t$. The updates of the wave speed for the $q^{\rm th}$ window are computed up to the maximum depth sensed by the waves  at $t \le t_{q,{\rm max}}$. This depth can be  estimated in practice using some conservative upper bound on $c(\bx)$. \bc{In our experience, the results are not very sensitive to this depth.}
The parameters for each inversion result displayed below are given in Appendix \ref{ap:A1}. 

\begin{figure}[t]
\begin{center}
\includegraphics[height=5cm]{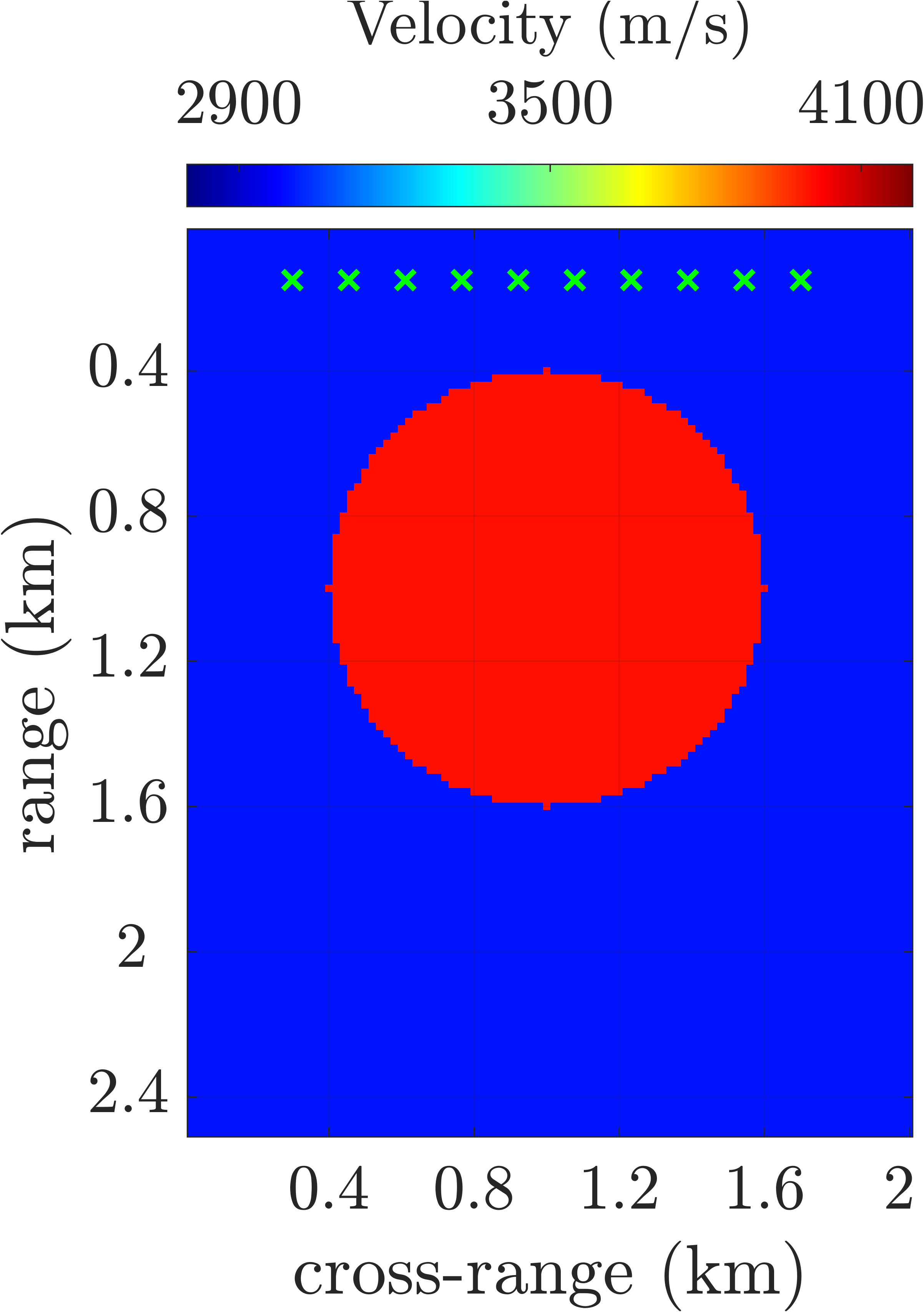}
\hfil
\includegraphics[height=5cm]{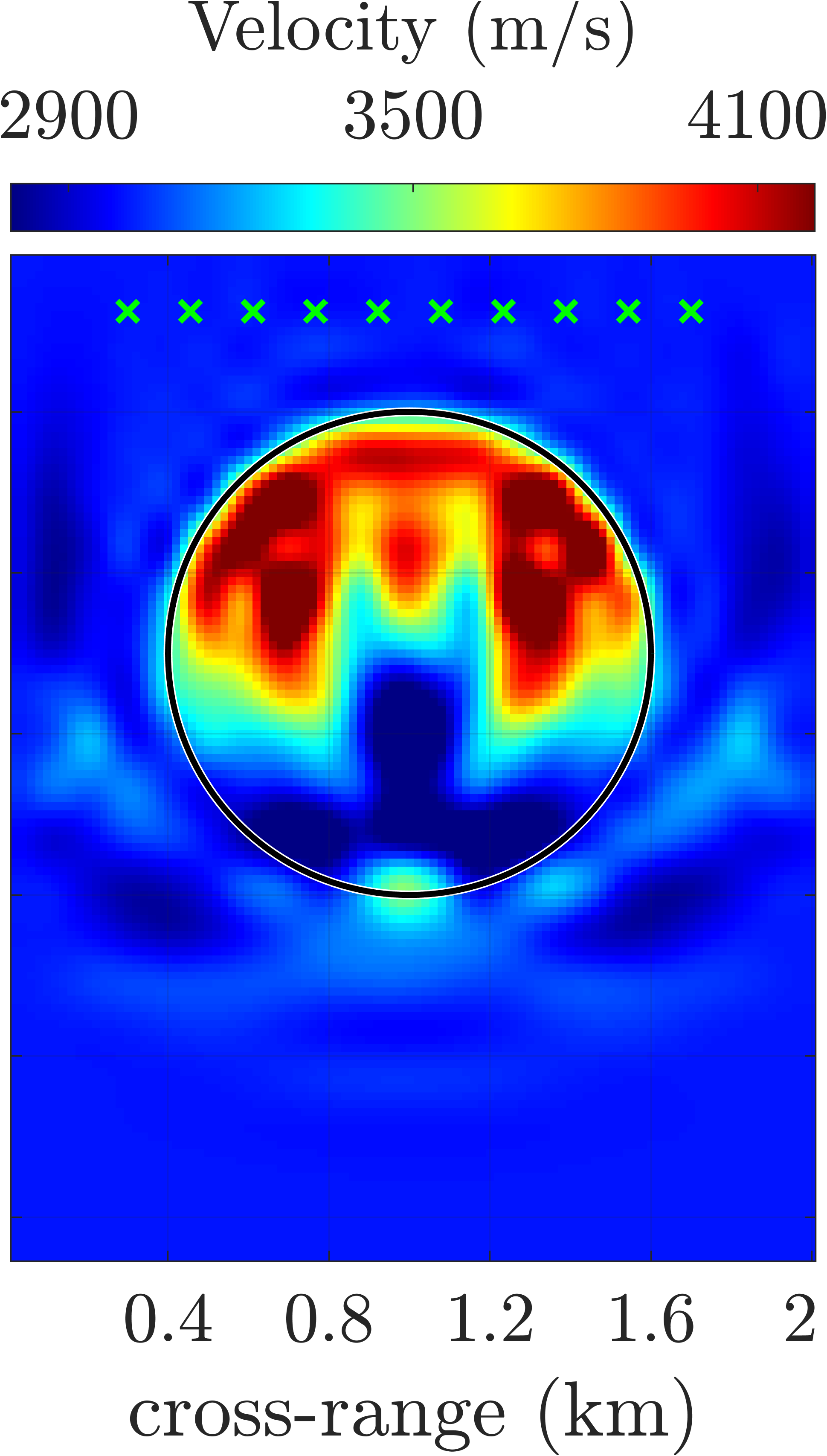}
\hfil
\includegraphics[height=5cm]{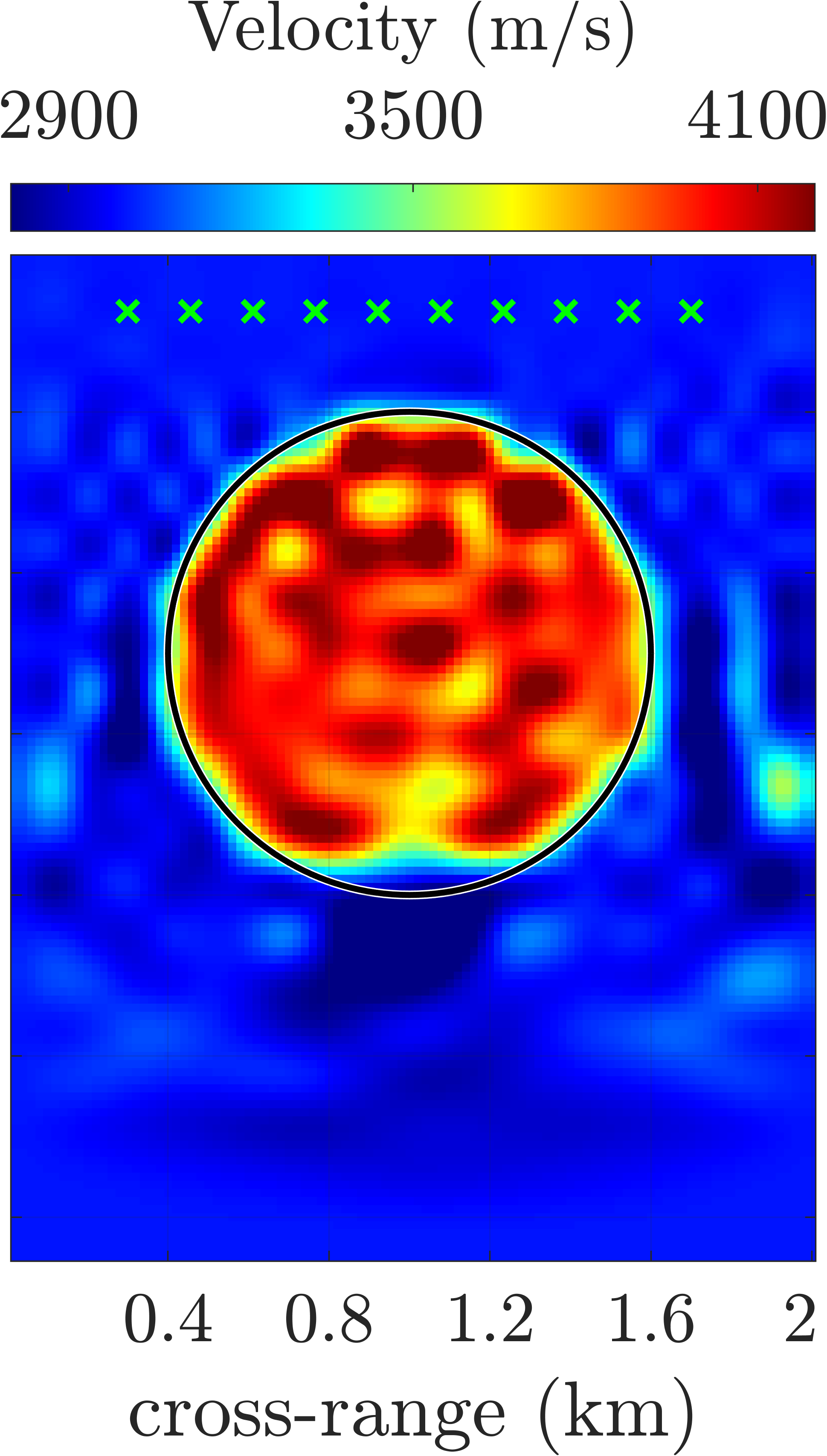}
\hfil
\includegraphics[height=5cm]{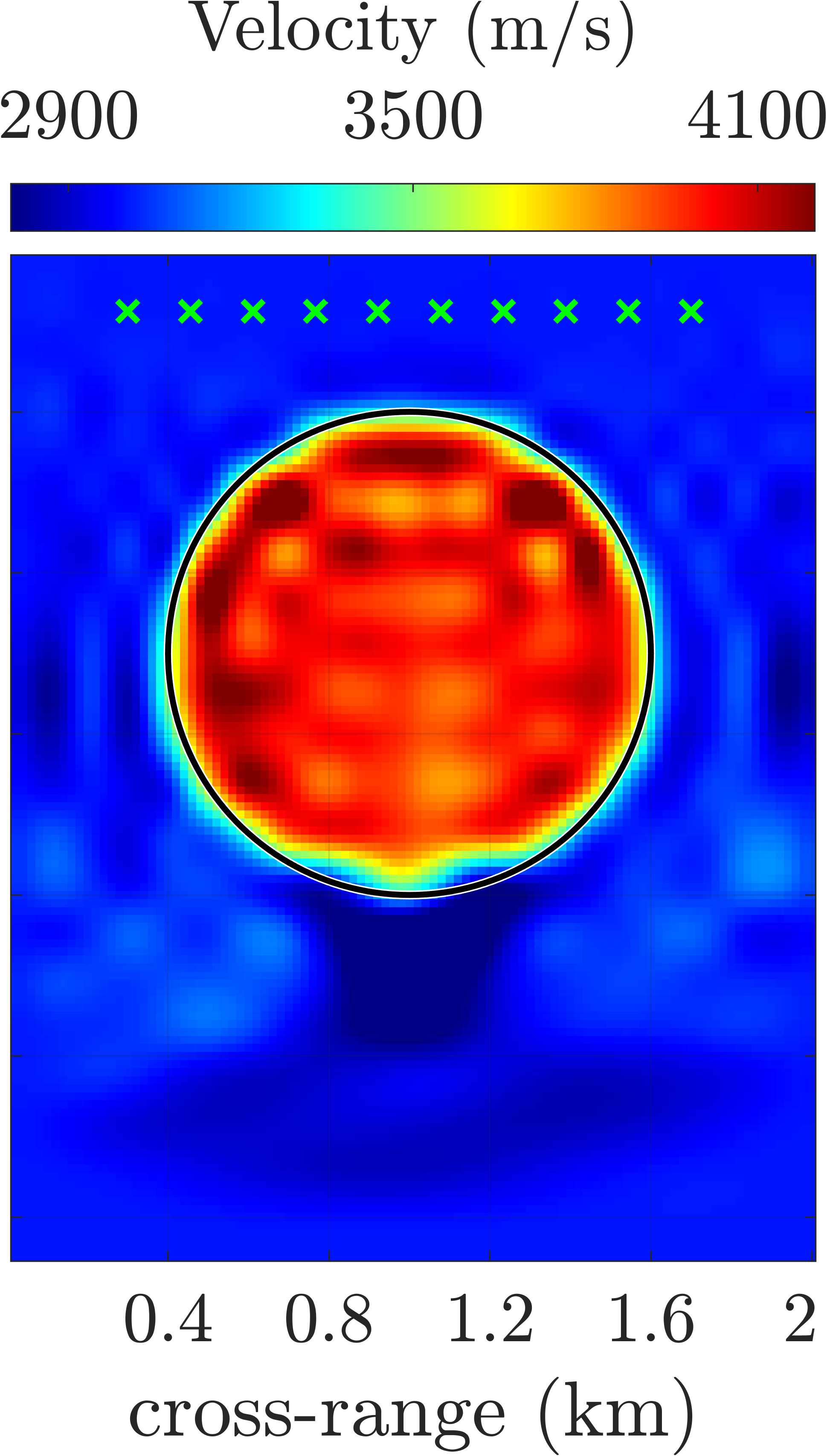} 
\end{center}
\caption{Inversion results: From left to right: True velocity $c(\bx)$, velocity estimated by 
FWI objective \eqref{eq:FWI}, by ROM objective  \eqref{eq:objR}, and by ROM objective \eqref{eq:objA}. 
The sources/receivers are shown with \bc{green $\times$}.}
\label{fig:Camembert_Inv}
\end{figure}

The first  results are for the Camembert model. We  show in Fig. \ref{fig:Camembert_Inv} the estimated wave speed after 60  Gauss-Newton iterations,   although the FWI approach stagnated after 30 iterations. The initial guess is the constant speed $\bar c = 3$km/s. Clearly, both ROM based methods give a good estimate of the fast inclusion, while FWI does not.

\begin{figure}[h]
\begin{center}
\includegraphics[width=0.32\textwidth]{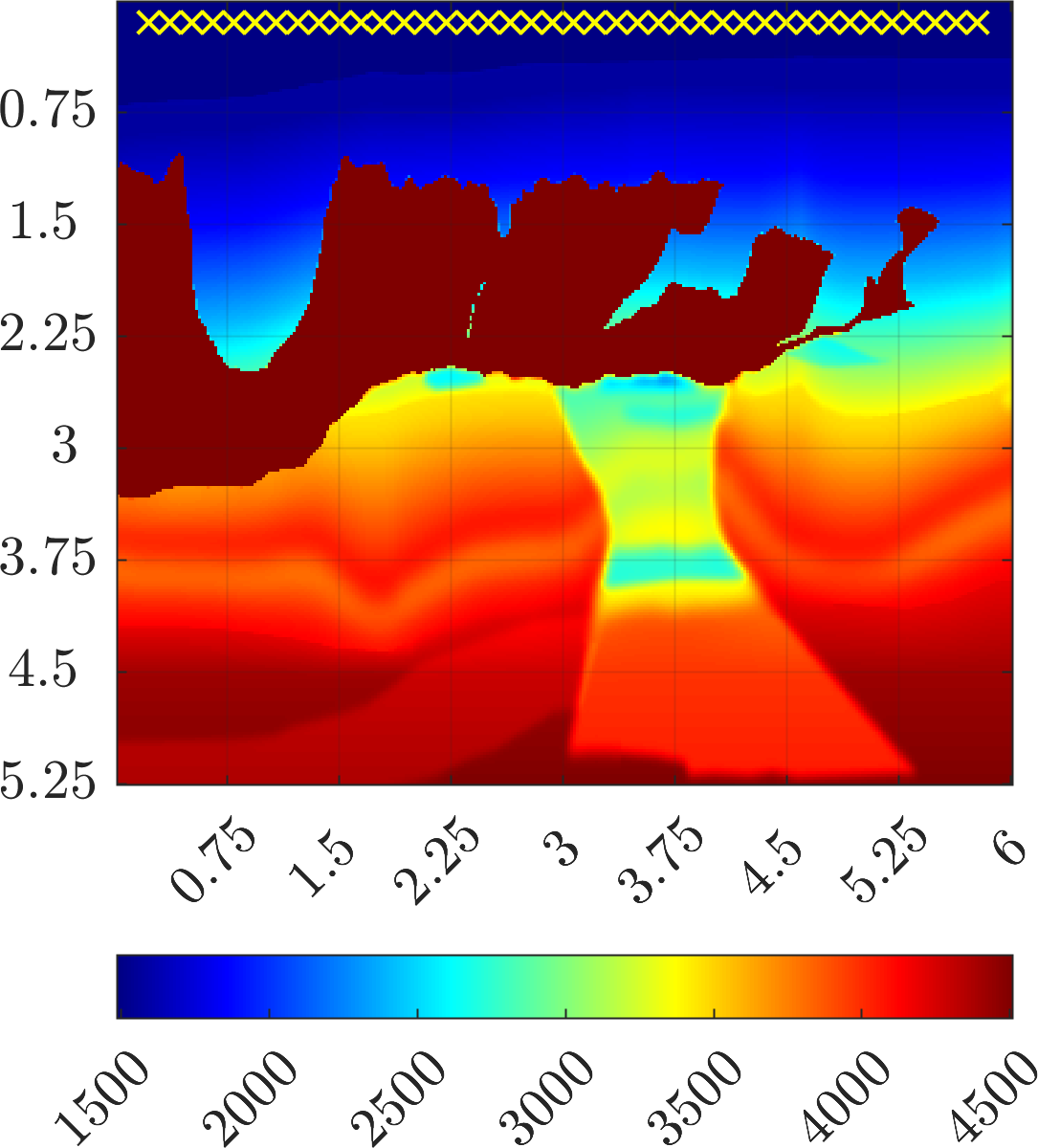}
\includegraphics[width=0.32\textwidth]{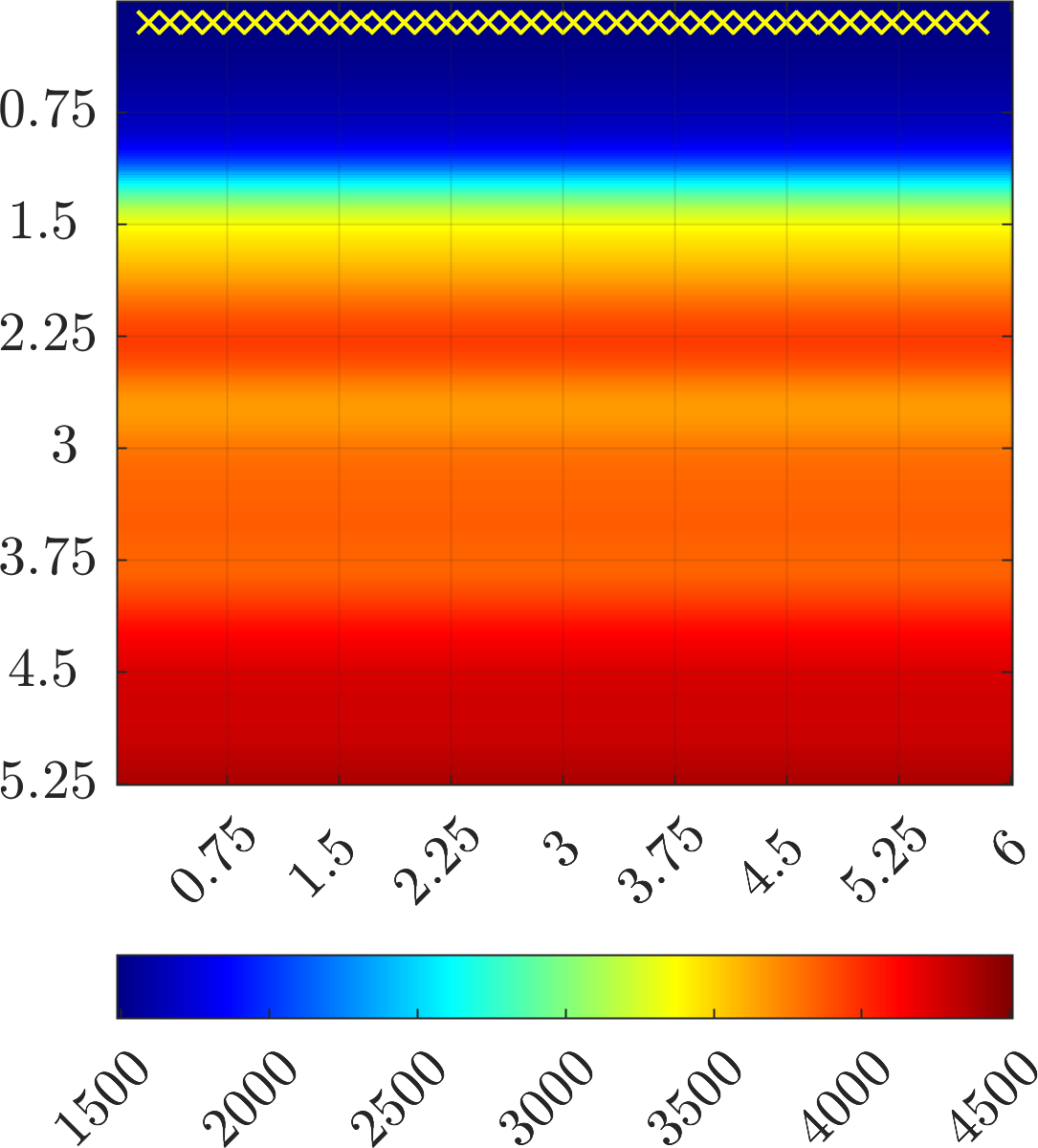} \\
\includegraphics[width=0.32\textwidth]{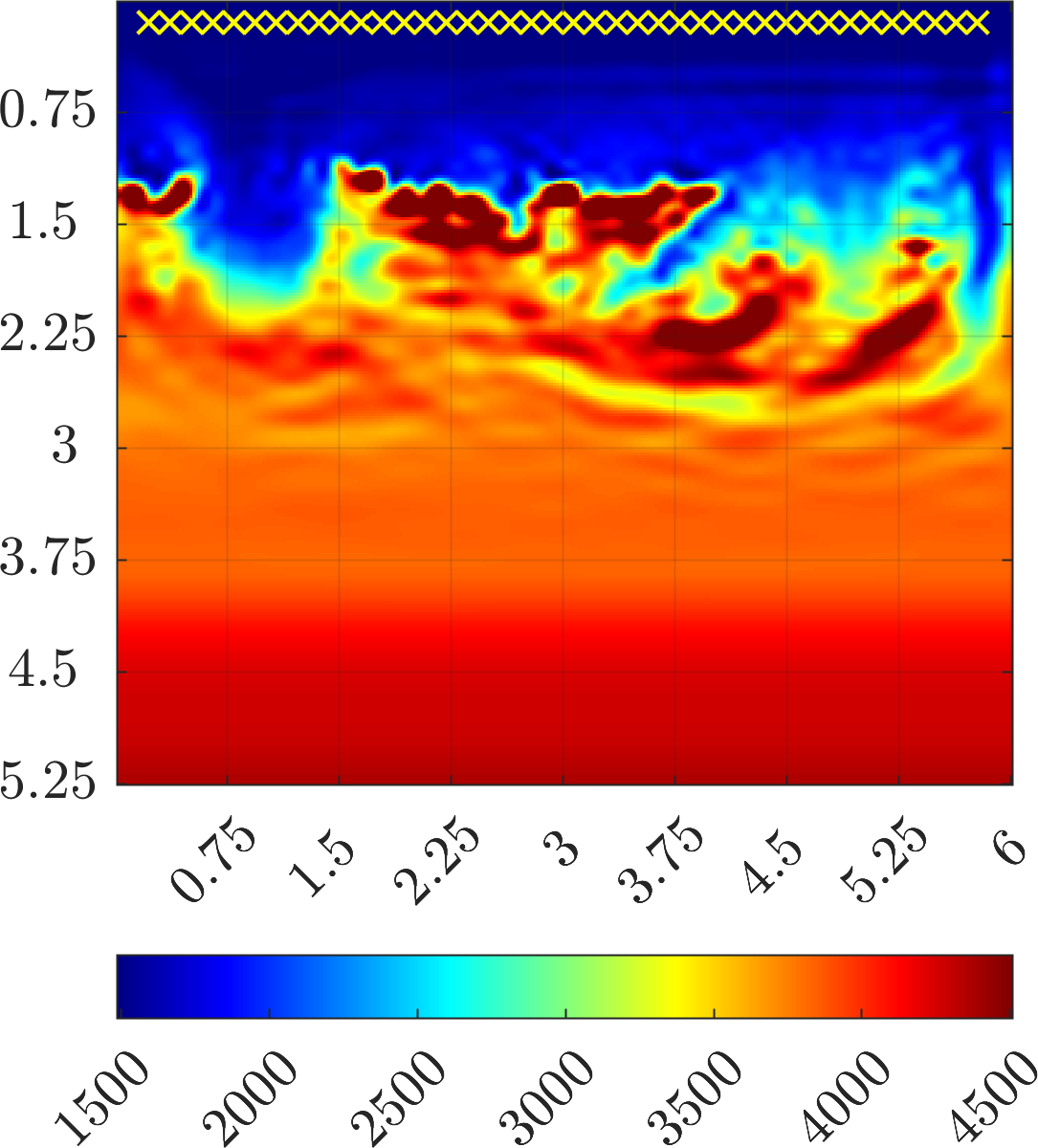}
\includegraphics[width=0.32\textwidth]{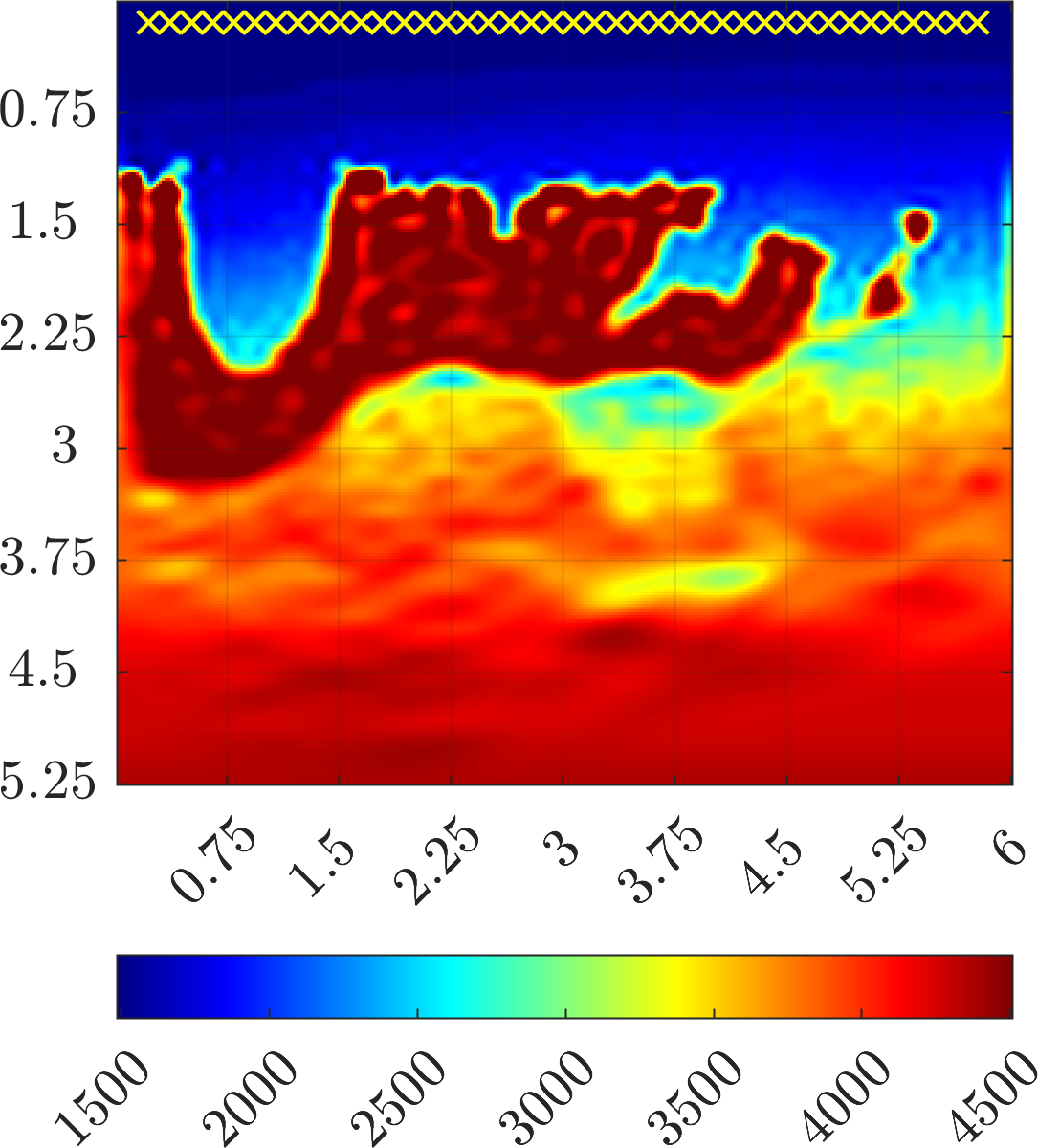}
\includegraphics[width=0.32\textwidth]{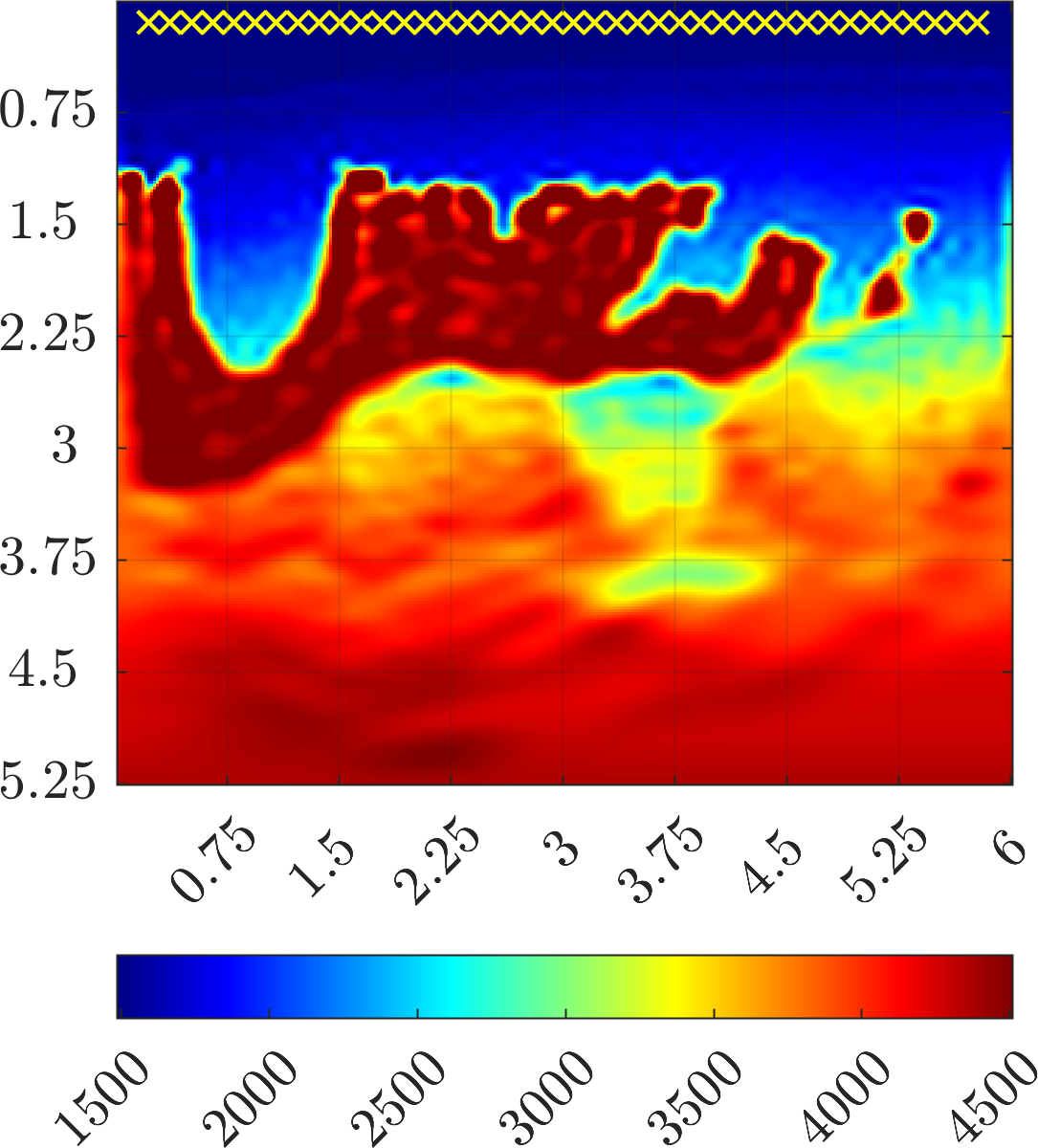} 
\end{center}
\caption{Inversion results: Top row, from left to right:  True $c(\bx)$ and the initial guess.
Bottom row, from left to right: velocity estimated by FWI objective \eqref{eq:FWI}, 
by ROM objective \eqref{eq:objR}, and by ROM objective \eqref{eq:objA}. 
The sources/receivers are shown with yellow $\times$.}
\label{fig:Salt_Inv}
\end{figure}

The second set of results is for the model displayed in the top left plot of Fig.~\ref{fig:Salt_Inv}. It was proposed in
\cite{billette2004} as a challenge problem, for imaging behind a salt body in the earth (the top fast structure in the figure, with wave speed $4.5$km/s). 
The initial guess is shown in the top right plot. The inversion results displayed in the bottom plots are after 35 Gauss-Newton iterations. They show that 
the ROM based methods perform similarly and image well. The FWI approach determines the top features of the salt body 
but then it gets stuck in a local minimum after 26 iterations.

\begin{figure}[h!]
\begin{center}
\begin{tabular}{cc}
True $c(\bx)$& Estimated wave speed \\
\includegraphics[width=0.39\textwidth]{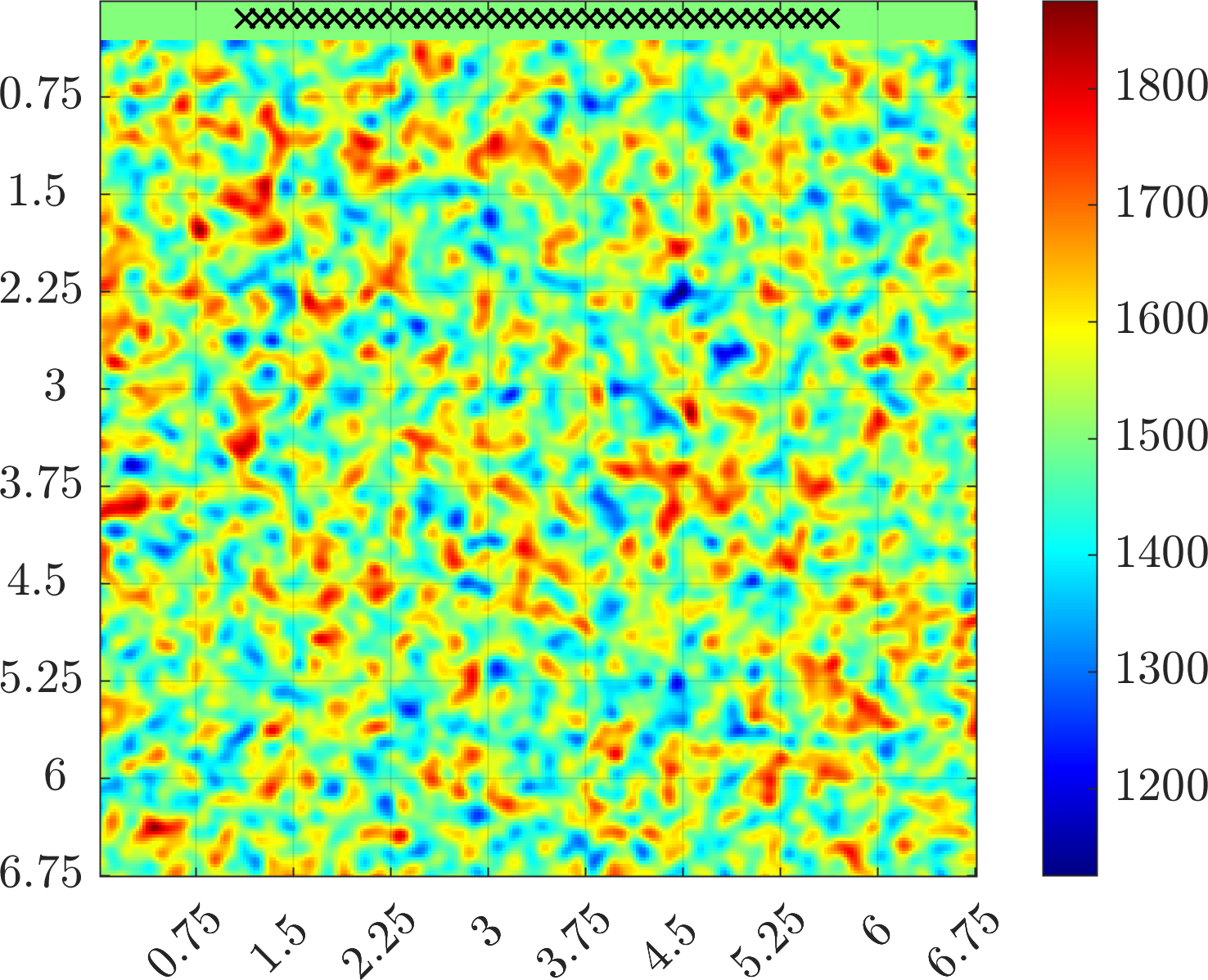} &
\includegraphics[width=0.39\textwidth]{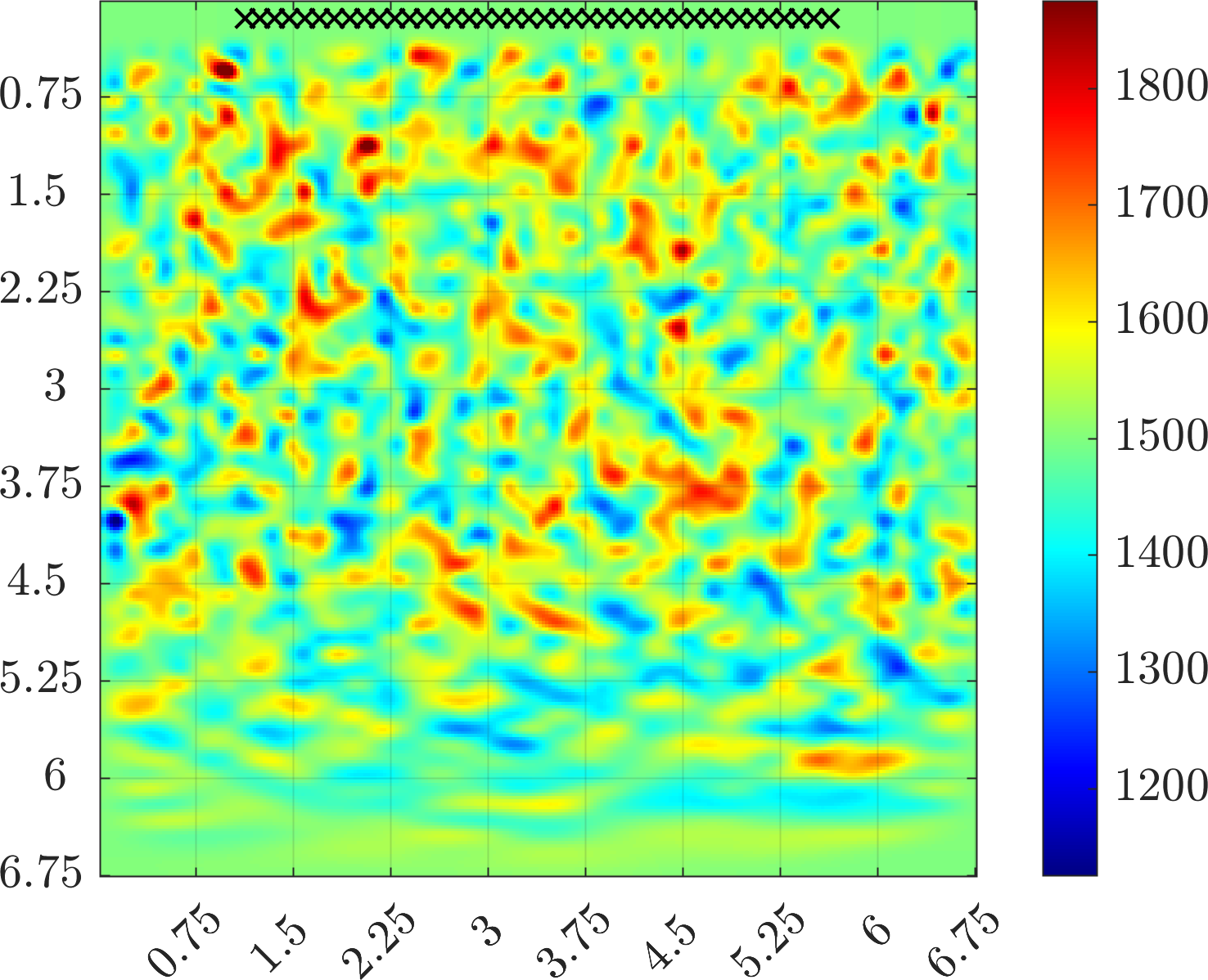}
\end{tabular}
\end{center}
\caption{Inversion result using the operator ROM. The sources/receivers are shown with black $\times$. The abscissa and ordinate are the cross-range and range in km and the wave speed values are in the colorbar, in m/s. }
\label{fig:Rand_Inv}
\end{figure}

\begin{figure}[h!]
\begin{center}
\begin{tabular}{ccc}
\includegraphics[width=0.305\textwidth]{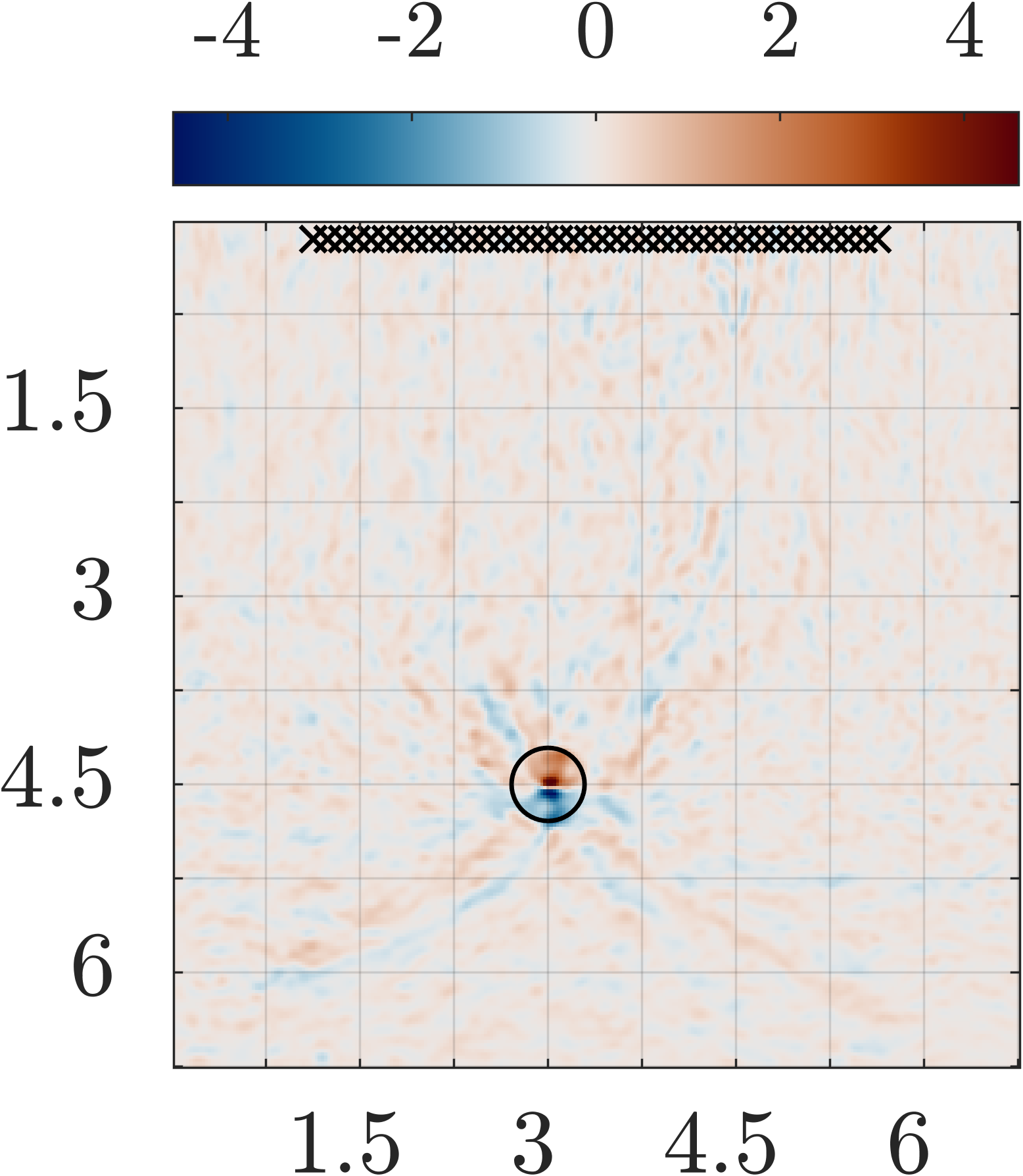} &
\includegraphics[width=0.305\textwidth]{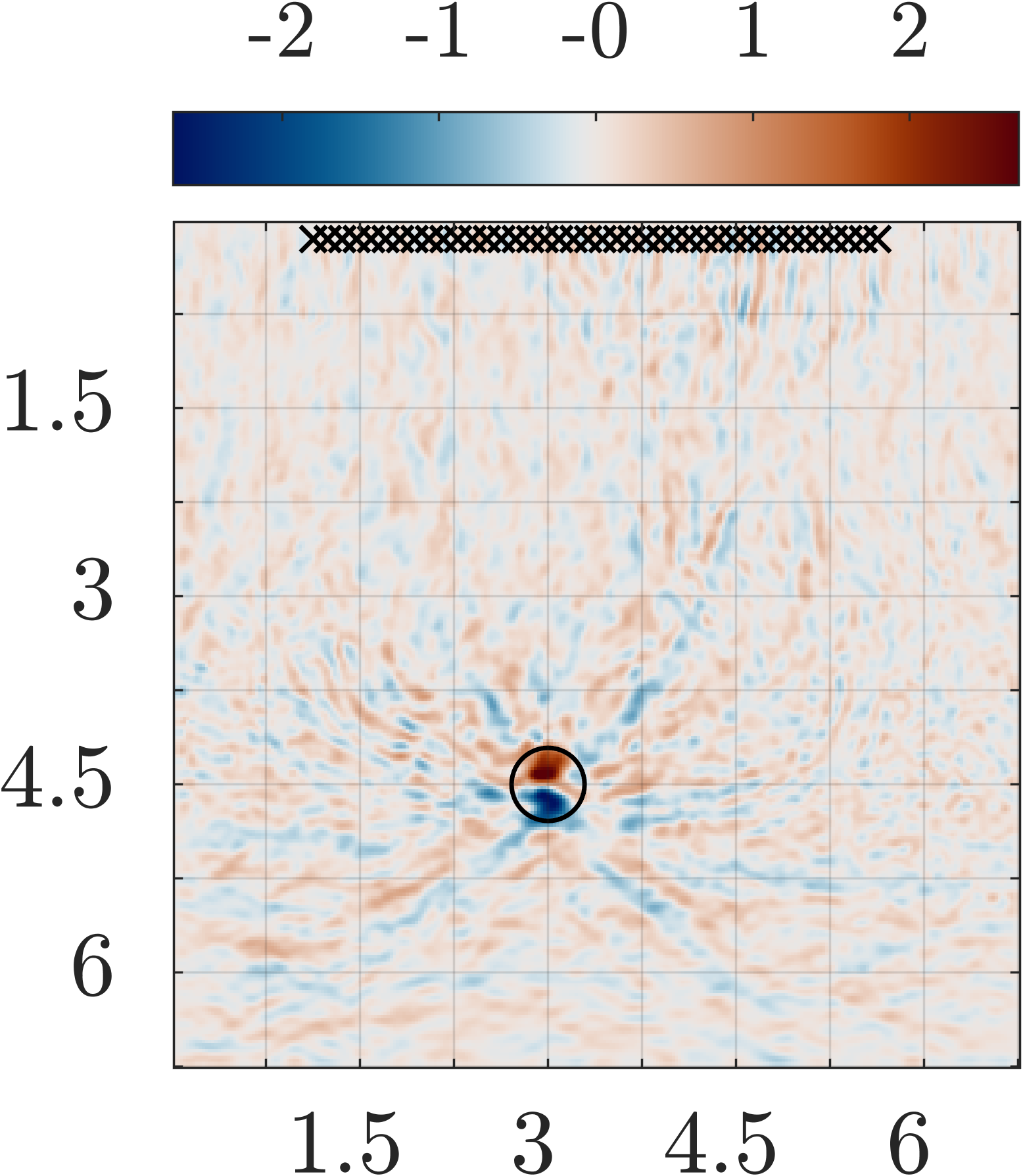} &
\includegraphics[width=0.305\textwidth]{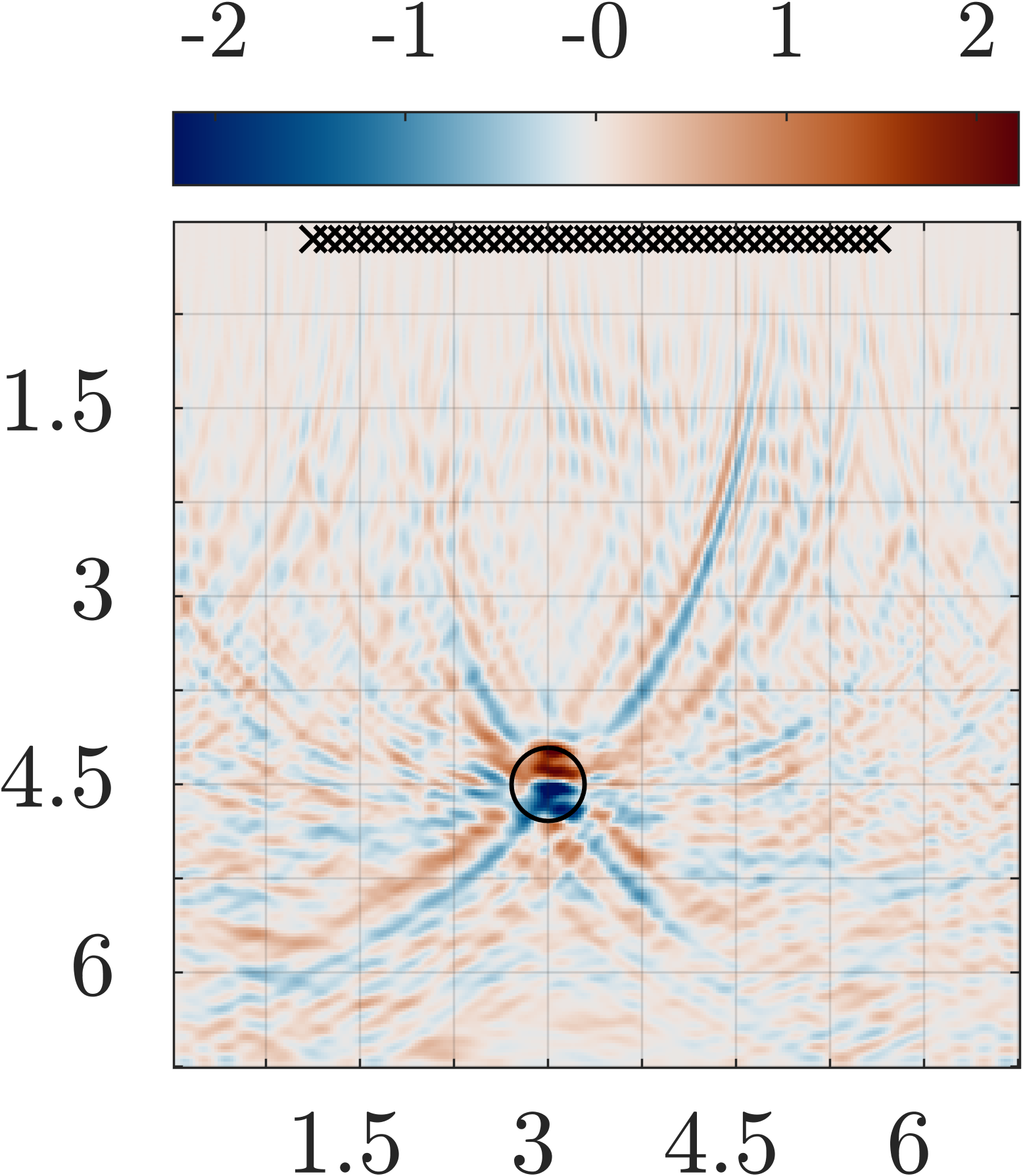}
\vspace{-2mm}
\end{tabular}
\end{center}
\caption{Time-reversal  focusing in the true medium (left), the estimated medium (middle) and the reference homogeneous medium (right) at a point labeled by a circle. The abscissa and ordinate are the cross-range 
and range in km. }
\label{fig:TR}
\end{figure}

The third (and last) results are for a realization of a smooth random medium (with Gaussian statistics and Gaussian covariance function), where the wave speed fluctuates about the reference value 
$\bar c = 1.5$km/s. The true medium is shown in the left plot and the reconstructed one is in the right plot of Fig. \ref{fig:Rand_Inv}. 
It is obtained after $135$ Gauss-Newton iterations for the operator ROM misfit \eqref{eq:objA}, with Tikhonov regularization. The initial guess is  the constant speed $\bar c = 1.5$km/s. 
We assess the quality of the reconstruction in two ways: First, we compute the correlation coefficient of the two images in Fig. \ref{fig:Rand_Inv} in the domain 
$[0.75\mbox{km}, 6\mbox{km}] \times [0.5\mbox{km}, 5.25\mbox{km}]$, with the MATLAB command ``corr2". A perfect reconstruction would give a correlation coefficient equal to $1$. 
The correlation coefficient for our reconstruction is $0.613$.  This may seem low, but we should remember that it is impossible to find the 
true medium due to physical constraints on the resolution limit.
The second  and better way to compare is to use a time-reversal experiment \cite{fink99,bookfouque}. 
\bc{A time-reversal experiment consists of two steps. 
In the first step waves transmitted by a controlled source and scattered by an unknown, complex medium are recorded by an array of sensors used as receivers.
In the second step the array of sensors used as sources transmit the time-reversed recorded signals. One observes a refocusing of the wave at the original source location due to the time reversibility of the wave equation. Moreover, time-reversal refocusing is sensitive to any relevant change in the medium between the two steps \cite{alfaro04,bal04}. 
This means that changes in the medium that affect the wave scattering process have a strong impact on the wave refocusing. 
A time-reversal experiment can be used, therefore, as a test to check whether an estimation of a wave speed by an imaging method is correct.
Here, we solve the wave equation with a point source at location $(3\mbox{km},4.5\rm{km})$ in the true medium and register the wave field at the array, over the duration $t \in [0, 5.8\mbox{s}]$. We time reverse these waves and send them back (computationally) into three media: the true one, the estimated one and the reference (homogeneous) medium. 
The back-propagated wave refocuses at the original source location in the true medium (left plot in Fig. \ref{fig:TR}). The more different the backpropagation medium is to the true one, the worse the refocusing is. Indeed, refocusing is quite poor in the reference medium (right plot in Fig. \ref{fig:TR}), but it is much improved in the estimated medium (middle plot in Fig. \ref{fig:TR}). This demonstrates that the estimated wave speed is similar to the one perceived by the waves in the true medium.}

\subsubsection{\bc{Computational cost}}
\label{sect:compnum}
\bc{Since both our inversion approaches use a Gauss-Newton iteration to minimize the objective functions \eqref{eq:objR}
or \eqref{eq:objA}, we compare their computational cost to that of the Gauss-Newton method for minimizing the FWI objective function \eqref{eq:FWI}. The same parameterization \eqref{eq:searchw} of the search velocity is assumed for all the approaches.
}

\bc{If the number $m$ of sensors and the number $2n-1$ of time steps is not too large, then the cost of each Gauss-Newton
iteration is dominated by the computation of the Jacobian of the objective function. For the FWI objective function, there is an
efficient way to compute the Jacobian, using the so-called adjoint formula \cite{virieux2009overview, symes2008migration}.
The analogues of such a formula for the objective functions \eqref{eq:objR} and \eqref{eq:objA} have not been derived, yet. 
Thus, our calculation of the Jacobian using finite differences, is the bottleneck of the computations. Aside from the Jacobian, 
the main computational burden of our approaches is due to the block-Cholesky factorization of the mass matrix, 
which involves $O\big((nm)^3\big)$ operations. }

\bc{If $nm \gg 1$, as it happens in three dimensional simulations, the dominant computational cost is  in solving the normal equations for the Gauss-Newton updates.  This  can be handled using iterative methods, such as conjugate gradient. To manage the extra computational burden of our approaches, one can  compute the  objective functions \eqref{eq:objR} and \eqref{eq:objA} from data gathered by sub-arrays, and then sum them in the optimization. This idea has been  used successfully for a different problem in \cite{borcea2014model}. In addition, one can exploit the algebraic structure of the 
wave operator ROM to define a new objective function that quantifies the misfit of a few block diagonals of $\bAR$   \cite{borcea2022waveform}.}

\section{Passive data acquisition}
\label{sect:passive}
So far, we assumed knowledge of the array response matrix  $ \boldsymbol{\mathcal{M}}(t)$, gathered by co-located sources and receivers, which we then transformed to the new data matrix $ \bD(t)$, as stated in Lemma \ref{lem.1}. 
In this section we consider  uncontrolled, ambient noise sources that emit stationary random signals and show how the cross correlations of the generated waves measured at a passive array of receivers can give directly the  matrix $ \bD(t)$. 
This opens the door to new applications where the use of controlled active sources is neither possible nor desired.

The fact that the cross correlations of the signals measured by a passive array  can give the same information as 
 an active array  is not new \cite{lobkis2001} and it has been exploited for many problems in imaging and free space communications \cite{garnierpassive2009,garnierpassive2010,garnier2016passive,passivecomm}.
The mathematical statement is that the cross correlations of the signals are related to the symmetrized Green's function of the wave equation.
This relationship holds for different situations, in open media with radiation conditions and in bounded cavities. To stay consistent with the rest of the paper,  we address here the case of a bounded cavity with homogeneous Dirichlet boundary conditions.
 
For a well-posed mathematical formulation of the problem with  random noise sources, we need to introduce some dissipation. 
Thus, we consider the solution $p(t,\bx)$ of the damped wave equation
\begin{equation}
\label{eq:ondes1}
(T_a^{-1}+\partial_t)^2
 p - c^2(\bx) \Delta  p = s(t,\bx),\qquad t \in \RR, \quad \bx \in \Omega \subset \RR^d,
\end{equation} 
in dimension $d \in \{2,3\}$, where the forcing term $s (t,\bx)$ is a zero-mean, stationary in time random process. 
As we did in section \ref{sect:ROM.1}, we apply the similarity transformation 
(\ref{eq:Sim1}) from $p(t,\bx)$ to $P(t,\bx)$, that acts as the identity at the points in the receiver array. This transforms the wave operator to $(T_a^{-1}+\partial_t)^2+\cA$, with the self-adjoint operator $\cA$ defined by (\ref{eq:D4}). The source term becomes 
\begin{equation}
{S}(t,\bx) = \frac{\bar c}{c(\bx)} s(t,\bx),
\end{equation}
and we model its autocorrelation by 
\begin{equation}
\langle S(t_1,\by_1) S(t_2,\by_2) \rangle =  F(t_1-t_2) {K} ( \by_1) \delta(  \by_1 -\by_2 )  ,
\label{eq:corcavity}
\end{equation}
where $\langle \cdot \rangle $ denotes  a statistical average. In this equation, the function ${K}(\by)$ characterizes the spatial support of the noise sources, which is assumed to be contained in $\Omega$.
The sources are  delta-correlated in space, hence the Dirac delta in \eqref{eq:corcavity}, and 
the covariance in time depends only on the time offset $t_1-t_2$, because of stationarity.
The function $F(t)$ decays to $0$ at infinity, it is even and it belongs to $L^1$ (which gives it ergodic properties).
The Fourier transform $\hat{F}(\om)$ of $F(t)$ is the power spectral density of the noise sources, which is nonnegative by Bochner's theorem.

For any $\by \in \Omega$, 
the Green's function $(t,\bx) \mapsto G(t,\bx,\by)$ is the solution of 
\begin{equation}
(T_a^{-1}+\partial_t)^2G +\cA G  = \delta(t) \delta(\bx-\by)  ,\qquad t \in \RR, \quad \bx \in \Omega ,
\end{equation}
with the initial condition
$G(t,\bx,\by)=0$ for all $t <0$ and the boundary condition $G(t,\bx,\by)=0$ for $\bx \in \partial \Omega$. It can be written 
explicitly  in terms of the
eigenvalues $\{\theta_j\}_{j \geq 1}$ and orthonormal eigenfunctions $\{y_j(\bx)\}_{j \geq 1}$ 
of $\cA$, as follows
\begin{equation}
\label{eq:expandgreendelta}
G(t,\bx,\by) = H(t)  \exp(-t/T_a)
\sum_{j=1}^\infty \frac{\sin(\sqrt{\theta_j} t)}{\sqrt{\theta_j} } y_j(\bx) y_j(\by),
\end{equation}
where we recall that $H(t)$ is the Heaviside step function.

The  empirical cross correlation of the recorded waves at $\bx_r$ and $\bx_{r'}$ is defined by
\begin{equation}
C_T (\tau,\bx_r,\bx_{r'}) = \frac{1}{T} \int_0^T dt \, 
P(t,\bx_r) P(t+\tau,\bx_{r'}), \label{eq:empirical}
\end{equation}
and the statistical cross correlation 
\begin{equation}
C^{(1)} (\tau,\bx_{r},\bx_{r'}) = \big< C_T (\tau,\bx_{r},\bx_{r'}) \big>,
\end{equation}
 is independent of $T$ by stationarity of the noise sources.
The statistical stability of \eqref{eq:empirical} follows from the ergodicity of the noise sources \cite{garnier2016passive} and we have
\begin{equation}
C_T (\tau,\bx_{r},\bx_{r'}) \stackrel{T \to +\infty}{\longrightarrow} C^{(1)} (\tau,\bx_{r},\bx_{r'}), 
\end{equation}
in probability.
The relation between the statistical cross correlation and the Green's function
can be obtained by using the normal mode expansion   \eqref{eq:expandgreendelta}. The next proposition, stated for an idealized situation for simplicity, shows that the $\tau$-derivative of the cross correlation
is a smoothed and symmetrized version of the Green's function:

\begin{proposition}
\label{prop:normal}%
Consider an inhomogeneous cavity with homogeneous dissipation quantified by $1/T_a$. Suppose that 
the source distribution extends over the whole cavity i.e., $K(\by) ={\bf 1}_\Omega(\by)$ in (\ref{eq:corcavity}).
Then, for any $r,r'=1,\ldots,m$, we have 
\begin{align}
{\partial_\tau} C^{(1)}(\tau,\bx_{r},\bx_{r'}) =
-\frac{T_a}{4} 
F (\tau) \star_\tau \big[  {\rm sgn}(\tau) G( |\tau| ,\bx_{r},\bx_{r'}) 
 \big] ,
\label{eq:kernwhite0}
\end{align}
where  $G$ is the Green's function (\ref{eq:expandgreendelta}) and ${\rm sgn}(t)=H(t)-H(-t)$ is the sign function.
\end{proposition}

The proof can be found in \cite[Section 2.5]{garnier2016passive}.
Note that $C^{(1)}$ is an even function in $\tau$, and therefore its $\tau$-derivative is an odd function.

We can now give the relation between the cross correlation matrix and the data matrix $\bD(t)$ exploited in the ROM procedure. 

\begin{thm}
\label{thm.4}
We have, for any $r,r'=1,\ldots,m,$ that 
\begin{align}
\frac{1}{T_a} \partial_\tau^2 C^{(1)} (\tau,\bx_{r},\bx_{r'}) 
\stackrel{T_a \to +\infty}{\longrightarrow}
 -\frac{1}{4} D_{r,r'}(\tau) ,
\end{align}
where $\bD(t)$ is the matrix defined in Lemma \ref{lem.1}, with a signal $f(t)$  whose Fourier transform satisfies $|\hat{f}(\om)|=\hat{F}^\12(\om)$.
\end{thm}

\textbf{Proof:} 
From Proposition \ref{prop:normal} we obtain
\begin{align*}
\frac{1}{T_a} \partial_\tau^2 C^{(1)} (\tau,\bx_{r},\bx_{r'}) 
\stackrel{T_a \to +\infty}{\longrightarrow}
-\frac{1}{4}  \sum_{j=1}^\infty  \int_\RR dt \, F(\tau-t) \cos \big(\sqrt{\theta_j} t\big)  \,
y_j (\bx_{r}) y_j(\bx_{r'}) .
\end{align*}
Now use that 
\begin{align*}
\int_\RR dt \, F(\tau-t) \cos (\sqrt{\theta_j} t )  = \int_\RR d \om \,  \hat F(\om) e^{i \om \tau} \hspace{-0.05in} \int_{\RR} \frac{dt}{2 \pi} \cos (\sqrt{\theta_j} t ) e^{i \om t} = 
\hat{F}(\sqrt{\theta_j}) \cos (\sqrt{\theta_j} \tau) 
\end{align*}
and obtain the desired result after deducing from equations \eqref{eq:D6} and \eqref{eq:D6p} that 
 \begin{align*}
   D_{r,r'}(\tau)  = \sum_{j=1}^\infty  \hat{F} (\sqrt{\theta_j}) \cos (\sqrt{\theta_j}\tau) 
y_j (\bx_{r}) y_j(\bx_{r'}). \qquad \Box
\end{align*}

Theorem \ref{thm.4} says that we can approximate $\bD(t)$  when $\tau \ll T_a$ i.e., in the presence of weak attenuation.
It shows that our ROM procedure is natural in the passive framework, since the cross correlation of the noise signals recorded by an array of receivers gives directly 
the data matrix $\bD(t)$ corresponding to  a virtual active array. That is to say, the 
virtual sources and receivers are naturally co-located and the signals are even (because cross correlations are even). This is what we need for the ROM-based procedure.

\begin{rem}
Here we considered a cavity $\Omega$ with uniformly distributed noise sources in $\Omega$. In this idealized case, no hypothesis regarding the geometry of the cavity  is needed 
to obtain the result in  Proposition \ref{prop:normal}.
However, this result also holds when the noise source distribution is non-homogeneous or spatially
limited, provided the cavity possesses some ergodic properties and the attenuation 
time $T_a$ is larger than the critical time necessary to reach ergodicity
 \cite{bardos2008,colin2009}.
 \end{rem}

\section{Concluding remarks and open questions}
\label{sec:conc}%
Our goal in this paper was to show how computationally efficient tools from numerical linear algebra and reduced order modeling 
can be used to improve the existing inverse scattering methodology. We presented in a unified way two  projection type reduced order models (ROMs) that
capture wave propagation in complementary ways and  have the following important properties for solving waveform inversion problems: 
(1) They can be computed from measurements of the waves that are available in most application setups. These include data acquisitions 
 with controlled sources that emit probing pulses and nonconventional ones with uncontrolled noise sources. (2) The mapping from the data  to the ROMs is nonlinear and yet, the ROM computation can be carried out with known linear algebra algorithms, in a non-iterative fashion. (3) This computation can be done so the causal physics of wave propagation is captured by a special algebraic structure of the ROM matrices. This is essential for the success of the inversion. 

We described \bc{two ideas}  for using the ROMs for inversion.  They are quite a departure from the data fitting approaches found in the literature. The inversion results based on the ROMs are promising, but the methodology  is young and there are many open questions. Here are a few examples: 

1.  The dependence of the projection basis $\bV(\bx)$ on $c(\bx)$ is understood only in a few media \cite[Appendix A]{borcea2021reduced}. This basis also depends on the choice of $\tau$ and the separation between the sensors in the array. A deeper
understanding of $\bV(\bx)$ would lead to a more rigorous foundation of waveform inversion based on $\bAR$ and of ROM based 
imaging methods like   \cite{borcea2021reduced,druskin2018nonlinear,DtB}.

2. How can one extend the ROM methodology to setups where the sources are still controlled by the user, but are placed in different locations than the receivers?
The gathered data matrices are no longer symmetric in such cases and a different projection methodology should be used, of the Petrov-Galerkin type. It is not difficult to get ROMs that do a good job in terms of approximating the forward map in such settings. However, it is not yet known how to get ROMs that are useful for inversion.

3. The proposed ROM construction relies on having a self-adjoint wave operator, which is no longer the case if there is attenuation. If the attenuation coefficient $T_a$ is constant  and known, as assumed in  section \ref{sect:passive}, then it is possible to remove attenuation effects  by multiplying the measurements  with $\exp(t/T_a)$. This is feasible if the attenuation is weak. The case of variable attenuation requires a fundamental rethinking of the methodology. So is the case of dispersive media.

4. The ROM based inversion methodology extends to media with variable mass density, to anisotropic media and to electromagnetic waves. Inversion with elastic waves remains largely unexplored. The main difficulty there is due to the 
multiple wave modes that propagate at different wave speeds.  

We hope that these questions and likely many others will motivate the applied and computational mathematics community to explore this new way of solving inverse wave scattering problems.

\section* {Acknowledgements} We thank  Vladimir Druskin, Mike Zaslavsky and  Andy Thaler for their foundational work on  the data driven ROM methodology. This material is based upon research supported in part by the U.S. Office of Naval Research 
under award number N00014-21-1-2370 (program manager Dr. Reza Malek-Madani) to Borcea and Mamonov; by the AFOSR award number FA9550-22-1-0077 (program manager Dr. Arje Nachman) to Borcea and Garnier,  and by the National Science Foundation under Grant No. 2110265 to Zimmerling.\\

We thank Hongyu Zhai for finding the mistake in the statement of Theorem 3 and Corollary 1, which is now corrected.

\appendix
 \section{Setup for the numerical simulations}
 \label{ap:A}
 All the simulations are in a two-dimensional rectangular domain with sound soft boundary.  
 The data are computed  with a time-domain solver for the wave equation \eqref{eq:I1}, 
 with a five point stencil discretization of the Laplacian  on a uniform grid with mesh size 
$20\mbox{m}$ for Fig. \ref{fig:Camembert_Inv},
$18.75\mbox{m}$ for Fig. \ref{fig:Salt_Inv} and
$25\mbox{m}$ for Fig. \ref{fig:Rand_Inv}.
The second time derivative is approximated by a three point finite difference scheme, on a time grid 
with step $\Delta t= \tau / 20$.  The data matrices are computed using definition \ref{eq:D5}. 

We now explain how we compute the 
second derivative data matrices. We begin by extending the finely sampled data evenly in discrete time
to get $\bD_{e,j}$, $j = -N,\ldots,N$, with $\bD_j = \bD_{e,\pm j}$, $j=0,1,\ldots,N$. 
Then, we take the discrete Fourier transform of $( \bD_{e,j} )_{j=-N}^{N}$ and differentiate 
in the Fourier domain, after using a sharp cutoff low-pass filter intended to stabilize the calculation. 
The cutoff frequency is at $(\om_o+4 B)/(2 \pi)$. We take the  inverse Fourier transform to 
obtain $\ddot{\bD}_{e,j}$, at $j = -N,\ldots,N$, the finely sampled second derivative data. 
Finally,  we sub-sample $\bD_{e,j}$ and $\ddot{\bD}_{e,j}$ to get 
\begin{equation}
\bD_j = \bD_{e,20j}, \quad \ddot{\bD}_j = \ddot{\bD}_{e,20j}, \quad j = 0,1,\ldots,2n
-1.
\end{equation}

\section{Parameters for the inversion simulations}
\label{ap:A1}
All the simulations use the  probing signal
 $
 f(t) = \cos(\om_o t) \exp(-B^2 t^2/2),
 $
 with $\om_o/(2 \pi) = 6$Hz and $B/(2 \pi)  = 4$Hz.  
The array of $m$ sensors (sources/receivers) is near the top boundary, at depth equal to the central wavelength.

For Fig. \ref{fig:Camembert_Inv}:  $\Omega = [0,2 \mbox{km}] \times [0,2.5 \mbox{km}]$ and $\Omega_{\rm inv} = [95 \mbox{m},1905 \mbox{m}] \times [119 \mbox{m},2381 \mbox{m}]$. The reference speed is $\bar c = 3$km/s, which gives  $\lambda_o = 300$m. 
We use an array of $m = 10$ sensors with aperture $1400\mbox{m}$, $\tau = 0.0435 s$ and $n = 16$. 
We use $N_t = 6$ time windows in the inversion, and parametrize the search speed with $N = 20 \times 20 = 400$ Gaussian basis functions, with standard deviation $\sigma^\perp = 55.5 \mbox{m}$ and $\sigma = 69.4 \mbox{m}$. 

For Fig. \ref{fig:Salt_Inv}: $\Omega = [0, 6 \mbox{km}] \times [0, 5.25 \mbox{km}]$ and 
$\Omega_{\rm inv} = [105 \mbox{m}, 5895 \mbox{m}] \times [ 92 \mbox{m}, 5158 \mbox{m}]$. 
The reference speed is $\bar c = 1.5$km/s, which gives $\lambda_o = 150$m. 
We use an array of $m = 40$ sensors with aperture of $5550 \mbox{m}$, $\tau = 0.0333 \mbox{s}$ and $n = 49$. 
We use $N_t = 6$ time windows in the inversion, and parametrize the search speed with 
$N = 55 \times 55 = 3025$ Gaussian basis functions, with standard deviation 
$\sigma^\perp = 61.2 \mbox{m}$ and $\sigma = 53.6 \mbox{m}$.

For Fig. \ref{fig:Rand_Inv}: $\Omega = [0,6.5 \mbox{km}] \times [0,6.5 \mbox{km}]$ and 
$\Omega_{\rm inv} = [114 \mbox{m}, 6386 \mbox{m}] \times [300 \mbox{m}, 6392 \mbox{m}]$. 
We have $\bar{c} = 1.5$km/s, which gives $\lambda_o = 150$m.   The fluctuations of $c(\bx)$  
have a contrast of $25\%$ and they occur on a length scale of the order $\lambda_o$. 
We use an array of $m = 40$ sensors with aperture $30\lambda_o = 4.5$km. 
We use $N_t = 9$ time windows and parametrize the search space with 
$N = 55 \times 55 = 3025$ Gaussian basis functions, with standard deviation 
$\sigma^\perp = 66.6 \mbox{m}$ and $\sigma = 63.4 \mbox{m}$.

\bibliographystyle{siam} 
\bibliography{biblio.bib}

\begin{thebibliography}{10}

\bibitem{alfaro04}
{\sc D.~G. Alfaro~Vigo, J.-P. Fouque, J.~Garnier, and A.~Nachbin}, {\em
  Robustness of time reversal for waves in time-dependent random media},
  Stochastic processes and their applications, 113 (2004), pp.~289--313.

\bibitem{antoulas2001survey}
{\sc A.~Antoulas, D.~Sorensen, and S.~Gugercin}, {\em A survey of model
  reduction methods for large-scale systems}, Contemporary mathematics, 280
  (2001), pp.~193--220.

\bibitem{arridge2012imaging}
{\sc S.~Arridge and O.~Scherzer}, {\em Imaging from coupled physics}, Inverse
  Problems, 28 (2012), p.~080201.

\bibitem{asvadurov2003optimal}
{\sc S.~Asvadurov, V.~Druskin, M.~Guddati, and L.~Knizhnerman}, {\em {On
  optimal finite-difference approximation of PML}}, SIAM Journal on Numerical
  Analysis, 41 (2003), pp.~287--305.

\bibitem{bal2013reconstruction}
{\sc G.~Bal and G.~Uhlmann}, {\em Reconstruction of coefficients in scalar
  second-order elliptic equations from knowledge of their solutions},
  Communications on Pure and Applied Mathematics, 66 (2013), pp.~1629--1652.

\bibitem{bal04}
{\sc G.~Bal and R.~Ver{\'a}stegui}, {\em Time reversal in changing
  environments}, SIAM Multiscale Modeling \& Simulation, 2 (2004),
  pp.~639--661.

\bibitem{bardos2008}
{\sc C.~Bardos, J.~Garnier, and G.~Papanicolaou}, {\em Identification of
  {G}reen's functions singularities by cross correlation of noisy signals},
  Inverse Problems, 24 (2008), p.~015011.

\bibitem{belishev2007recent}
{\sc M.~Belishev}, {\em Recent progress in the boundary control method},
  Inverse problems, 23 (2007), p.~R1.

\bibitem{benner2015survey}
{\sc P.~Benner, S.~Gugercin, and K.~Willcox}, {\em A survey of projection-based
  model reduction methods for parametric dynamical systems}, SIAM Review, 57
  (2015), pp.~483--531.

\bibitem{billette2004}
{\sc F.~Billette and S.~Brandsberg-Dahl}, {\em {The 2004 BP velocity benchmark:
  67th Annual EAGE Meeting, EAGE}}, in Expanded Abstracts, vol.~305, 2004.

\bibitem{bleistein2001multidimensional}
{\sc N.~Bleistein, J.~Stockwell, and J.~Cohen}, {\em Multidimensional seismic
  inversion}, Springer, 2001.

\bibitem{borcea2008electrical}
{\sc L.~Borcea, V.~Druskin, and F.~Guevara-Vasquez}, {\em Electrical impedance
  tomography with resistor networks}, Inverse Problems, 24 (2008), p.~035013.

\bibitem{borcea2005continuum}
{\sc L.~Borcea, V.~Druskin, and L.~Knizhnerman}, {\em On the continuum limit of
  a discrete inverse spectral problem on optimal finite difference grids},
  Communications on Pure and Applied Mathematics, 58 (2005), pp.~1231--1279.

\bibitem{borcea2010circular}
{\sc L.~Borcea, V.~Druskin, and A.~Mamonov}, {\em Circular resistor networks
  for electrical impedance tomography with partial boundary measurements},
  Inverse Problems, 26 (2010), p.~045010.

\bibitem{borcea2020reduced}
{\sc L.~Borcea, V.~Druskin, A.~Mamonov, M.~Zaslavsky, and J.~Zimmerling}, {\em
  Reduced order model approach to inverse scattering}, SIAM Journal on Imaging
  Sciences, 13 (2020), pp.~685--723.

\bibitem{borcea2014model}
{\sc L.~Borcea, V.~Druskin, A.~V. Mamonov, and M.~Zaslavsky}, {\em A model
  reduction approach to numerical inversion for a parabolic partial
  differential equation}, Inverse Problems, 30 (2014), p.~125011.

\bibitem{DtB}
\leavevmode\vrule height 2pt depth -1.6pt width 23pt, {\em Untangling the
  nonlinearity in inverse scattering with data-driven reduced order models},
  Inverse Problems, 34 (2018), p.~065008.

\bibitem{borcea2019robust}
\leavevmode\vrule height 2pt depth -1.6pt width 23pt, {\em Robust nonlinear
  processing of active array data in inverse scattering via truncated reduced
  order models}, Journal of Computational Physics, 381 (2019), pp.~1--26.

\bibitem{borcea2021reduced}
{\sc L.~Borcea, J.~Garnier, A.~Mamonov, and J.~Zimmerling}, {\em Reduced order
  model approach for imaging with waves}, Inverse Problems, 38 (2021),
  p.~025004.

\bibitem{borcea2022waveform}
\leavevmode\vrule height 2pt depth -1.6pt width 23pt, {\em Waveform inversion
  via reduced order modeling}, Geophysics, 88 (2022), pp.~1--91.

\bibitem{borcea2022internal}
\leavevmode\vrule height 2pt depth -1.6pt width 23pt, {\em Waveform inversion
  with a data driven estimate of the internal wave}, SIAM Journal on Imaging
  Sciences, 16 (2023), pp.~280--312.

\bibitem{PublishedSirev}
{\sc L.~Borcea, J.~Garnier, A.~V. Mamonov, and J.~Zimmerling}, {\em When data
  driven reduced order modeling meets full waveform inversion}, SIAM Review, 66
  (2024), pp.~501--532.

\bibitem{borges2017high}
{\sc C.~Borges, A.~Gillman, and L.~Greengard}, {\em High resolution inverse
  scattering in two dimensions using recursive linearization}, SIAM Journal on
  Imaging Sciences, 10 (2017), pp.~641--664.

\bibitem{brunton2019data}
{\sc S.~Brunton and J.~Kutz}, {\em {Data-driven science and engineering:
  Machine learning, dynamical systems, and control}}, Cambridge University
  Press, 2019.

\bibitem{brunton2016discovering}
{\sc S.~Brunton, J.~Proctor, and J.~Kutz}, {\em Discovering governing equations
  from data by sparse identification of nonlinear dynamical systems},
  Proceedings of the National Academy of Sciences, 113 (2016), pp.~3932--3937.

\bibitem{bunks1995multiscale}
{\sc C.~Bunks, F.~Saleck, S.~Zaleski, and G.~Chavent}, {\em Multiscale seismic
  waveform inversion}, Geophysics, 60 (1995), pp.~1457--1473.

\bibitem{chen1997inverse}
{\sc Y.~Chen}, {\em {Inverse scattering via Heisenberg's uncertainty
  principle}}, Inverse problems, 13 (1997), p.~253.

\bibitem{cheney2009fundamentals}
{\sc M.~Cheney and B.~Borden}, {\em Fundamentals of radar imaging}, SIAM, 2009.

\bibitem{colin2009}
{\sc Y.~Colin~de Verdi\`ere}, {\em Semiclassical analysis and passive imaging},
  Nonlinearity, 22 (2009), pp.~R45--R75.

\bibitem{corte2020deep}
{\sc G.~C{\^o}rte, J.~Dramsch, H.~Amini, and C.~MacBeth}, {\em {Deep neural
  network application for 4D seismic inversion to changes in pressure and
  saturation: Optimizing the use of synthetic training datasets}}, Geophysical
  Prospecting, 68 (2020), pp.~2164--2185.

\bibitem{curlander1991synthetic}
{\sc J.~Curlander and R.~McDonough}, {\em Synthetic aperture radar}, vol.~11,
  Wiley, 1991.

\bibitem{ding2022coupling}
{\sc W.~Ding, K.~Ren, and L.~Zhang}, {\em Coupling deep learning with full
  waveform inversion}, arXiv preprint arXiv:2203.01799,  (2022).

\bibitem{druskin2016direct}
{\sc V.~Druskin, A.~V. Mamonov, A.~E. Thaler, and M.~Zaslavsky}, {\em Direct,
  nonlinear inversion algorithm for hyperbolic problems via projection-based
  model reduction}, SIAM Journal on Imaging Sciences, 9 (2016), pp.~684--747.

\bibitem{druskin2018nonlinear}
{\sc V.~Druskin, A.~V. Mamonov, and M.~Zaslavsky}, {\em A nonlinear method for
  imaging with acoustic waves via reduced order model backprojection}, SIAM
  Journal on Imaging Sciences, 11 (2018), pp.~164--196.

\bibitem{druskin2002three}
{\sc V.~Druskin and S.~Moskow}, {\em {Three-point finite-difference schemes,
  Pad{\'e} and the spectral Galerkin method. I. One-sided impedance
  approximation}}, Mathematics of computation, 71 (2002), pp.~995--1019.

\bibitem{EngquistFroese}
{\sc B.~Engquist and B.~Froese}, {\em {Application of the Wasserstein metric to
  seismic signals}}, {Communications in Mathematical Sciences}, 12 (2014),
  pp.~979--988.

\bibitem{engquist2022optimal}
{\sc B.~Engquist and Y.~Yang}, {\em Optimal transport based seismic inversion:
  Beyond cycle skipping}, Communications on Pure and Applied Mathematics, 75
  (2022), pp.~2201--2244.

\bibitem{fink99}
{\sc M.~Fink}, {\em Time-reversed acoustics}, Scientific American, 281 (1999),
  pp.~91--97.

\bibitem{bookfouque}
{\sc J.-P. Fouque, J.~Garnier, G.~Papanicolaou, and K.~S{\o}lna}, {\em Wave
  Propagation and Time Reversal in Randomly Layered Media}, Springer, New York,
  2007.

\bibitem{passivecomm}
{\sc J.~Garnier}, {\em Passive communication with ambient noise}, SIAM J. Appl.
  Math., 81 (2021), pp.~814--833.

\bibitem{garnierpassive2009}
{\sc J.~Garnier and G.~Papanicolaou}, {\em Passive sensor imaging using cross
  correlations of noisy signals in a scattering medium}, SIAM J. Imaging
  Sciences, 2 (2009), pp.~396--437.

\bibitem{garnierpassive2010}
\leavevmode\vrule height 2pt depth -1.6pt width 23pt, {\em Resolution analysis
  for imaging with noise}, Inverse Problems, 26 (2010), p.~074001.

\bibitem{garnier2016passive}
\leavevmode\vrule height 2pt depth -1.6pt width 23pt, {\em Passive imaging with
  ambient noise}, Cambridge University Press, 2016.

\bibitem{gauthier1986two}
{\sc O.~Gauthier, J.~Virieux, and A.~Tarantola}, {\em Two-dimensional nonlinear
  inversion of seismic waveforms: Numerical results}, Geophysics, 51 (1986),
  pp.~1387--1403.

\bibitem{gill2019practical}
{\sc P.~Gill, W.~Murray, and M.~Wright}, {\em Practical optimization}, SIAM,
  2019.

\bibitem{gilman2017transionospheric}
{\sc M.~Gilman, E.~Smith, and S.~Tsynkov}, {\em Transionospheric synthetic
  aperture imaging}, Birkh\"auser, Cham, 2017.

\bibitem{golubVanLoan}
{\sc G.~H. Golub and C.~F. Van~Loan}, {\em Matrix computations}, John Hopkins
  University Press, Baltimore, 2013.

\bibitem{gradshteyn2014table}
{\sc I.~Gradshteyn and I.~Ryzhik}, {\em Table of integrals, series, and
  products}, Academic press, 2014.

\bibitem{herkt2013convergence}
{\sc S.~Herkt, M.~Hinze, and R.~Pinnau}, {\em {Convergence analysis of Galerkin
  POD for linear second order evolution equations}}, Electron. Trans. Numer.
  Anal, 40 (2013), pp.~321--337.

\bibitem{hesthaven2022reduced}
{\sc J.~Hesthaven, C.~Pagliantini, and G.~Rozza}, {\em Reduced basis methods
  for time-dependent problems}, Acta Numerica, 31 (2022), pp.~265--345.

\bibitem{huang2018source}
{\sc G.~Huang, R.~Nammour, and W.~Symes}, {\em Source-independent extended
  waveform inversion based on space-time source extension: Frequency-domain
  implementation}, Geophysics, 83 (2018), pp.~R449--R461.

\bibitem{yu2001global}
{\sc O.~Imanuvilov and M.~Yamamoto}, {\em Global uniqueness and stability in
  determining coefficients of wave equations}, Communications in Partial
  Differential Equations, 26 (2001), pp.~1409--1425.

\bibitem{jollife2016principal}
{\sc I.~Jollife and J.~Cadima}, {\em Principal component analysis: A review and
  recent developments}, Philosophical Transactions of the Royal Society A:
  Mathematical, Physical and Engineering Sciences, 374 (2016), p.~20150202.

\bibitem{karhunen1947}
{\sc K.~Karhunen}, {\em Uber lineare methoden in der
  wahrscheinlichkeitsrechnung}, Ann. Acad. Sci. Fennicae. Ser. A. I.
  Math.-Phys., 37 (1947), pp.~1--79.

\bibitem{kato}
{\sc T.~Kato}, {\em Perturbation theory for linear operators}, vol.~132,
  Springer Science \& Business Media, 2013.

\bibitem{kunisch2010optimal}
{\sc K.~Kunisch and S.~Volkwein}, {\em {Optimal snapshot location for computing
  POD basis functions}}, ESAIM: Mathematical Modelling and Numerical Analysis,
  44 (2010), pp.~509--529.

\bibitem{lieu2006reduced}
{\sc T.~Lieu, C.~Farhat, and M.~Lesoinne}, {\em Reduced-order fluid/structure
  modeling of a complete aircraft configuration}, Computer methods in applied
  mechanics and engineering, 195 (2006), pp.~5730--5742.

\bibitem{lobkis2001}
{\sc O.~I. Lobkis and R.~L. Weaver}, {\em On the emergence of the {G}reen's
  function in the correlations of a diffuse field}, J. Acoustic. Soc. Am., 110
  (2001).

\bibitem{loeve1978}
{\sc M.~Loeve}, {\em Probability theory. Vol. II, 4th ed.}, Springer-Verlag,
  1978.

\bibitem{mahankali2020convexity}
{\sc S.~Mahankali and Y.~Yang}, {\em The convexity of optimal transport-based
  waveform inversion for certain structured velocity models}, arXiv preprint
  arXiv:2009.00708,  (2020).

\bibitem{nachman2007conductivity}
{\sc A.~Nachman, A.~Tamasan, and A.~Timonov}, {\em Conductivity imaging with a
  single measurement of boundary and interior data}, Inverse Problems, 23
  (2007), p.~2551.

\bibitem{schilders2008model}
{\sc W.~Schilders, H.~Van~der Vorst, and J.~Rommes}, {\em Model order
  reduction: theory, research aspects and applications}, vol.~13, Springer,
  2008.

\bibitem{stefanov2005stable}
{\sc P.~Stefanov and G.~Uhlmann}, {\em Stable determination of generic simple
  metrics from the hyperbolic dirichlet-to-neumann map}, International
  Mathematics Research Notices, 2005 (2005), pp.~1047--1061.

\bibitem{symes2008migration}
{\sc W.~Symes}, {\em Migration velocity analysis and waveform inversion},
  Geophysical prospecting, 56 (2008), pp.~765--790.

\bibitem{symes2022error}
\leavevmode\vrule height 2pt depth -1.6pt width 23pt, {\em Error bounds for
  extended source inversion applied to an acoustic transmission inverse
  problem}, Inverse Problems, 38 (2022), p.~115002.

\bibitem{herrmann2013}
{\sc T.~van Leeuwen and F.~Herrmann}, {\em Mitigating local minima in
  full-waveform inversion by expanding the search space}, Geophysical Journal
  International, 195 (2013), pp.~661--667.

\bibitem{virieux2009overview}
{\sc J.~Virieux and S.~Operto}, {\em An overview of full-waveform inversion in
  exploration geophysics}, Geophysics, 74 (2009), pp.~WCC1--WCC26.

\bibitem{warner2016}
{\sc M.~Warner and L.~Guasch}, {\em Adaptive waveform inversion: theory},
  Geophysics, 81 (2016), pp.~R429--R445.

\bibitem{yang2018application}
{\sc Y.~Yang, B.~Engquist, J.~Sun, and B.~Hamfeldt}, {\em {Application of
  optimal transport and the quadratic Wasserstein metric to full-waveform
  inversion}}, Geophysics, 83 (2018), pp.~R43--R62.

\end{thebibliography}

\end{document}